\documentclass[12pt]{amsart}
\usepackage{amsmath,amssymb,amscd,array, mathrsfs }

\usepackage{amsmath,amscd,amssymb,amsthm,array}
\usepackage{color}

\usepackage{amsmath,amssymb,mathrsfs,amsthm, tikz-cd,mathrsfs}
\usepackage{comment}
\def\url#1{\expandafter\s

\tring\csname #1\endcsname}

\def\mmat #1,#2,#3,#4,{\text{\small\arraycolsep=3pt $
\begin{pmatrix}#1&#2\\#3&#4\end{pmatrix}$}}

\usepackage{hyperref}
\usepackage{capt-of}

\usepackage{multirow}
\usepackage[all]{xy}

\usepackage{comments}
\newComments\SBe{Said}{blue}
\newComments\SBo{Sofiane}{blue}
\newComments\AM{Nacer}{blue}
\newComments\DL{DL}{red}
\newComments\QEh{QEh}{blue}

\def\mmat #1,#2,#3,#4,{\text{\small\arraycolsep=3pt $
\begin{pmatrix}#1&#2\\#3&#4\end{pmatrix}$}}

\usepackage{lscape}
\usepackage{tikz-cd}
\usepackage{enumerate}

\usepackage{DLdef1}
\usepackage{multicol}

\def\mmat #1,#2,#3,#4,{\text{\small\arraycolsep=3pt $
\begin{pmatrix}#1&#2\\#3&#4\end{pmatrix}$}}

\renewcommand {\ssbegin}[2][*]
 {\refstepcounter{subsection}%
\if#1*
\addcontentsline{toc}{subsection}{\thesubsection.\hskip 1pc #2}%
\else
\addcontentsline{toc}{subsection}{\thesubsection.\hskip 1pc #2. #1}%
\fi
 \def \secno {\gdef \secno {}{\ssecfont
\thesubsection.\hskip 2ex}%
 }%
 \begin{#2}}

\renewcommand {\sssbegin}[2][*]
 {\refstepcounter{subsubsection}
\if#1*
\addcontentsline{toc}{subsubsection}{\thesubsubsection.\hskip 1pc #2}%
\else
\addcontentsline{toc}{subsubsection}{\thesubsubsection.\hskip 1pc #2. #1}
\fi
 \def \secno {\gdef \secno {}{\ssecfont \thesubsubsection.\hskip 2ex}%
 }%
 \begin{#2}}

\renewcommand {\parbegin}[2][*]
 {\refstepcounter{paragraph}
\if#1*
\addcontentsline{toc}{paragraph}{\theparagraph.\hskip 1pc #2}%
\else
\addcontentsline{toc}{paragraph}{\theparagraph.\hskip 1pc #2. #1}
\fi
 \def \secno {\gdef \secno {}{\ssecfont \theparagraph.\hskip 2ex}%
 }%
 \begin{#2}}

\rmnameii{etr}{etr}
\rmnameii{evv}{ev}

\DeclareMathOperator{\K}{\mathbb{K}}
\newcommand{\Z}{\mathbb{Z}}

\newcommand {\A}{{\cal{A}}}

\newcommand{\Ll}{{\mathrm{L}}}
\newcommand{\Rr}{{\mathrm{R}}}
\newcommand{\too}{\longrightarrow}
    \newcommand{\Om}{\Omega}

\newcommand{\al}{\alpha}
\newcommand{\ga}{\gamma}

\newcommand{\la}{\lambda}

\newcommand{\de}{\delta}  
\newcommand{\prs}{\langle\;,\;\rangle}
\def\br{[\;,\;]}
\newcommand{\esp}{\quad\mbox{and}\quad}
\begin{document}

\title[Flat pseudo-Euclidean Leibniz superalgebras ]{Flat pseudo-Euclidean Leibniz superalgebras }

\author{Sa\"id Benayadi}
\address {Universit\'e de Lorraine, Laboratoire LMAM, CNRS-UMR 7122,\\ Ile du Saulcy, F-57045 Metz
		cedex 1, France.}
\email{said.benayadi@univ-lorraine.fr}

\author{Sofiane Bouarroudj}

\address{Division of Science and Mathematics, New York University Abu Dhabi, P.O. Box 129188, Abu Dhabi, United Arab Emirates.}
\email{sofiane.bouarroudj@nyu.edu}

\author{Hamza El Ouali}
\address{Universit\'e Cadi-Ayyad,
	Facult\'e des sciences et techniques,
	B.P. 549, Marrakech, Maroc.}
\email{eloualihamza11@gmail.com}

\thanks{S.B. was supported by the grant NYUAD-065}

\keywords{Leibniz superalgebra, pre-left Leibniz superalgebra, pseudo-Euclidean left-Leibniz superalgebra, flat pseudo-Euclidean left-Leibniz superalgebra, Levi-Civita product, $T^*$-extension, double extension}

 \subjclass[2020]{17A32; 17A70; 17A60; 17D25}

\begin{abstract}
In this paper, we introduce pre-Lie and pre-Leibniz superalgebras, which generalize pre-Lie and pre-Leibniz algebras  to the super setting. Additionally, we define a Levi-Civita product associated with a symmetric non-degenerate bilinear form on a non-associative superalgebra. This leads to the definition of flat pseudo-Euclidean left Leibniz superalgebras as those whose Levi-Civita product induces a pre-Leibniz structure. We study the structure of flat pseudo-Euclidean left Leibniz superalgebras and provide a characterization theorem.
In the second part, we focus on quadratic Leibniz superalgebras and show that such a superalgebra is flat if and only if it is symmetric Leibniz and 2-step nilpotent. We further study the structure of quadratic 2-step nilpotent symmetric Leibniz superalgebras. Finally, we introduce the notion of double extension for flat pseudo-Euclidean (resp. Lie) left Leibniz superalgebras and prove that any flat pseudo-Euclidean non-Lie left Leibniz superalgebra can be obtained by a sequence of double extensions starting from a flat pseudo-Euclidean Lie superalgebra.

\end{abstract}


\maketitle

\thispagestyle{empty}
\setcounter{tocdepth}{2}
\tableofcontents

\section{Introduction}
\subsection{On a Milnor's problem} A pseudo-Riemannian Lie group is a Lie group $G$ equipped with a left-invariant pseudo-Riemannian metric $\mu$. The Lie algebra $\mathfrak{g} = T_e G$, endowed with the scalar product $\prs = \mu_e$, is called a pseudo-Riemannian Lie algebra. The Levi-Civita connection $\nabla$ associated with the pseudo-Riemannian metric $\mu$ defines a bilinear product $(u, v) \mapsto u \star v$ on $\mathfrak{g}$, called the Levi-Civita product, given by the  Koszul’s formula:
$$
2 \langle u \star v, w \rangle = \langle [u, v], w \rangle + \langle [w, u], v \rangle + \langle [w, v], u \rangle
\quad \text{for all } u,v,w \in \mathfrak{g}.
$$
This product is uniquely determined by the metric and satisfies
$$u\star v-v\star u=[u,v],\esp \langle u\star v, w\rangle=-\langle v, u\star w\rangle, \; \text{for all }u,v,w\in \mathfrak{g}.$$
The curvature of $\mathfrak{g}$ at the identity element $e$ is given by (where  $  v\mapsto
 \Ll_u(v):= u \star v
$ is the left multiplication):
$$
K(u, v) := \Ll_{[u, v]} - [\Ll_u, \Ll_v],
$$
for all $u,v\in \mathfrak{g}$. If the curvature $K$ vanishes, then $(G, \mu)$ is called a flat pseudo-Riemannian Lie group, and its Lie algebra $(\mathfrak{g}, \prs)$ is called a flat pseudo-Riemannian Lie algebra.
 The vanishing of the curvature is equivalent to the fact that $\mathfrak{g}$, endowed with the Levi-Civita product, is a left-symmetric algebra, see Definition \ref{left-sym}. Thus, a flat pseudo-Riemannian Lie algebra can be viewed as a left-symmetric algebra equipped with a symmetric non-degenerate bilinear form such that the left multiplications are skew-symmetric with respect to this bilinear form.

In 1976, J. Milnor \cite{Milnor1} characterized flat Euclidean Lie algebras by proving that a Lie algebra $\mathfrak{g}$ is a flat Euclidean Lie algebra if and only if it splits orthogonally as $\mathfrak{g} = \mathfrak{b} \oplus \mathfrak{u}$, where $\mathfrak{u}$ is an abelian ideal, $\mathfrak{b}$ is an abelian subalgebra, and $\operatorname{ad}_b$ is skew-symmetric with respect to the inner product for every $b \in \mathfrak{b}$, where $ \operatorname{ad}$ is the adjoint operator on $\mathfrak{g}$. A. Aubert and A. Medina initiated the study of  flat {\it Lorentzian} nilpotent Lie algebras in \cite{Aub}. Later, they introduced a construction method known as the double extension of flat pseudo-Euclidean Lie algebras, which is an analogue of the quadratic double extension method developed earlier by A. Medina and P. Revoy \cite{Medina2}.  They showed that flat Lorentzian nilpotent Lie algebras can be obtained by this method, and that they are at most 3-step nilpotent. Furthermore, M. Guediri  \cite{guediri} classified 2-step nilpotent flat Lorentzian Lie groups and proved that their Lie algebras are a trivial extension of the 3-dimensional Heisenberg algebra $H_3$, and that the restriction of the metric to $H_3$ is Lorentzian with a degenerate center. More recently,  M. Ait Ben Haddou, M. Boucetta and H. Lebzioui \cite{ABL} studied flat Lorentzian Lie algebras. They showed that flat Lorentzian nilpotent Lie algebras with degenerate center are obtained by the double extension method. Additionally,  they  proved in \cite{BL1} that non-unimodular flat Lorentzian nilpotent Lie algebras can be constructed using this method. Later, M. Boucetta and H. Lebzioui \cite{BL2} studied flat pseudo-Euclidean Lie algebras of signature $(2, n -2)$. They showed that every nilpotent flat pseudo-Euclidean Lie algebra with such a signature can be obtained as a double extension of a nilpotent flat Lorentzian Lie algebra.

\subsection{Leibniz (super)algebras} 
Leibniz algebras were originally introduced by A. Bloh under the name of ``D-algebras" in the mid-1960s (see \cite{Bloh1, Bloh2, Bloh3}). In the early 1990s, these algebras were rediscovered and developed by Jean-Louis Loday while studying  cohomology groups of Lie algebras (see \cite{Loday1, Loday2}). Leibniz algebras can be viewed as analogues of Lie algebras without the anticommutativity condition. This gives rise to two types of Leibniz algebras: left Leibniz algebras and right Leibniz algebras. The name ``left" (resp. ``right") Leibniz algebra comes from the fact that the left (resp. right) multiplication is a derivation. It is worth noting that by adding the anticommutativity condition to a (left or right) Leibniz algebra, a Leibniz algebra turns into a  Lie algebra. Therefore, the class of Leibniz algebras is a generalization of the class of Lie algebras. A non-associative algebra is called a symmetric Leibniz algebra if it is both a left and a right Leibniz algebra \cite{Mason}. On the other hand, in \cite{Tang}, the authors introduced the notion of a pre-Leibniz algebra, originally referred to as a Leibniz-dendriform algebra. A pre-Leibniz algebra naturally induces a Leibniz algebra structure on its underlying vector space. 

Leibniz algebras appear nowadays in various fields such as differential geometry, homological algebra, K-theory, classical algebraic topology, noncommutative geometry, and physics. In particular, real left Leibniz algebras are the infinitesimal versions of Lie racks. In 2004, Kinyon \cite{Kinyon} proved that if $(X, e)$ is a pointed Lie rack, then $T_e X$ carries a structure of a left Leibniz algebra.

A. S. Dzhumadil’daev introduced the notion of Leibniz superalgebra in \cite{Dz}, which is both a non-anti-supercommutative analogue of Lie superalgebras and an extension of Leibniz superalgebras. Recently, left and right Leibniz superalgebras have been extensively studied (see, for example, \cite{Omirov1, Omirov2, BCON, CaSa, CT, kh, Liu1, Liu2, Omirov3}).

There is also the notion of Zinbiel superalgebras; see, for instance, \cite{BM2} and references therein, but their study is outside the scope of this paper.


\subsection{Double and $T^*$-extensions} Quadratic non-associative algebras are non-associative algebras equipped with a pseudo-Euclidean bilinear form $\prs$ that is invariant (also called associative). It is well known that the existence of quadratic non-associative (super)algebras is a powerful tool for studying their structure (for examples the Killing form on semisimple Lie algebras, the super-Killing form on certain classes of simple Lie superalgebras, and the Albert form on semisimple Jordan algebras). Quadratic Lie (super)algebras have been studied in \cite{Benayadi6, Benayadi2, Benayadi5} and references therin. In \cite{BenFa0, BenFa}, the authors studied even quadratic Leibniz superalgebras, and in \cite{BenFa1}, they considered the odd quadratic case. On the other hand, in \cite{Medina2}, it was shown that a quadratic Lie algebra is flat if and only if it is 2-step nilpotent.
 In this paper, we extend this result by proving that a quadratic left Leibniz superalgebra is flat if and only if it is a 2-step nilpotent symmetric Leibniz superalgebra. 
 
 In \cite{Bordemann}, M. Bordemann introduced the notion of $T^*$-extension for quadratic non-associative algebras. He proved that a non-associative algebra $\A$, endowed with a quadratic structure, is either a $T^*$-extension (or a non-degenerate ideal of codimension one in a $T^*$-extension) of another non-associative algebra $\mathfrak{h}$ if and only if $\A$ contains a totally isotropic ideal of dimension $\left\lfloor \frac{\dim \A}{2} \right\rfloor$, where $\left\lfloor \cdot \right\rfloor$ denotes the integer part. In the context of even quadratic Leibniz superalgebras, this notion was introduced in \cite{BenFa0}. In that paper, the authors provided a complete description of such quadratic solvable Leibniz superalgebras via the $T^*$-extension construction.

 \subsection{On the results}
In this paper, we introduce the concept of pre-left (resp. pre-right) Leibniz superalgebras $(\A, \star, \circ)$, which generalizes the concept of pre-Leibniz algebras (see Definition \ref{pre-Leibniz}). We show that a pre-left (resp. pre-right) Leibniz superalgebra $(\mathcal{A}, \star, \circ)$ naturally induces a left (resp. right) Leibniz superalgebra structure on the underlying $\mathbb{Z}_2$-graded vector superspace. In particular, if $\star = -\circ$, then the pre-Leibniz superalgebra reduces to a left-symmetric superalgebra,  showing that the class of left-symmetric superalgebras is a subclass of pre-Leibniz superalgebras.

We introduce an analogue of the Levi-Civita product in the context of pseudo-Euclidean non-associative superalgebras (see Proposition \ref{LVCT}). Recall that a  superalgebra $(\mathcal{A}, \bullet)$ is called pseudo-Euclidean if it is equipped with a symmetric non-degenerate bilinear form $\prs$. Since left Leibniz superalgebras generalize both Lie superalgebras and Lie algebras, we introduce a notion of curvature that extends the classical curvature concept defined in the Lie algebra setting (see Definition~\ref{courbures}). 
If the curvature of a pseudo-Euclidean left Leibniz superalgebra $(\A, \bullet, \prs)$ vanishes, then $(\A, \bullet, \prs)$ is called a flat pseudo-Euclidean left Leibniz superalgebra. This vanishing is equivalent to the fact that $(\A, \bullet, \prs)$, endowed with the Levi-Civita products, is a pre-left Leibniz superalgebra. Moreover, we study flat pseudo-Euclidean left Leibniz superalgebra and we show that every flat pseudo-Euclidean left Leibniz superalgebra is a symmetric Leibniz superalgebra. We provide a characterization of flat pseudo-Euclidean left Leibniz superalgebra (see Theorem \ref{chara}).

Next, in this paper, we study flat even (resp. odd) quadratic Leibniz superalgebras. We show that an  even (resp. odd) quadratic Leibniz superalgebra is flat if and only if it is a 2-step nilpotent symmetric Leibniz superalgebra. We also provide a characterization of flat even (resp. odd) quadratic Leibniz superalgebras (see Theorem \ref{chara-2-step}). Furthermore, the notion of $T^*$-extension of even quadratic Leibniz superalgebras, introduced by S. Benayadi and M. Fahmi in \cite{BenFa0}, is extended by introducing the concept of odd quadratic Leibniz superalgebras. We prove that every even (resp. odd) reduced 2-step nilpotent symmetric quadratic Leibniz superalgebra can be obtained as a $T^*$-extension (resp. $\Pi(T^*)$-extension) of a trivial superalgebra by means  a non-degenerate bilinear form. We also provide an explicit description of quadratic 2-step nilpotent symmetric Leibniz superalgebras.

Finally, we introduce the notion of double extension of flat even (resp. odd) pseudo-Euclidean left Leibniz superalgebras, and we describe the particular case of Lie superalgebras within this framework. We prove that every flat even (resp. odd) pseudo-Euclidean non-Lie left Leibniz superalgebra can be obtained by a sequence of double extensions starting from a flat even (resp. odd) pseudo-Euclidean Lie superalgebra.

$$
$$

\section{Backgrounds}
Let $\mathbb{K}$ be an arbitrary field of characteristic zero. The group of integers modulo $2$ is denoted by $\Z_{2}$. 

 Let $V=V_{\bar 0}\oplus V_{\bar 1}$ be a $\Z_{2}$-graded space defined over an arbitray field ${\mathbb K}$. The degree of a homogeneous element $v\in V_{\bar{i}}$ is denoted by $|v|:=\bar{i}$. The element $v$ is called \textit{even} if $v\in V_{\bar 0}$ and \textit{odd} if $v\in V_{\bar 1}$. Throughout the text, all elements are supposed to be homogeneous unless otherwise stated. A linear map $\varphi:V\rightarrow W$ between $\Z_{2}$-graded space is called \textit{even} if $\varphi(V_{\bar{i}})\subset W_{\bar{i}}$ and \textit{odd} if $\varphi(V_{\bar{i}})\subset W_{\overline{i}+\overline{1}}$.

Let $V=V_{\bar 0}\oplus V_{\bar 1}$ be a $\Z_2$-graded vector space. We denote by $\Pi$ the \textit{change of parity functor} $\Pi: V\mapsto \Pi (V)$, where $\Pi(V)$ is another copy of $V$ such that $ \Pi(V)_{\bar 0}:=V_{\bar 1};~~\Pi(V)_{\bar 1}:=V_{\bar 0}$.  Elements of $\Pi(V)$ shall be denoted by $\Pi(v),~\forall v\in V$. In fact, $|\Pi(v)|=|v|+{\bar 1}$ for every homogeneous $v\in V$.

Let $U$ and $V$ be two $\mathbb{Z}_2$-graded vector spaces. Given a morphism $G: U \to V$, we define the morphism
$
G^\Pi : \Pi(U) \to \Pi(V)
$
defined by $G^\Pi(\Pi(u))=\Pi(G(u))$. It is straightforward to verify that this construction makes the assignment
$
U \mapsto \Pi(U)
$ into a functor satisfying $\Pi^2 = \mathrm{id}$.

Moreover, given a morphism $G: U \to V$, one can construct morphisms $
\Pi(G) : U \to \Pi(V), $ defined by $ \Pi(G)(u) = \Pi(G(u)),
$ and
$
G\Pi : \Pi(U) \to V,$ defined by  $(G\Pi)(\Pi(u)) = G(u).
$
Clearly, we have the factorization $
G^\Pi = \Pi(G\Pi)$, see
\cite{Youri}.

\subsection{Bilinear forms on a superspace} Let $V$ and $W$ be two superspaces. Let ${\cal B}\in \text{Bil}(V,W)$ be a homogeneous bilinear form. Recall that the Gram matrix $B=(B_{ij})$ associated to $\cal B$ is given by the following formula:
\begin{equation*}\label{martBil}
B_{ij}=(-1)^{p({\cal B})p(v_i)}{\cal B}(v_{i}, w_{j})\text{~~for the basis vectors $v_{i}\in V$ and $w_{j}\in W$.}
\end{equation*}
This definition allows us to identify a~bilinear form $B(V, W)$ with an element of $\mathrm{Hom}(V, W^*)$. 
Consider the \textit{upsetting} of bilinear forms
$u:\text{Bil} (V, W)\rightarrow \text{Bil}(W, V)$ given by the formula 
\begin{equation*}\label{susyB}
u({\cal B})(w, v)=(-1)^{p(v)p(w)}{\cal B}(v,w)\text{~~for any $v \in V$ and $w\in W$.}
\end{equation*}
In terms of the Gram matrix $B$ of ${\cal B}$, we have 
\begin{equation*}
u(B)=
\left ( \begin{array}{cc}  R^{t} & (-1)^{p(B)}T^{t} \\ (-1)^{p(B)}S^{t} & -U^{t} \end{array} \right ),
\text{ for $B=\left (  \begin{array}{cc}R & S \\ T & U \end{array} \right )$}.
\end{equation*}
Suppose now that $V=W$. Following \cite{L}, we say that: 
\begin{itemize} \item[(i)] the form
$\cal B$ is  \textit{symmetric} if  and only if $u(B)=B$;

\item[(ii)] the form $\cal B$ is \textit{anti-symmetric} if and only if $u(B)=-B$.

\item[(iii)] the form $\cal B$ is non-degenerate if $u\in V$ satisfies $\cal B(v,u)=0,$ for any $v\in V$, then $u = 0$,
\item[(iv)] the form $\cal B$ is even if $\cal B(V_{\bar{0}}, V_{\bar{1}})=\cal B(V_{\bar{1}},V_{\bar{0}})=\{0\}$,
 \item[(v)] the form $\cal B$ is odd if  ${\cal B}(V_{\al}, V_{\al}) = \{0\}$, for all $\al\in\mathbb{Z}_2,$ 

\end{itemize}

From now on, we will denote a bilinear form with $\prs$ instead. 
\sssbegin{Definition}
Let $(V,\prs)$ be a pseudo-Euclidean $\mathbb{Z}_2$-graded vector space, where $\prs$ is a non-degenerate symmetric bilinear form on $V$. Let $f$ be a homogeneous endomorphism of $V$ of degree $\al \in \{\bar{0},\bar{1}\}$.
\begin{itemize}
  \item[(i)] The adjoint of $f$ with respect to $\prs$, denoted by $f^*$, is the unique homogeneous endomorphism of degree $\al$ satisfying
  \[
  \langle f(u), v\rangle = (-1)^{\al|u|} \langle u, f^*(v)\rangle, \quad \text{for all } u,v\in V.
  \]
\item[(ii)] The endomorphism $f$ is said to be:
  \begin{itemize}
    \item Symmetric (or self-adjoint) with respect to $\prs$ if $f^* = f$, i.e.,
    \[
    \langle f(u), v\rangle = (-1)^{\al|u|} \langle u, f(v)\rangle, \quad \text{for all } u,v\in V.
    \]

    \item Anti-symmetric  (or anti-self-adjoint) with respect to $\prs$ if $f^* = -f$, i.e.,
    \[
     \langle f(u), v\rangle = -(-1)^{\al|u|} \langle u, f(v)\rangle, \quad \text{for all } u,v\in V.
    \]
  \end{itemize}
\end{itemize}
\end{Definition}
\subsection{Leibniz superalgebras}
A superalgebra $(\A, \bullet)$  over a field $\mathbb{K}$ consists of a $\Z_{2}$-graded space $\A=\A_{\bar 0}\oplus \A_{\bar 1}$ and a binary operation satisfying: $\A_\alpha \bullet \A_\beta \subset \A_{\alpha+\beta}, \forall \alpha, \beta \in \mathbb{Z}_2$. 

Recall that there are two versions of skew- or anti-commutativity that are synonyms only in the non-super case; for two homogeneous elements $a$ and $b$ of a superalgebra we call the following conditions
\[
\begin{array}{lcll}
b \bullet a & = & (-1)^{|b||a|}a \bullet b & \text{super commutativity,} \\[2mm]
b \bullet a & = & -(-1)^{|b||a|}a \bullet b & \text{super anti-commutativity, }\\[2mm]
b \bullet a & = & (-1)^{(|b|+1)(|a|+1)} a \bullet b & \text{super skew-commutativity, }\\[2mm]
b \bullet a & = &  -(-1)^{(|b|+1)(|a|+1)}a \bullet b & \text{super antiskew-commutativity.}
\end{array}
\]
 For an element $u \in \A$, we denote by $\Ll_u$ and $\Rr_u$, respectively, the left and the right multiplication operators in $\A$. More precisely $\Ll_u$ and $\Rr_u$ are defined as  follows: $\Ll^\bullet_u(v):=u \bullet v$ and $\Rr^\bullet_u(v):=(-1)^{|u||v|}v\bullet u,$ for all $v \in \A$. 

\sssbegin{Definition}
     A  pseudo-Euclidean superalgebra $( \A, \bullet, \prs)$ is a superalgebra $(\A, \bullet)$ endowed with a
non-degenerate, symmetric bilinear form $\prs$.
\end{Definition}

\sssbegin{Definition}
    Let $(\A, \bullet)$ be a  superalgebra. A bilinear form $\prs$ on $\A$ is invariant if \[
    \langle u\bullet v, w\rangle = \langle u, v\bullet w\rangle,
    \]
    for all $u,v,w\in \A$.
\end{Definition}

\sssbegin{Definition}
  A disuperalgebra over $\K$ is a $\mathbb{Z}_2$-graded vector space $\A$ endowed with
two binary operations $\star, \circ:\, \A\times\A \too\A$, called left and right products, satisfying: 
$$\A_\alpha \star \A_\beta\subset \A_{\alpha+\beta}, \quad  \A_\alpha \circ \A_\beta \subset \A_{\alpha+\beta},\;\; \forall \alpha,\beta\in \mathbb{Z}_2.$$ 
\end{Definition}

\sssbegin{Definition}
    Let $(\A, \bullet)$ be a superalgebra.
    \begin{enumerate}
  \item[(i)]  $(\mathcal{A}, \bullet)$ is called a left Leibniz superalgebra if, for any $u, v, w \in \mathcal{A}$,
$$
u \bullet(v \bullet w)=(u \bullet v) \bullet w+(-1)^{|u||v|}v \bullet(u \bullet w).
$$
\item[(ii)] $(\mathcal{A}, \bullet)$ is called a right Leibniz superalgebra if, for any $u, v, w \in \mathcal{A}$,
$$
u\bullet(v \bullet w) =(u \bullet v) \bullet w-(-1)^{|v||w|} (u \bullet w) \bullet v.
$$
\item[(iii)] If $(\mathcal{A}, \bullet)$ is both left and right Leibniz then it is called a symmetric Leibniz superalgebra.  
 \end{enumerate}
     \end{Definition}

Let $(\mathcal{A}, \bullet)$ be a left (resp. right) Leibniz superalgebra. Then the left (resp. right) multiplication operators satisfy the identity:
\begin{equation}
\mathrm{L}^\bullet_{u \bullet v} = \left[\mathrm{L}^\bullet_u, \mathrm{L}^\bullet_v\right]
\quad (\text{resp. } \mathrm{R}^\bullet_{u \bullet v} = -\left[\mathrm{R}^\bullet_u, \mathrm{R}^\bullet_v\right]),
\label{mul}\end{equation}
for all $u, v \in \A$.

It is obvious that $(\mathcal{A}, \bullet)$ is a left Leibniz superalgebra if and only if $\left(\mathcal{A}, \bullet_{\text {opp }}\right)$ is a right Leibniz superalgebra where  $u \bullet_{\text {opp }} v:=(-1)^{|u||v|}v \bullet u$. It is also obvious that any Lie superalgebra is symmetric Leibniz superalgebra. However, the class of Leibniz superalgebras is much larger than the class of Lie superalgebras, and many of the results on the structure theory of Lie superalgebra theory can be extended to left Leibniz superalgebras.

Let $(\mathcal{A}, \bullet)$ be a left (resp. right) Leibniz superalgebra. We denote
$$
\operatorname{Leib}(\mathcal{A})=\operatorname{span}\{u \bullet v+(-1)^{|u||v|}v \bullet u, \text{ where } u, v \in \mathcal{A}\}.
$$
According to identity \eqref{mul}, we have
$
\mathrm{L}^\bullet_{u \bullet v} = -(-1)^{|u||v|} \mathrm{L}^\bullet_{v \bullet u}.
$ Therefore, $
\mathrm{L}^\bullet_{u \bullet v + (-1)^{|u||v|} v \bullet u} = 0,
$ which implies that $\mathrm{L}^\bullet_w = 0$,  (resp $\Rr_w^\bullet=0$), for any $w \in \mathrm{Leib}(\mathcal{A})$. It follows that $\mathrm{Leib}(\mathcal{A})$ is a graded ideal of the Leibniz superalgebra $(\mathcal{A}, \bullet)$, and $\mathrm{Lie}(\A) := \A/\mathrm{Leib}(\A)$ is a Lie superalgebra. Moreover, $\A$ is a Lie superalgebra if and only if $\mathrm{Leib}(\A) = \{0\}$. 

We call $\mathrm{Leib}(\A)$ the Leibniz ideal of $\A.$

\sssbegin{Definition}\label{left-sym}
 A left-symmetric superalgebra (LSSA for short) is a superalgebra $(\A , \bullet)$ over a field $\mathbb{K}$, equipped with a bilinear product $\bullet$, such that the associator
$
(u,v,w) := (u\bullet v) \bullet w - u \bullet (v \bullet w)
$
is symmetric in the first two arguments. That is,
\begin{align}
(u,v,w) = (-1)^{|u||v|}(v, u, w),\label{leftsym}\end{align}
for all $u,v,w \in \A$
\end{Definition}
Let $(\mathcal{A}, \bullet)$ be a left-symmetric superalgebra. Then, it is clear that the supercommutator
$$
[u,v] := u \bullet v - (-1)^{|u||v|} v \bullet u
$$
defines a Lie superalgebra structure on $\mathcal{A}$. Moreover, the  condition
\eqref{leftsym} is equivalent to
$$
\mathrm{L}^\bullet_{[u,v]} = [\mathrm{L}^\bullet_u, \mathrm{L}^\bullet_v],
$$
for all $u, v \in \mathcal{A}$.
\sssbegin{Definition}
 Let $(\A, \bullet)$ be a non-associative superalgebra, $V$ be a $\mathbb{Z}_2$-graded vector space and $r, l$ : $\A \rightarrow \operatorname{End}(V)$ be two even linear maps. If $(\A, \bullet)$ is a left (resp. right) Leibniz superalgebra, then we say that $(r, l)$ is a left (resp. right) representation of $\A$ in $V$ if for all $u,v\in \A$ :
$$
\begin{gathered}
l(u\bullet v)=[l(u), l(v)] ; \quad r(u\bullet v)=l(u) r(v)+(-1)^{|u||v|} r(v) r(u) ; \quad r(u\bullet v)=[l(u), r(v)] \\
\text { (resp. } \quad l(u\bullet v)=(-1)^{|u||v|}[r(v), l(u)] ; \quad l(u\bullet v)=l(u) l(v)+(-1)^{|u||v|} r(v) l(u) \\
\left.r(u\bullet v)=(-1)^{|u||v|}[r(v), r(u)]\right)
\end{gathered}
$$
We say that $(r, l)$ is a representation of $\A$ in $V$ if $(r, l)$ is a left and a right representation of $\A$ in $V$.
\end{Definition}
\sssbegin{Proposition}[\cite{BenFa}]
Let $(\A, \bullet)$ be a left (resp. right) Leibniz superalgebra, and let $
r,\, l : L \to \operatorname{End}(V)
$ be two even linear maps. Then, the $\mathbb{Z}_2$-graded vector space $\A \oplus V$, endowed with the product
$$
(u+ x)\bullet (v + y) := u\bullet v + l(u)(y) + (-1)^{|u||v|} r(v)(w), \quad \forall\, u + w,\, v + y \in \A_{|u|} \oplus V_{|u|},
$$
is a left (resp. right) Leibniz superalgebra if and only if $(r, l)$ is a left (resp. right) representation of $(\A, \bullet)$ on $V$.
\end{Proposition}

\sssbegin{Proposition}[\cite{BenFa}]\label{represantaion} 
Let $(\A, \bullet)$ be a left (respectively, right) Leibniz superalgebra, and let $(r,l)$ denote its right (respectively, left) adjoint representation. Consider the even linear maps
$$
l^*(u)(f) = (-1)^{|f||u|} f \circ r(u), \quad r^*(u)(f) = (-1)^{|f||u|} f \circ l(u),
$$
for all  $u \in \A_{|u|}$, $f \in \A^*_{|f|}$, where $\A^*$ is the dual superspace of $\A$.
Then, $(r^*, l^*)$ is a left (respectively, right) representation of $(\A, \bullet)$ in $\A^*$ if and only if $(\A, \bullet)$ is a right (respectively, left) Leibniz superalgebra.
\end{Proposition}
\sssbegin{Definition}
Let $(\A, \bullet)$ be a left (respectively, right) Leibniz superalgebra, $V$ a $\mathbb{Z}_2$-graded vector space, and $(r, l)$ a left (respectively, right) representation of $\A$ on $V$. Let
$
\Om : \A\times \A \to V
$
be an even bilinear map. Then $\Om$ is called an even 2-cocycle of the left (respectively, right) Leibniz superalgebra $\A$ relative to $(r, l)$ if the following identity holds for all  $u,v,w \in \A$ 
\[
\begin{aligned}
&\Omega(u\bullet v, w) - \Omega(u, v\bullet w) + (-1)^{|u||v|} \Omega(v, u\bullet w) \\
&\quad - l(u)(\Omega(v, w)) + (-1)^{|u||v|}l(v)(\Omega(u, w)) + (-1)^{|w|(|u| + |v|)} r(w)(\Omega(u, v)) = 0 \\
&\text{(respectively, } 
 \Omega(u\bullet v, w) + \Omega(u, v\bullet w) + (-1)^{|v||w|} \Omega(u\bullet w, v) \\
&\quad - l(u)(\Omega(v, w)) - (-1)^{|u||v|} r(v)(\Omega(u, w)) + (-1)^{|w|(|u| + |v|)} r(w)(\Omega(u, v)) = 0\text{).}
\end{aligned}
\]
In both cases, such a map $\Om$ is called an even Leibniz 2-cocycle of $(\A, \bullet)$ relative to the representation $(r, l)$.
\end{Definition}

\section{Flat pseudo-Euclidean left-Leibniz superalgebras }

\subsection{Pre-Leibniz superalgebras}

The notion of a pre-Leibniz algebra was first introduced in \cite{Tang} where it was referred to as a Leibniz-dendriform algebra. We now introduce the notion of a pre-Leibniz superalgebra.
\sssbegin{Definition}\label{pre-Leibniz}
\begin{enumerate}
    \item[(i)] A disuperalgebra $(\A, \star, \circ)$ is called a pre-left-Leibniz
superalgebra (PLLA for short) if the following identities are satisfied (for all $u,v,w\in \A$)
\begin{align}\label{id1}
    (u\star v)\star w-u\star(v\star w)&=-(-1)^{|u||v|}((v\circ u)\star w+v\star(u\star w)),\\ \label{id2}
    (u\star v)\circ w-u\star(v\circ w)&=-(-1)^{|u||v|}((v\circ u)\circ w-v\circ(u\circ w)),\\ \label{id3}
    u\circ (v\circ w)&=-u\circ(v\star w).
\end{align}

\item[(ii)] A disuperalgebra $(\A, \star, \circ)$ is called a pre-right-Leibniz
superalgebra (PRLA for short) if the following identities are satisfied (for all $u,v,w\in \A$)
\begin{align}
\label{id4}
(u \star v) \circ w + u \circ (v \circ w) &= -(-1)^{|u||v|} \left((v \circ u) \circ w - v \circ (u \circ w)\right), \\
\label{id5}
(u \star v) \star w - u \star (v \star w) &= -(-1)^{|u||v|} \left((v \circ u) \star w - v \circ (u \star w)\right), \\
\label{id6}
 u \star (v \star w) &= -u \star (v\circ w), \end{align}
\item[(iii)] If $(\A, \star, \circ)$ is both pre-left and pre-right Leibniz then it is called a pre-symmetric Leibniz superalgebra.
\end{enumerate}
\end{Definition}
\sssbegin{Remark}\rm{ \begin{enumerate}
    \item[(i)] If $(\A, \star, \circ)$  is a pre-left (resp. right)-Leibniz superalgebra, then $\A$ can be endowed with the new products \[
    u\star^0 v: =(-1)^{|u||v|}u\circ v \text{ and  } u\circ^0 v : =(-1)^{|u||v|}u\star v.
    \] Then $(\A, \star^0,\circ^0)$ is a pre-right (resp. left)-Leibniz superalgebra.
\item[(ii)] If   $\star=-\circ$, then  $(\A, \star)$ is a left-symmetric superalgebra, i.e.,  $$(u\star v)\star w-u\star(v\star w)=(-1)^{|u||v|}((v\star u)\star w-v\star(u\star w)),$$
for all $u,v,w\in \A.$ 

\end{enumerate}}\end{Remark}
\sssbegin{Definition}
  Let $(\A, \star,\circ)$ be a disuperalgebra. On the underlying $\Z_2$-graded vector space $\A$, we
define a new product as follows:  
\begin{equation} \label{mix-prod}
    u\bullet v :=u\star v + (-1)^{|u||v|}v\circ u,\; u,v\in\A.
    \end{equation}
 The disuperalgebra $(\A, \star,\circ)$ is called left (right)-Leibniz-admissible disuperalgebra, if $(\A, \bullet)$ is a left (right)-Leibniz superalgebra.
\end{Definition}

\sssbegin{Proposition}\label{prLeibniz}
  Let $(\A, \star,\circ)$ be a pre-left \textup{(}resp.  right\textup{)}-Leibniz superalgebra. Then $(\A,\bullet)$ is left \textup{(}resp.  right\textup{)}-Leibniz superalgebra where the product $\bullet$ is defined as in \textup{(}\ref{mix-prod}\textup{)}.  In
particular, pre-left \textup{(}resp.  right\textup{)}-Leibniz superalgebras are left-Leibniz-admissible superalgebras. 
\end{Proposition}
\begin{proof}
 Assume that $(\A, \star, \circ)$ is a pre-left-Leibniz superalgebra.  For any $u,v,w\in \A$, we have 
   \begin{align*}
       (u&\bullet v)\bullet w-u\bullet (v\bullet w)+(-1)^{|u||v|} v\bullet (u\bullet w)\\=& (u\star v+ (-1)^{|u||v|} v\circ u)\star w +(-1)^{|w|(|u|+|v|)} w\circ (u\star v+ (-1)^{|u||v|} v\circ u)- u\star (v\star w+(-1)^{|v||w|} w\circ v) \\& -(-1)^{|u|(|v|+|w|) }(v\star w+(-1)^{|v||w|} w\circ v)\circ u+(-1)^{|u||v|}v\star (u\star w+(-1)^{|u||w|} w\circ u) \\& +(-1)^{|v|(|u|+|w|)+|u||v| }(u\star w+(-1)^{|u||w|} w\circ u)\circ v\\=& \left( (u\star v)\star w-u\star(v\star w)+(-1)^{|u||v|}(v\circ u)\star w+v\star(u\star w))\right)-\left( (-1)^{|u|(|v|+|w|) }(v\star w)\circ u-v\star (w\circ u)\right)\\&- (-1)^{|u|(|v|+|w|) }\left((-1)^{|v|w|}((w\circ v)\circ u-w\circ (v\circ u))\right)+ (-1)^{|w|(|u|+|v|) }\left((w\circ u)\circ v+w\circ((u\star v)\right)\\&+ (-1)^{|u||w| }((u\star w)\circ v-u\star(v\star w))=0.
\end{align*}
Therefore, $(\A, \bullet)$ is a left-Leibniz superalgebra.  

The case where $(\A, \star, \circ)$ is a pre-right-Leibniz superalgebra can be treated similarly by verifying the identity corresponding to right-Leibniz superalgebras.
\end{proof}
\sssbegin{Remark}
Let $(\A, \star, \circ)$ be a pre-left-Leibniz superalgebra. Then, the product defined by 
\[
u \bullet v := v \circ u + (-1)^{|u||v|} u \star v, \quad \text{for all } u,v \in \A,
\]
is a right-Leibniz superalgebra. Furthermore, if $(\A, \star, \circ)$ is both a pre-left-Leibniz and a pre-right-Leibniz superalgebra, then $(\A, \bullet)$ is a symmetric Leibniz superalgebra.
\end{Remark}

\sssbegin{Proposition}
  Let $(\A, \star,\circ)$ be a pre-left-Leibniz superalgebra. Then $(\Ll^\star, \Ll^\circ)$ is a left representation of $\A$ in $\A$\end{Proposition}
  \begin{proof}
   Let $u,v,w\in \A$. According to \eqref{id1},\eqref{id2} and \eqref{id3}, we have
   \begin{align*}
     \Ll_{u\bullet v}^\star =[\Ll_u^\star, \Ll_v^\star],\quad  \Ll_{u\bullet v}^\circ =\Ll_u^\star\circ\Ll_v^\circ+(-1)^{|u||v|}\Ll_v^\circ\circ\Ll_u^\circ,\quad \Ll_{u\bullet v}^\circ =[\Ll_u^\star,\Ll_v^\circ].
   \end{align*}
   Therefore, $(\Ll^\star, \Ll^\circ)$  is a left representation of $\A$ in $\A$.
\end{proof}
\subsection{The Levi-Civita product}
In this section, we introduce the notion of the Levi-Civita product in the setting of pseudo-Euclidean non-associative superalgebras. This construction extends the classical Levi-Civita product developed for pseudo-Euclidean Lie superalgebras in \cite{BenFa}, and for pseudo-Euclidean  Lie color algebras in \cite{color}. Our approach is also a generalization of the Levi-Civita product introduced for pseudo-Euclidean non-associative algebras in \cite{AbBoEl}. 

We now introduce the notion of the Levi-Civita product in the setting of pseudo-Euclidean non-associative superalgebras.
\sssbegin{Proposition}\label{LVCT} Let $(\A,\bullet,\prs)$ be a pseudo-Euclidean non-associative superalgebra. Then there exists a unique products $(\star, \circ)$ on $\A$ satisfying
	\begin{align} \label{torsion}
		u\bullet v&=u\star v+(-1)^{|u||v|}v\circ u
	\\ \label{compatibl} \langle u\star v,w\rangle&=(-1)^{|u||v|}\langle v, u\circ w\rangle.\end{align}More precisely, the product $\star$ is defined by
	\begin{equation}\label{lv}
		2\langle u\star v,w\rangle=\langle u\bullet v,w\rangle -(-1)^{|u||v|+|w||u|}\langle v\bullet w,  u\rangle+(-1)^{|v||w|+|u||w|}\langle w\bullet u,  v\rangle.
	\end{equation}
	Moreover, if $\bullet$ is super anti-commutative (resp. super commutative) then $\star=-\circ$ (resp. $\star=\circ$). 
	The product $(\star,\circ)$ is called a Levi-Civita product associated with  $(\A,\bullet,\prs)$.	
\end{Proposition}

\begin{proof} Let $u,v,w\in \A$ we have
	\begin{align*}
		\langle u\star v,w\rangle&=\langle u\bullet v,w\rangle -(-1)^{|u||v|}\langle v\circ u,w\rangle\\
        &=\langle u\bullet v,w\rangle -(-1)^{|u||v|+|w||u|}\langle v\star w,  u\rangle\\
        &=\langle u\bullet v,w\rangle -(-1)^{|u||v|+|w||u|}\langle v\bullet w,  u\rangle+(-1)^{|u||v|+|w||u|+|v||w|}\langle w\circ v,  u\rangle\\&=\langle u\bullet v,w\rangle -(-1)^{|u||v|+|w||u|}\langle v\bullet w,  u\rangle+(-1)^{|v||w|+|u||w|}\langle w\star u,  v\rangle\\ &=\langle u\bullet v,w\rangle -(-1)^{|u||v|+|w||u|}\langle v\bullet w,  u\rangle+(-1)^{|v||w|+|u||w|}\langle w\bullet u,  v\rangle-(-1)^{|v||w|}\langle u\circ w,  v\rangle  \\ &=\langle u\bullet v,w\rangle -(-1)^{|u||v|+|w||u|}\langle v\bullet w,  u\rangle+(-1)^{|v||w|+|u||w|}\langle w\bullet u,  v\rangle-(-1)^{|u||v|}\langle v, u\circ w\rangle  \\ &=\langle u\bullet v,w\rangle -(-1)^{|u||v|+|w||u|}\langle v\bullet w,  u\rangle+(-1)^{|v||w|+|u||w|}\langle w\bullet u,  v\rangle-\langle u\star v,  w\rangle
	\end{align*}
	This proves the uniqueness and gives  \eqref{lv}. On the other hand,
	\begin{align*}
		2\langle u\star v+(-1)^{|u||v|}v\circ u,w\rangle=& 2\langle u\star v,w\rangle+ (-1)^{|u||v|}2\langle v\circ u,w\rangle\\=& 2\langle u\star v,w\rangle+ (-1)^{|u||v|+|w|(|u|+|v|)}2\langle w,v\circ u\rangle\\=& 2\langle u\star v,w\rangle +(-1)^{|u||v|+|w||u|}2\langle v\star w, u\rangle\\=& \langle u\bullet v,w\rangle -(-1)^{|u||v|+|w||u|}\langle v\bullet w,  u\rangle+(-1)^{|v||w|+|u||w|}\langle w\bullet u,  v\rangle\\&+  (-1)^{|u||v|+|w||u|}\langle v\bullet w,u\rangle -(-1)^{|u||v|+|w||u|+|v||w|+|u||v|}\langle w\bullet u,  v\rangle\\&+(-1)^{2|u||v|+2|u||w|}\langle u\bullet v,  w\rangle\\=& 2\langle u\bullet v,w\rangle
		\end{align*}

	The case where $\bullet$ is super anti-commutative. We have
	\begin{align*}
		& 2\langle u\star v,w\rangle +
		2(-1)^{|u||v|}\langle v, u\star w\rangle=2\langle u\star v,w\rangle +
		2(-1)^{|v||w|}\langle  u\star w,v\rangle\\&=
		 \langle u\bullet v,w\rangle -(-1)^{|u||v|+|w||u|}\langle v\bullet w,  u\rangle+(-1)^{|v||w|+|u||w|}\langle w\bullet u,  v\rangle+ (-1)^{|v||w|}\langle u\bullet w,v\rangle \\& -(-1)^{|v||w|+|u||w|+|u||v|}\langle w\bullet v,  u\rangle +(-1)^{|v||w|+|u||v|+|w||v|}\langle v\bullet u,  w\rangle\\&=
		 \langle u\bullet v,w\rangle -(-1)^{|u||v|+|w||u|}\langle v\bullet w,  u\rangle+(-1)^{|v||w|+|u||w|}\langle w\bullet u,  v\rangle\\&- (-1)^{|v||w|+|w||u|}\langle w\bullet u,v\rangle +(-1)^{|u||w|+|u||v|}\langle v\bullet w,  u\rangle-\langle u\bullet v,  w\rangle=0.
	\end{align*}
	The case where $\bullet$ is super commutative.  We have
	\begin{align*}
		& 2\langle u\star v,w\rangle -
		2(-1)^{|u||v|}\langle v, u\star w\rangle=2\langle u\star v,w\rangle -
		2(-1)^{|v||w|}\langle  u\star w,v\rangle\\&=
		 \langle u\bullet v,w\rangle -(-1)^{|u||v|+|w||u|}\langle v\bullet w,  u\rangle+(-1)^{|v||w|+|u||w|}\langle w\bullet u,  v\rangle\\&- (-1)^{|v||w|}\langle u\bullet w,u\rangle +(-1)^{|v||w|+|u||w|+|u||v|}\langle w\bullet v,  u\rangle -(-1)^{|v||w|+|u||v|+|w||v|}\langle v\bullet u,  w\rangle\\&=
		 \langle u\bullet v,w\rangle -(-1)^{|u||v|+|w||u|}\langle v\bullet w,  u\rangle+(-1)^{|v||w|+|u||w|}\langle w\bullet u,  v\rangle\\&- (-1)^{|v||w|+|w||u|}\langle w\bullet u,u\rangle +(-1)^{|u||w|+|u||v|}\langle v\bullet w,  u\rangle-\langle u\bullet v,  w\rangle=0. \qed
	\end{align*}
    \noqed
		\end{proof}

Let $(\A, \bullet, \prs)$ be a pseudo-Euclidean nonassociative superalgebra, and let $(\star,\circ)$ denote the Levi-Civita products associated with $(\A, \bullet, \prs)$. According to relations \eqref{compatibl} and \eqref{torsion}, we have, for any $u \in \A$, that $\mathrm{L}^\circ_u$ is the adjoint of $\Ll_u^\star$ with respect to $\prs$, and
\(
\Ll^\bullet_u = \mathrm{L}^\star_u + \mathrm{R}^\circ_u\), and  \(\Rr^\bullet_u = \Rr^\star_u + \Ll^\circ_u .
\) In particular:
\begin{itemize}
    \item[(i)] If $(\A, \bullet)$ is a superanticommutative superalgebra, then for any $u \in \A$,  $\mathrm{L}^\star_u$ is skew-symmetric with respect to $\prs$, and
    \(
    \Ll_u^\bullet = \mathrm{L}^\star_u - \mathrm{R}^\star_u.
    \)
    \item[(ii)] If $(\A, \bullet)$ is a supercommutative superalgebra, then for any $u \in \A$,  $\mathrm{L}^\star_u$ is symmetric with respect to $\prs$, and
    \(
    \Ll_u^\bullet = \mathrm{L}^\star_u + \mathrm{R}^\star_u.
    \)
\end{itemize}

\sssbegin{Definition}\label{courbures}
Let $(\A , \bullet, \prs)$ be a pseudo-Euclidean left-Leibniz superalgebra and let $(\star,\circ)$ denote its Levi-Civita products. The curvatures tensors of $(\A , \bullet,\prs)$ are defined as follows (for all $u,v,w\in\A$):
\begin{equation}
    \begin{cases}
        K_1(u,v):= \Ll_{u\bullet v}^\star-[\Ll_u^\star, \Ll_v^\star],\\
        K_2(u,v):= -\Ll_{u\bullet v}^\circ+\Ll_u^\star\circ\Ll_v^\circ+(-1)^{|u||v|}\Ll_v^\circ\circ\Ll_u^\circ,\\
        K_3(u,v):=\Ll_{u\bullet v}^\circ-[\Ll_u^\star,\Ll_v^\circ].
    \end{cases}
\end{equation}
\end{Definition}
\sssbegin{Remark}
Let $(\A, \br, \prs)$ be a pseudo-Euclidean Lie superalgebra and let $\star$ denote its Levi-Civita product. The curvature tensor of $(\A , \br, \prs)$ is given by $K:=K_1=K_2=-K_3.$  
\end{Remark}
\sssbegin{Proposition}
 Let $(\A , \bullet, \prs)$ be a pseudo-Euclidean left-Leibniz superalgebra and let $(\star,\circ)$ denote its Levi-Civita products. Then, for all  $u,v,w\in \A$, we have
 $$(u\bullet v)\bullet w-u\bullet (v\bullet w)+(-1)^{|u||v|}v\bullet (u\bullet w)=K_1(u,v)w+(-1)^{|u||v|+|u||w|}K_2(v,w)u+(-1)^{|v||w|}K_3(w,u)v.$$
\end{Proposition}
\begin{proof}
 Let $u,v,w\in \A$, we have 
 \begin{align*}
  (u\bullet & v)\bullet w-u\bullet (v\bullet w)+(-1)^{|u||v|}v\bullet (u\bullet w)\\=& (u\bullet v)\star w+(-1)^{|u||w|+|v||w|}w\circ (u\star v+(-1)^{|u||v|} v\circ u) \\& -u\star(v\star w+(-1)^{|v||w|} w\circ v)-(-1)^{|u||v|+|u||w|}(v\bullet w)\circ u\\& +(-1)^{|u||v|} v\star(u\star w+(-1)^{|u||w|} w\circ u)+(-1)^{2|u||v|+|v||w|}(u\bullet w)\circ v\\=& ((u\bullet v)\star w-u\star (v\star w)+(-1)^{|u||v|}v\star(u\star w))\\&+ (-1)^{|u||v|+|u||w|}(-(v\bullet w)\circ u   +v\star(w\circ u)  +(-1)^{|v||w|} w\circ (v\circ u))\\&+  (-1)^{|v||w|}((u\bullet w)\circ v   -u\star(w\circ v)  +(-1)^{|u||w|} w\circ (u\star v)))\\=& K_1(u,v)w+ (-1)^{|u||v|+|u||w|}K_2(v,w)u+ (-1)^{|u||w|} K_3(u,w)v.
 \end{align*}
\end{proof}
\sssbegin{Remark}
 Let $(\A , \bullet, \prs)$ be a pseudo-Euclidean Lie superalgebra and let $\star$ denote its Levi-Civita product. Then, for all  $u,v,w\in \A$, we have
 $$K(u,v)w+(-1)^{|u||v|+|u||w|}K(v,w)u+(-1)^{|v||w|+|u||w|}K(w,u)v=0,$$
since $K(u,v)=-(-1)^{|u||v|}K(v,u)$.
\end{Remark}
\sssbegin{Definition}
  A pseudo-Euclidean left Leibniz superalgebra $(\A, \bullet,\prs)$ is called \textit{flat} if $K_1=K_2=K_3=0$. In particular, if $(\A, \bullet, \prs)$ is pseudo-Euclidean Lie superalgebra, then it is called a flat pseudo-Euclidean Lie superalgebra, if $K=0$. 
\end{Definition}

\sssbegin{Remark}
Let $(\A, \star, \circ)$ be a pre-left-Leibniz superalgebra equipped with a non-degenerate symmetric bilinear form $\prs$ satisfying
\[
\langle u \star v, w \rangle = (-1)^{|u||v|} \langle v, u \circ w \rangle,
\]
for all $u, v, w \in \A$. According to Proposition \ref{prLeibniz}, $(\A, \star, \circ)$ is a left-Leibniz-admissible superalgebra, and the associated product $\bullet$ on $\A$ is given by
\[
u \bullet v = u \star v + (-1)^{|u||v|} v \circ u, \quad \forall u,v\in \A.
\]
Since the bilinear form $\prs$ also satisfies the compatibility condition
\[
\langle u \star v, w \rangle = (-1)^{|u||v|} \langle v, u \circ w \rangle,
\]
the uniqueness of the Levi-Civita product implies that the pair $(\star, \circ)$ coincides with the Levi-Civita product associated with the pseudo-Euclidean left-Leibniz superalgebra $(\A, \bullet, \prs)$. In this case, $(\A, \bullet, \prs)$ is flat.

Conversely, if $(\A, \bullet, \prs)$ is a flat pseudo-Euclidean left-Leibniz superalgebra, then its Levi-Civita product defines a pre-left-Leibniz superalgebra structure on $\A$.
\end{Remark}

\sssbegin{Example}{\rm 
Let \( \A = \A_{\bar 0} \oplus \A_{\bar 1} \) be a six-dimensional left-Leibniz superalgebra, where \( \A_0 = \mathrm{span}\{e_1, e_2, e_3, e_4\} \) and \( \A_1 = \mathrm{span}\{f_1, f_2\} \). The left-Leibniz product \( \bullet \) is defined on the basis elements by:
\[
\begin{aligned}
    &e_1 \bullet e_1 = e_2 + e_3, \quad e_4 \bullet e_4 = e_2, \quad f_2 \bullet e_4 = -e_4 \bullet f_2 = 2f_1, \\
    &f_2 \bullet f_2 = e_2 - 2e_3, \quad e_1 \bullet e_4 = e_4 \bullet e_1 = e_1 + e_3, \quad e_1 \bullet f_2 = f_2 \bullet e_1 = f_1.
\end{aligned}
\]

Consider the even symmetric bilinear form \( \prs : \A \times \A \to \mathbb{K} \) defined by:
\[
\langle e_1, e_2 \rangle = 1, \quad \langle e_3, e_4 \rangle = 1, \quad \langle f_1, f_2 \rangle = 1,
\]
with all other pairings between basis elements being zero. 

A straightforward computation shows that the Levi-Civita products associated with \( (\A, \bullet, \prs) \) are given by:
\[
\begin{aligned}
    &e_1 \star e_1 = \frac{1}{2}(e_2 + e_3), \quad e_4 \star e_4 = \frac{1}{2}e_2, \quad f_2 \star e_4 = -e_4 \star f_2 = f_1, \\
    &f_2 \star f_2 = \frac{1}{2}e_2 - e_3, \quad e_1 \star e_4 = e_4 \star e_1 = \frac{1}{2}(e_1 + e_3), \quad e_1 \star f_2 = f_2 \star e_1 = \frac{1}{2}f_1,
\end{aligned}
\]
and
\[
\begin{aligned}
    &e_1 \circ e_1 = \frac{1}{2}(e_2 + e_3), \quad e_4 \circ e_4 = -\frac{1}{2}e_2, \quad f_2 \circ e_4 = -e_4 \circ f_2 = f_1, \\
    &f_2 \circ f_2 = -\frac{1}{2}e_2 + e_3, \quad e_1 \circ e_4 = e_4 \circ e_1 = \frac{1}{2}(e_1 + e_3), \quad e_1 \circ f_2 = f_2 \circ e_1 = \frac{1}{2}f_1.
\end{aligned}
\]
Hence, by direct computation, we conclude that \( (\A, \bullet, \prs) \) is a flat pseudo-Euclidean left-Leibniz superalgebra.
}
\end{Example}

\sssbegin{Example}{\rm 
Let \( \mathcal{A} = \A_{\bar 0} \oplus \A_{\bar 1} \) be a four-dimensional left-Leibniz superalgebra, where \( \A_0 = \mathrm{span}\{e_1, e_2\} \) and \( \A_1 = \mathrm{span}\{f_1, f_2\} \). The left-Leibniz product \( \bullet \) is defined on the basis elements by:
\[
\begin{aligned}
    e_2\bullet e_2=e_1,\quad e_2 \bullet f_1 =f_1\bullet e_2= f_2, \quad f_1 \bullet f_1 = e_1,
\end{aligned}
\]

Consider the odd symmetric bilinear form \( \prs : \mathcal{A} \times \mathcal{A} \to \mathbb{K} \) defined by:
\[
\langle e_1, f_1 \rangle = 1, \quad \langle e_2, f_2 \rangle = 1, 
\]
 A straightforward computation shows that the Levi-Civita products associated with \( (\A, \bullet, \prs) \) are given by:
\[
\begin{aligned}
    e_2\star e_2=\frac{1}{2}e_1,\quad e_2 \star f_1 =f_1\star e_2= \frac{1}{2}f_2, \quad f_1 \star f_1 = \frac{1}{2}e_1,
\end{aligned}
\]
and
\[
\begin{aligned}
    e_2\circ e_2=\frac{1}{2}e_1,\quad e_2 \circ f_1 =f_1\circ e_2= \frac{1}{2}f_2, \quad f_1 \circ f_1 = \frac{1}{2}e_1,
\end{aligned}
\]
Hence, by direct computation, we conclude that \( (\mathcal{A}, \bullet, \prs) \) is a flat pseudo-Euclidean left-Leibniz superalgebra.
}
\end{Example}

\sssbegin{Proposition}\label{flat}
	Let $(\A, \bullet, \prs)$ be a  pseudo-Euclidean left Leibniz superalgebra. Then $\A$ is flat if and only if 
	\begin{equation}\label{equivf}
		\Ll_{u\bullet v}^\star=[\Ll_u^\star,\Ll_v^\star], \quad \text{ and }\quad \Ll_u^{\star}\circ \Ll_v^{\star}=-\Ll_u^{\star}\circ \Ll_v^{\circ}=-\Ll_u^{\circ}\circ\Ll_v^{\star}=\Ll_u^{\circ}\circ \Ll_v^{\circ},\end{equation}
	for any $u,v\in \A$.
\end{Proposition}

\begin{proof}We suppose that $\A$ is flat that is $K_1=K_2=K_3=0$.  In this case, the following identities hold:
\begin{align}
\label{re1} \Ll_{u\bullet v}^\star &= [\Ll_u^\star, \Ll_v^\star], \\
\label{re2} \Ll_{u\bullet v}^\circ &= \Ll_u^\star \circ \Ll_v^\circ +(-1)^{|u||v|} \Ll_v^\circ \circ \Ll_u^\circ, \\
\label{re3} \Ll_{u\bullet v}^\circ &= [\Ll_u^\star, \Ll_v^\circ].
\end{align}
Combining equations \eqref{re2} and \eqref{re3}, we deduce that
$$
\Ll_v^\circ \circ \Ll_u^\star = - \Ll_v^\circ \circ \Ll_u^\circ.
$$
Now, applying the adjoint of equation \eqref{re2}, we obtain for any homogeneous elements $u,v,w,z \in \A$:
\begin{align*}
&\left\langle \left(\Ll_{u \bullet v}^\circ - \Ll_u^\star \circ \Ll_v^\circ - (-1)^{|u||v|} \Ll_v^\circ \circ \Ll_u^\circ\right)(w), z \right\rangle \\
&= (-1)^{|w|(|u|+v|)} \left\langle w, \Ll_{u \bullet v}^\star(z) \right\rangle
- (-1)^{|u|(|w|+|v|) + |v||w|} \left\langle w, \Ll_v^\star \circ \Ll_u^\circ(z) \right\rangle
- (-1)^{|v||w| + |u||w|} \left\langle w, \Ll_u^\star \circ \Ll_v^\star(z) \right\rangle \\
&= (-1)^{|w|(|u|+|v|)} \left\langle w, \left( \Ll_{u \bullet v}^\star - (-1)^{|u||v|} \Ll_v^\star \circ \Ll_u^\circ - \Ll_u^\star \circ \Ll_v^\star \right)(z) \right\rangle.
\end{align*}
Hence, by the non-degeneracy of the bilinear form, we conclude that
$$
\Ll_{u \bullet v}^\star = (-1)^{|u||v|} \Ll_v^\star \circ \Ll_u^\circ + \Ll_u^\star \circ \Ll_v^\star.
$$
Substituting this into \eqref{re1}, we obtain
$$
\Ll_v^\star \circ \Ll_u^\circ = - \Ll_v^\star \circ \Ll_u^\star.
$$
Taking adjoints again yields
$$
\Ll_u^\star \circ \Ll_v^\circ = - \Ll_u^\circ \circ \Ll_v^\circ.
$$
From these relations, we deduce that for any homogeneous elements $u,v \in \A$, the following identity holds:
$$
\Ll_u^\star \circ \Ll_v^\star = - \Ll_u^\star \circ \Ll_v^\circ = - \Ll_u^\circ \circ \Ll_v^\star = \Ll_u^\circ \circ \Ll_v^\circ.
$$
Conversely, it is obvious that the equations \eqref{equivf} implies that  $K_1=K_2=K_3=0$ which completes the proof.
\end{proof}
\sssbegin{Proposition}\label{symmetricLeibniz}
    Let $(\A, \bullet, \prs)$ be a flat pseudo-Euclidean left-Leibniz superalgebra, then $(\A, \bullet)$ is symmetric Leibniz superalgebra.
\end{Proposition}
\begin{proof}
Assume that $(\A, \bullet, \prs)$ is a flat pseudo-Euclidean left-Leibniz superalgebra. According to Proposition \ref{flat}, we have the following identities:
\begin{equation*}
\Ll_{u\bullet v}^\star = [\Ll_u^\star, \Ll_v^\star] \quad \text{and} \quad \Ll_u^{\star} \circ \Ll_v^{\star} = -\Ll_u^{\star} \circ \Ll_v^{\circ} = -\Ll_u^{\circ} \circ \Ll_v^{\star} = \Ll_u^{\circ} \circ \Ll_v^{\circ},
\end{equation*}
for all homogeneous elements $u, v \in \A$.
From these identities, we obtain:
\begin{align*}
\Ll_{u\bullet v}^\star &= \Ll_u^\star \circ \Ll_v^\star - (-1)^{|u||v|} \Ll_v^\star \circ \Ll_u^\star = \Ll_u^\star \circ \Ll_v^\star + (-1)^{|u||v|} \Ll_v^\circ \circ \Ll_u^\star = -\Ll_u^\circ \circ \Ll_v^\star + (-1)^{|u||v|} \Ll_v^\circ \circ \Ll_u^\star.
\end{align*}
On the other hand, taking the adjoint of the second identity yields:
$$
\Ll_{u\bullet v}^\circ = -[\Ll_u^\circ, \Ll_v^\circ],
$$
which shows that $(\A, \star, \circ)$ also satisfies the pre-right-Leibniz superalgebra condition.
Thus, by Proposition \ref{prLeibniz}, $(\A, \bullet)$ is a right-Leibniz superalgebra. Since it is also left-Leibniz superalgebra, we conclude that $(\A, \bullet)$ is a symmetric Leibniz superalgebra.
\end{proof}

Let $(\mathcal{A}, \bullet, \prs)$ be a pseudo-Euclidean nonassociative superalgebra, and let $(\star, \circ)$ denote the Levi-Civita products associated with $(\mathcal{A}, \bullet, \prs)$. Define the underlying superanticommutative and supercommutative superalgebras of $\A$ by
$$
[u,v] := \tfrac{1}{2}(u \bullet v - (-1)^{|u||v|} v \bullet u), \quad \{u,v\} := \tfrac{1}{2}(u \bullet v + (-1)^{|u||v|} v \bullet u),
$$
for all $u,v \in \A$, so that
\begin{equation}\label{pralg}
u \bullet v = [u,v] + \{u,v\}.
\end{equation}
We denote by $\mathcal{A}^- := (\mathcal{A}, \br)$ the superanticommutative superalgebra associated with $(\mathcal{A}, \bullet)$, and by $\mathcal{A}^+ := (\mathcal{A}, \{\,,\,\})$ the supercommutative superalgebra associated with $(\mathcal{A}, \bullet)$.

Since $(\mathcal{A}^-, \prs)$ is a pseudo-Euclidean superanticommutative superalgebra and $(\mathcal{A}^+, \prs)$ is a pseudo-Euclidean supercommutative superalgebra, it follows from Proposition~\ref{LVCT} that there exists a unique Levi-Civita product $\rhd$ (respectively, $\lhd$) associated with the superalgebra $(\mathcal{A}^-, \prs)$ (respectively, $(\mathcal{A}^+, \prs)$).

For any $u, v, w \in \A$, we have
\begin{align*}
2\langle u \rhd v + u \lhd v, w \rangle &= \langle [u,v], w \rangle - (-1)^{|u||v|+|u||w|} \langle [v,w], u \rangle + (-1)^{|u||w|+|v||w|} \langle [w,v], u \rangle \\
&\quad + \langle \{u,v\}, w \rangle - (-1)^{|u||v|+|u||w|} \langle \{v,w\}, u \rangle + (-1)^{|u||w|+|v||w|} \langle \{w,v\}, u \rangle \\
&= \langle u \bullet v, w \rangle - (-1)^{|u||v|+|u||w|} \langle v \bullet w, u \rangle + (-1)^{|u||w|+|v||w|} \langle w \bullet v, u \rangle \\
&= 2\langle u \star v, w \rangle.
\end{align*}
Hence, we conclude that
$$
u \rhd v + u \lhd v = u \star v,
$$
and thus the left and right multiplication operators of the Levi-Civita product satisfy
$$
\Ll_u^{\star} = \Ll_u^{\rhd} + \Ll_u^{\lhd},\quad \Ll_u^{\circ} = -\Ll_u^{\rhd} + \Ll_u^{\lhd} \quad \Rr_u^{\star} = \Rr_u^{\rhd} + \Rr_u^{\lhd}, \quad\Rr_u^{\circ} = -\Rr_u^{\rhd} + \Rr_u^{\lhd}\quad \forall u \in \A.
$$

\sssbegin{Theorem}\label{chara}

Let $(\A, \bullet, \prs)$ be a pseudo-Euclidean left-Leibniz superalgebra. Denote by $\rhd$ \textup{(}resp. $\lhd$\textup{)} the Levi-Civita product associated with the super-anticommutative (resp. supercommutative) superalgebra $(\A^-, \prs)$ (resp. $(\A^+, \prs)$). Then $(\A, \bullet, \prs)$ is flat if and only if  $(\A^-, \prs)$ is flat pseudo-Euclidean Lie superalgebra and the following relations are satisfied for all $u, v \in \A$:
$$
\Ll^{\lhd}_u \circ \Ll^{\lhd}_v = \Ll^{\lhd}_u \circ \Ll^{\rhd}_v = \Ll^{\rhd}_u \circ \Ll^{\lhd}_v = 0,
\quad \text{and} \quad
\Ll^{\rhd}_{\{u , v\}} = \Ll^{\lhd}_{\{u , v\}} = \Ll^{\lhd}_{[u,v]} = 0.
$$
\end{Theorem}

\begin{proof}
Assume that $(\A, \bullet, \prs)$ is a flat pseudo-Euclidean left-Leibniz superalgebra. According to~\eqref{equivf}, we have, for all $u, v \in \A$,
\begin{align*}
0 &= \Ll^{\star}_u \circ \Ll^{\star}_v + \Ll^{\star}_u \circ \Ll^{\circ}_v
= \Ll^{\star}_u \circ (\Ll^{\star}_v + \Ll^{\circ}_v)
= 2(\Ll^{\rhd}_u \circ \Ll^{\lhd}_v + \Ll^{\lhd}_u \circ \Ll^{\lhd}_v), \\
0 &= \Ll^{\star}_u \circ \Ll^{\star}_v + \Ll^{\circ}_u \circ \Ll^{\star}_v
= (\Ll^{\star}_u + \Ll^{\circ}_u) \circ \Ll^{\star}_v
= 2(\Ll^{\lhd}_u \circ \Ll^{\rhd}_v + \Ll^{\lhd}_u \circ \Ll^{\lhd}_v), \\
0 &= \Ll^{\circ}_u \circ \Ll^{\circ}_v + \Ll^{\star}_u \circ \Ll^{\circ}_v
= (\Ll^{\circ}_u + \Ll^{\star}_u) \circ \Ll^{\circ}_v
= 2(\Ll^{\lhd}_u \circ \Ll^{\lhd}_v - \Ll^{\lhd}_u \circ \Ll^{\rhd}_v).
\end{align*}
These relations imply that
\begin{equation}
\Ll^{\lhd}_u \circ \Ll^{\lhd}_v = \Ll^{\lhd}_u \circ \Ll^{\rhd}_v = \Ll^{\rhd}_u \circ \Ll^{\lhd}_v = 0, \quad \forall u,v \in \A.
\label{T1}\end{equation}
On the other hand, observe that
$$
\Ll_{\{u , v\}}^{\star} =\frac{1}{2}(\Ll_{u\bullet v}^\star+(-1)^{|u||v|} \Ll_{u\bullet v}^\star)=\frac{1}{2}([\Ll_u^\star, \Ll_v^\star] +(-1)^{|u||v|} [\Ll_v^\star, \Ll_u^\star])=0. 
$$
Moreover, since $\Ll_{\{u , v\}}^{\star} = \Ll_{\{u , v\}}^{\rhd} + \Ll_{\{u , v\}}^{\lhd} = 0$, and  $\Ll^{\rhd}_u$ and $\Ll^{\lhd}_u$ are respectively anti-symmetric and symmetric with respect to $\prs$, it follows that
\[
\Ll_{\{u , v\}}^{\rhd} = \Ll_{\{u , v\}}^{\lhd} = 0.
\]
Next, we have
\begin{align*}
\Ll_{[u, v]}^\star &= \tfrac{1}{2}\left( \Ll_{u\bullet v}^\star - (-1)^{|u||v|} \Ll_{v\bullet u}^\star \right) \\
&= \tfrac{1}{2}\left( [\Ll_u^\star, \Ll_v^\star] - (-1)^{|u||v|}[\Ll_v^\star, \Ll_u^\star] \right) \\
&= [\Ll_u^\star, \Ll_v^\star] \\
&= [\Ll^{\rhd}_u, \Ll^{\rhd}_v] + [\Ll^{\rhd}_u, \Ll^{\lhd}_v] + [\Ll^{\lhd}_u, \Ll^{\rhd}_v] + [\Ll^{\lhd}_u, \Ll^{\lhd}_v] \\
&\overset{\eqref{T1}}{=} [\Ll^{\rhd}_u, \Ll^{\rhd}_v],
\end{align*}
Thus, we obtain:
\[
\Ll_{[u, v]}^\star = \Ll_{[u, v]}^\rhd + \Ll_{[u, v]}^\lhd = [\Ll^{\rhd}_u, \Ll^{\rhd}_v].
\]
Since both $\Ll_{[u, v]}^\rhd$ and $[\Ll^{\rhd}_u, \Ll^{\rhd}_v]$ are anti-symmetric with respect to $\prs$, and $\Ll_{[u, v]}^\lhd$ is symmetric with respect to $\prs$, we conclude that
\[
\Ll_{[u, v]}^\rhd = [\Ll^{\rhd}_u, \Ll^{\rhd}_v], \qquad \Ll_{[u, v]}^\lhd = 0.
\]
which shows that $(\A^-,  \prs)$ is a flat pseudo-Euclidean Lie superalgebra. This completes the proof.
\end{proof}

\sssbegin{Proposition}\label{ideal}
 Let $(\A, \bullet, \prs)$ be a flat pseudo-Euclidean left-Leibniz superalgebra. Then the following assertions hold true:
\begin{enumerate}
\item[\textup{(}i\textup{)}] The subspaces $\operatorname{Leib}(\A)$ and $\operatorname{Leib}(\A)^\perp$ are two-sided  ideals of the Levi-Civita products $(\star, \circ)$.
\item[\textup{(}ii\textup{)}] The subspaces \(\operatorname{Leib}(\A) \cap \operatorname{Leib}(\A)^\perp\) and \(\operatorname{Leib}(\A) + \operatorname{Leib}(\A)^\perp\) are also two-sided  ideals of \(( \star, \circ)\).
\item[\textup{(}iii\textup{)}] If \(\operatorname{Leib}(\A)\) is non-degenerate with respect to \(\prs\), then \((\A, \bullet, \prs)\) is a flat pseudo-Euclidean Lie superalgebra. In particular, if \((\A, \bullet)\) is a non-Lie left-Leibniz superalgebra, then \(\operatorname{Leib}(\A) \cap \operatorname{Leib}(\A)^\perp \neq \{0\}\).

\item[\textup{(}iv\textup{)}] If $(\A, \bullet)$ is a non-Lie left-Leibniz algebra and there exists a non-zero element $u \in \operatorname{Leib}(\A) \cap \operatorname{Leib}(\A)^\perp$, then the subspace $I = \K u$ is a totally isotropic two-sided ideal of $(\A, \star, \circ)$, and its orthogonal complement $I^\perp$ is also a two-sided ideal of $(\A, \star, \circ)$.

\end{enumerate}
\end{Proposition}

\begin{proof}
(i) According to Proposition~\ref{symmetricLeibniz}, the left-Leibniz superalgebra $(\A, \bullet)$ is symmetric, which implies that $\Ll_u^\bullet = \Rr_u^\bullet = 0$ for all $u \in \operatorname{Leib}(\A)$. Since $(\A, \bullet, \prs)$ is flat, we have
$     \Ll_{v \bullet w}^\star = [\Ll_v^\star, \Ll_w^\star],
    $
and, in particular, this yields $\Ll_u^\star = 0$ for all $u \in \operatorname{Leib}(\A)$. By adjoint, we also obtain $\Ll_u^\circ = 0$. From the torsion identity~\eqref{torsion}, we get
\[
\Ll_u^\bullet = \Ll_u^\star + \Rr_u^\circ \quad \text{and} \quad \Rr_u^\bullet = \Rr_u^\star + \Ll_u^\circ.
\]
Hence, we conclude that
\[
\Ll_u^\bullet = \Rr_u^\bullet = \Ll_u^\star = \Ll_u^\circ = \Rr_u^\star = \Rr_u^\circ = 0,
\]
for all \(u \in \operatorname{Leib}(\A)\). This shows that \(\operatorname{Leib}(\A)\) and \(\operatorname{Leib}(\A)^\perp\) are two-sided ideals of \((\A, \star, \circ)\). Let \(u \in \operatorname{Leib}(\A)\), \(v \in \operatorname{Leib}(\A)^\perp\), and \(w \in \A\). Using the invariance of the bilinear form \(\prs\), we have:
\begin{align*}
    \langle v \star w, u \rangle &= (-1)^{|v||w|} \langle w, v \circ u \rangle = 0, &
    \langle w \star v, u \rangle &= (-1)^{|v||w|} \langle v, w \circ u \rangle = 0, \\
    \langle v \circ w, u \rangle &= (-1)^{|v||w|} \langle w, v \star u \rangle = 0, &
    \langle w \circ v, u \rangle &= (-1)^{|v||w|} \langle v, w \star u \rangle = 0.
\end{align*}
Therefore, \(\operatorname{Leib}(\A)^\perp\) is a two-sided ideal of \((\A, \star, \circ)\).

(ii) Follows directly from (i), since the intersection and sum of two-sided ideals remain two-sided ideals.

(iii)  Assume that \(\operatorname{Leib}(\A)\) is non-degenerate with respect to \(\prs\). Then, we have a direct sum decomposition:
\[
\A = \operatorname{Leib}(\A) \oplus \operatorname{Leib}(\A)^\perp.
\]
Let \(u, v \in \operatorname{Leib}(\A)^\perp\). Since \(\operatorname{Leib}(\A)^\perp\) is a two-sided ideal, we know that \(u \bullet v + (-1)^{|u||v|} v \bullet u \in \operatorname{Leib}(\A) \cap \operatorname{Leib}(\A)^\perp\). Hence,
\[
u \bullet v = -(-1)^{|u||v|} v \bullet u.
\]
In addition, since \(\Ll_w^\bullet = \Rr_w^\bullet = 0\), for any $w\in \operatorname{Leib}(\A)$ we deduce that \(w \bullet z=z\bullet w = 0\), for any $z\in \A$ and therefore the Leibniz kernel is trivial, 
\(
\operatorname{Leib}(\A) = \{0\}.
\)
This implies that \((\A, \bullet)\) is a Lie superalgebra.

(iv) Since $(\A, \bullet)$ is a non-Lie left-Leibniz superalgebra, there exists a non-zero element $u \in \operatorname{Leib}(\A) \cap \operatorname{Leib}(\A)^\perp$. By the proof of Part (i), we have $\Ll_u^\star = \Ll_u^\circ = \Rr_u^\star = \Rr_u^\circ = 0$. It follows that the subspace $I := \K u$ is a totally isotropic two-sided ideal of $(\A, \star, \circ)$. Moreover, it is straightforward to verify that its orthogonal complement $I^\perp$ is also a two-sided ideal of $(\A, \star, \circ)$.
\end{proof}

\section{Flat quadartic Leibniz superalgebras}

\subsection{Characterization of flat quadratic Leibniz superalgebras}
In this paragraph, we study of flat quadratic Leibniz superalgebra.  We show that a quadratic Leibniz superalgebra is flat if and only if it is a 2-step nilpotent symmetric Leibniz superalgebra. Furthermore, we provide a characterization of flat quadratic Leibniz superalgebra. 
\sssbegin{Definition}
    A triple $(\A, \bullet, \prs)$ is called a quadratic (resp. odd-quadratic) Leibniz superalgebra if
$(\A, \bullet)$ is a left (or right) Leibniz superalgebra equipped with a non-degenerate, symmetric, and invariant bilinear form $\prs$ which is even (resp. odd). In this case, the bilinear form $\prs$  is called an invariant (resp. odd-invariant) scalar product on $(\A, \bullet)$.
\end{Definition}
\sssbegin{Proposition}\label{LCquadratic}
  Let $(\A,\bullet, \prs)$ be a quadratic \textup{(}resp. odd-quadratic\textup{)} Leibniz superalgebra. Then the Levi-Civita product $(\star,\circ)$ is given by 
  $$u\star v=(-1)^{|u||v|}v\circ u=\frac{1}{2}u\bullet v,$$
  for all $u,v\in \A.$
\end{Proposition}
\begin{proof}
 Let $u,v,w\in \A$. We have 
 \begin{align*}
     2\langle u\star v,w\rangle=&\langle u\bullet v,w\rangle -(-1)^{|u||v|+|w||u|}\langle v\bullet w,  u\rangle+(-1)^{|v||w|+|u||w|}\langle w\bullet u,  v\rangle\\=& \langle u\bullet v,w\rangle -(-1)^{|u||v|+|w||u|}\langle v\bullet w,  u\rangle+(-1)^{2|v||w|+|u||w|+|v||u|}\langle v, w\bullet u\rangle \\=& \langle u\bullet v,w\rangle -(-1)^{|u||v|+|w||u|}\langle v\bullet w,  u\rangle+(-1)^{|u||w|+|v||u|}\langle v\bullet w,  u\rangle\\=& \langle u\bullet v,w\rangle.
 \end{align*}
  Since $\prs$ is non-degenerate, then we have $u\star v=\frac{1}{2}u\bullet v$. It follows that
  \[
  u\star v=\frac{1}{2}u\bullet v=\frac{1}{2}(u\star v+(-1)^{|u||v|} v\circ u),
  \]
  and hence $u\star v=(-1)^{|u||v|}v\circ u=\frac{1}{2}u\bullet v$.
\end{proof}
\sssbegin{Proposition}\label{pr1}
  Let $(\A,\bullet, \prs)$ be a quadratic (resp. odd-quadratic) Leibniz superalgebra. $(\A,\bullet, \prs)$ is flat if and only if  $(\A, \bullet)$ is 2-nilpotent symmetric Leibniz superalgebra.  
\end{Proposition}
\begin{proof}
According to Proposition \eqref{LCquadratic}, we have
$
u \star v = (-1)^{|u||v|} v \circ u = \frac{1}{2} u \bullet v, \quad \text{for all } u,v \in \A,
$ which implies that $\Ll_u^\star = \frac{1}{2} \Ll_u^\bullet$ and $\Ll_u^\circ = \frac{1}{2} \Rr_u^\bullet$. Assume that $(\A, \bullet, \prs)$ is flat, that is, its curvatures tensors satisfy $K_1 = K_2 = K_3 = 0$. Then we have:
$$
0 = \Ll_{u \bullet v}^\star - [\Ll_u^\star, \Ll_v^\star] = \frac{1}{2} \Ll_{u \bullet v}^\bullet - \frac{1}{4} [\Ll_u^\bullet, \Ll_v^\bullet].
$$
Thus, we get that $
\Ll_{u \bullet v}^\bullet = [\Ll_u^\bullet, \Ll_v^\bullet] = 0.
$  Since $\prs$ is an invariant (resp. odd-invariant) scalar product, we have $\Rr_u^\bullet = (\Ll_u^\bullet)^*$. It follows that
$
\Rr_{u \bullet v}^\bullet = \left( \Ll_{u \bullet v}^\bullet \right)^* = 0
$. Consequently, for all $u,v,w \in \A$, we have $$
(u \bullet v) \bullet w = w \bullet (u \bullet v) = 0.$$
Therefore, the superalgebra $(\A, \bullet)$ is symmetric Leibniz and 2-step nilpotent.

Conversely, assume that $(\A, \bullet)$ is symmetric Leibniz and 2-step nilpotent. Then, for all $u,v,w \in \A$, we have
$$
(u \bullet v) \bullet w = u \bullet (v \bullet w) = 0,
$$
which implies that all compositions of left and right multiplication operators vanish. In particular, using the identities $\Ll_u^\star = \frac{1}{2} \Ll_u^\bullet$ and $\Ll_u^\circ = \frac{1}{2} \Rr_u^\bullet$, one can deduce that
$
K_1 = K_2 = K_3 = 0,
$ and hence $(\A, \bullet, \prs)$ is flat.
\end{proof}

Let $(\mathcal{A}, \bullet)$ be a 2-step nilpotent symmetric Leibniz superalgebra. A straightforward computation shows that the anticommutative superalgebra $\mathcal{A}^-$ associated with $(\mathcal{A}, \bullet)$ is a 2-step nilpotent Lie superalgebra, and the commutative superalgebra $\mathcal{A}^+$ is a 2-step nilpotent associative superalgebra satisfying the identities
\begin{equation}
[u, \{v, w\}] = \{u, [v, w]\} = 0, \quad \text{for all } u, v, w \in \mathcal{A}.
\label{idsymmetric}
\end{equation}

Conversely, let $(\mathcal{A}, \br)$ be a 2-step nilpotent Lie superalgebra, and let $\mu : \mathcal{A} \times \mathcal{A} \to Z(\mathcal{A})$ be an even symmetric bilinear map taking values in the center $Z(\mathcal{A})$ of $(\mathcal{A}, \br)$. Suppose that $\mu$ satisfies the identities
$$
\mu([u, v], w) = \mu(\mu(u, v), w) = 0, \quad \text{for all } u, v, w \in \mathcal{A}.
$$
Then define a new product
$$
u \bullet v := [u, v] + \mu(u, v), \quad \forall u,v\in \A
$$
One checks immediately that $(\A, \bullet)$ is a 2-step nilpotent symmetric Leibniz superalgebra. In this case, we call  $(\A, \br + \mu)$ as the 2-step nilpotent symmetric Leibniz superalgebra  associated with the Lie superalgebra $(\A, \br)$ by
means of $\mu$.
\sssbegin{Proposition}\label{cara}
Let $(\mathcal{A}, \bullet, \prs)$  be a pseudo-Euclidean symmetric Leibniz superalgebra. Then $(\mathcal{A}, \bullet, \prs)$ is a flat  quadratic (resp. odd-quadratic) Leibniz superalgebra if and only if $(\mathcal{A}^-, \prs)$ is a quadratic \textup{(}resp. odd-quadratic\textup{)} 2-step nilpotent Lie superalgebra and $(\mathcal{A}^+, \prs)$ is a quadratic \textup{(}resp. odd-quadratic\textup{)} 2-step nilpotent associative superalgebra.
\end{Proposition}

\begin{proof}
Assume that $(\mathcal{A}, \bullet, \prs)$ is a flat  quadratic (resp. odd-quadratic) Leibniz superalgebra. In particular, $(\mathcal{A}, \bullet)$ is a 2-step nilpotent symmetric Leibniz superalgebra. For any $u, v, w \in \mathcal{A}$, we have:
\begin{align*}
\langle [u,v], w\rangle &= \frac{1}{2} \langle u \bullet v - (-1)^{|u||v|} v \bullet u, w\rangle = \frac{1}{2} (\langle u, v \bullet w\rangle - (-1)^{|u||v|} \langle v, u \bullet w \rangle)\\
&= \frac{1}{2} (\langle u, v \bullet w\rangle - (-1)^{|u||v| + |v|(|u|+|w|)} \langle u \bullet w, v\rangle) = \frac{1}{2} (\langle u, v \bullet w\rangle - (-1)^{|v||w|} \langle u, w \bullet v\rangle) \\
&= \langle u, [v, w]\rangle.
\end{align*}
This identity shows that $(\mathcal{A}^-, \prs)$ is a quadratic (resp. odd-quadratic) Lie superalgebra. In a similar manner, one can verify that $(\mathcal{A}^+, \prs)$ is a quadratic (resp. odd-quadratic) associative superalgebra. Since $(\mathcal{A}, \bullet)$ is 2-step nilpotent, we conclude that $(\mathcal{A}^-, \prs)$ is a quadratic (resp. odd-quadratic) 2-step nilpotent Lie superalgebra and $(\mathcal{A}^+, \prs)$ is a quadratic (resp. odd-quadratic) 2-step nilpotent associative superalgebra.

Conversely, suppose that $(\mathcal{A}^-, \prs)$ is a quadratic (resp. odd-quadratic) 2-step nilpotent Lie superalgebra and $(\mathcal{A}^+, \prs)$ is a quadratic (resp. odd-quadratic) 2-step nilpotent associative superalgebra. For any $u, v, w \in \mathcal{A}$, we have
$$
\langle u \bullet v, w\rangle = \langle [u, v] + \{u, v\}, w\rangle = \langle u, [v, w]\rangle + \langle u, \{v, w\}\rangle = \langle u, v \bullet w\rangle,
$$
which shows that $\prs$ is invariant (resp. odd-invariant) with respect to the product $\bullet$, and hence $(\mathcal{A}, \bullet, \prs)$ is a quadratic (resp. odd-quadratic) Leibniz superalgebra. Since both $\mathcal{A}^-$ and $\mathcal{A}^+$ are 2-step nilpotent, and the identities defining symmetry hold (cf. \eqref{idsymmetric}), it follows that $(\mathcal{A}, \bullet)$ is 2-step nilpotent. By Proposition~\ref{pr1}, we conclude that $(\mathcal{A}, \bullet, \prs)$ is flat. Therefore, $(\mathcal{A}, \bullet, \prs)$ is a flat quadratic (resp. odd-quadratic) Leibniz superalgebra.
\end{proof}

\sssbegin{Theorem}\label{chara-2-step}
 Let $(\A, \bullet, \prs)$ be a pseudo-Euclidean 2-step nilpotent symmetric Leibniz superalgebra obtained from a Lie algebra $(\A, \br)$ by means of   $\mu$. Denote by $Z(\A)$ the center of $(\A, \br)$. Then, the bilinear form $\prs$ is invariant \textup{(}resp. odd-invariant\textup{)} if and only if $(\A, \br, \prs)$ is a quadratic \textup{(}resp. odd-quadratic\textup{)} 2-step nilpotent Lie superalgebra, and there exist an even (resp. odd) symmetric trilinear form $T : \A \times \A \times \A \to \K$ and a totally isotropic subspace $J \subseteq Z(\A)$ such that
$$
T(J^\bot, \cdot, \cdot) = 0 \quad \text{and} \quad \langle \mu(u,v), w \rangle = T(u,v,w), \quad \text{for all }u,v,w \in \A.
$$
\end{Theorem}
\begin{proof}
Assume that $\prs$ is invariant (resp. odd-invariant) on $(\A, \bullet)$. According to Proposition \ref{cara} and the structural description of symmetric Leibniz superalgebras mentioned above, we have that $(\mathcal{A}, \br, \prs)$ is a quadratic (resp. odd-quadratic) 2-step nilpotent Lie superalgebra, and the symmetric bilinear map $\mu(u, v)$ satisfies $\mu(u, v) \in Z(\mathcal{A})$ for all $u, v \in \mathcal{A}$, with the relations 
$$
\mu([u, v], w) = \mu(\mu(u, v), w) = 0, \quad \text{and} \quad \langle \mu(u, v), w \rangle = \langle u, \mu(v, w) \rangle.
$$
Let $J := \operatorname{span}\{\mu(u, v) \mid u, v \in \mathcal{A} \} \subseteq Z(\mathcal{A})$. Then, for any $u, v, w, z \in \mathcal{A}$, we have
$$
\langle \mu(u, v), \mu(w, z) \rangle = \langle u, \mu(v, \mu(w, z)) \rangle = 0,
$$
since $\mu(v, \mu(w, z)) = 0$. Hence, $J$ is a totally isotropic subspace of $\mathcal{A}$.
Define a trilinear map $T : \mathcal{A} \times \mathcal{A} \times \mathcal{A} \rightarrow \mathbb{K}$ by
$$
T(u, v, w) := \langle \mu(u, v), w \rangle.
$$
We now show that $T$ is symmetric. Indeed, for all  $u, v, w \in \mathcal{A}$, we have:
\begin{align*}
T(u, v, w) &= \langle \mu(u, v), w \rangle = (-1)^{|u||v|} \langle \mu(v, u), w \rangle = (-1)^{|u||v|} T(v, u, w), \\
T(u, v, w) &= \langle u, \mu(v, w) \rangle = (-1)^{|u|(|v|+|w|)} \langle \mu(v, w), u \rangle = (-1)^{|u|(|v|+|w|)} T(v, w, u).
\end{align*}
Thus, $T$ is a symmetric trilinear form. Moreover, if $\prs$ is even (resp. odd), then it is clear that $T$ is also even (resp. odd).

On the other hand,  for any $u, v \in \mathcal{A}$ and $w \in J^\bot$, we have
$$
T(w, u, v) = \langle \mu(w, u), v \rangle = \langle w, \mu(u, v) \rangle = 0,
$$
so that $T(J^\bot, \cdot, \cdot) = 0$.

Conversely, suppose that $(\mathcal{A},\br, \prs)$ is a quadratic (resp. odd-quadratic) 2-step nilpotent Lie superalgebra, and there exist a even (resp. odd) symmetric trilinear form $T : \mathcal{A} \times \mathcal{A} \times \mathcal{A} \rightarrow \mathbb{K}$ and a totally isotropic subspace $J \subseteq Z(\mathcal{A})$ such that
$$
T(J^\bot, \cdot, \cdot) = 0 \quad \text{and} \quad \langle \mu(u, v), w \rangle = T(u, v, w), \quad \text{for all } u, v, w \in \mathcal{A}.
$$
Then for all $u, v \in \mathcal{A}$ and $w \in J^\bot$, we have
$$
T(w, u, v) = \langle \mu(w, u), v \rangle = \langle w, \mu(u, v) \rangle = 0,
$$
implying that $\mu(u,v) \in J$ , and hence $J \supseteq \operatorname{span}\{\mu(u, v) \mid u, v \in \mathcal{A} \}$. Since $J \subseteq Z(\mathcal{A})$ is totally isotropic, it follows that
$$
\langle \mu(u, v), \mu(w, z) \rangle = \langle \mu(\mu(u, v), w), z \rangle = 0,
$$
and
$$
\langle \mu([u, v], w), z \rangle = \langle [u, v], \mu(w, z) \rangle = \langle u, [v, \mu(w, z)] \rangle = 0.
$$
Thus, the bilinear product defined by $u \bullet v = [u, v] + \mu(u, v)$ endows $\mathcal{A}$ with the structure of a 2-step nilpotent symmetric Leibniz superalgebra. Since $T$ is symmetric trilinear form, then $\mu$ is invariant (resp. odd-invariant). According to Proposition \ref{cara}, the structure $(\mathcal{A}, \bullet, \prs)$ is a quadratic (resp. odd-quadratic) symmetric Leibniz superalgebra.
\end{proof}

\sssbegin{Remark}{\rm 
    The above Theorem provides a characterization of quadratic 2-step nilpotent symmetric Leibniz superalgebras based on quadratic 2-step nilpotent Lie superalgebras. 
    A similar approach is followed in \cite{chin}. More generally, quadratic Lie superalgebras play a significant role in the theory and have been the subject of numerous studies (see \cite{Benayadi6, Benayadi2,Benayadi5}, among others). This framework enables the straightforward construction of quadratic 2-step nilpotent symmetric Leibniz superalgebras, along with explicit examples.}
\end{Remark}

\sssbegin{Example}{\rm 
Let \( \mathcal{A} = \mathcal{A}_{\bar 0} \oplus \mathcal{A}_{\bar 1} \) be a four-dimensional Lie superalgebra, where \( \mathcal{A}_0 = \mathrm{span}\{e_1, e_2\} \) and \( \mathcal{A}_1 = \mathrm{span}\{e_3, e_4\} \). The Lie bracket \( \br \) is defined on the basis elements by:
\[
[e_2, e_4] = -[e_4, e_2] = e_3, \quad [e_4, e_4] = e_1.
\]
Consider the even symmetric bilinear form \( \prs : \mathcal{A} \times \mathcal{A} \to \mathbb{K} \) given by:
\[
\langle e_1, e_2 \rangle = 1, \quad \langle e_3, e_4 \rangle = 1,
\]
with all other pairings zero. A direct computation shows that this bilinear form is invariant under the Lie bracket.

Let us now determine the space \( J \) and the trilinear form \( T \) associated with a symmetric product \( \mu \). Since \( J \) is a totally isotropic subspace of the center \( Z(\mathcal{A}) \), we have \( J \subseteq \mathrm{span}\{e_1, e_3\} \). If we choose $J = \mathrm{span}\{e_3\}$, then $J^\perp = \mathrm{span}\{e_1, e_2, e_3\}$. The condition $T(J^\perp, \cdot, \cdot) = 0$ then implies that
$$
T(e_1, \cdot, \cdot) = T(e_2, \cdot, \cdot) = T(e_3, \cdot, \cdot) = 0.
$$
Moreover, since $e_4 \in \mathcal{A}_1$ and $T$ is symmetric, we also have $T(e_4, \cdot, \cdot) = 0$, which implies that $T = 0$, and consequently $\mu = 0$.

 Whether we choose \( J = \mathrm{span}\{e_1\} \) or \( J = \mathrm{span}\{e_1, e_3\} \), the condition \( T(J^\perp, \cdot, \cdot) = 0 \) implies that \( \mu(u, v) \in J \) for all \( u, v \in \mathcal{A} \). A straightforward computation shows that \( \mu(e_2, e_2) = \alpha e_1,\) where \( \alpha \neq 0 \). Therefore, the quadratic 2-step nilpotent symmetric Leibniz superalgebra constructed from the Lie superalgebra \( (\mathcal{A}, \br) \) by means of \( \mu \) is defined by:
\[
\begin{array}{c}
     e_2 \bullet e_2 = \alpha e_1, \quad e_2 \bullet e_4 = -e_4 \bullet e_2 = e_3, \quad e_4 \bullet e_4 = e_1, \text{ and } \\[1mm]
    \langle e_1, e_2 \rangle = 1, \quad \langle e_3, e_4 \rangle = 1.
\end{array}
\]
}
\end{Example}

\sssbegin{Example}{\rm 
Let $\mathcal{A} = \mathcal{A}_{\bar 0} \oplus \mathcal{A}_{\bar 1}$ be a six-dimensional Lie superalgebra, where $\mathcal{A}_0 = \mathrm{span}\{e_1, e_2\}$ and $\mathcal{A}_1 = \mathrm{span}\{e_3, e_4, e_5, e_6\}$. The non-zero Lie brackets are given by:

$$
[e_2, e_5] = e_3, \quad [e_2, e_6] = e_4, \quad [e_5, e_5] = e_1, \quad [e_6, e_6] = e_1.
$$
Consider the even symmetric bilinear form $\prs : \mathcal{A} \times \mathcal{A} \to \mathbb{K}$ defined by:
$$
\langle e_1, e_2 \rangle = 1, \quad \langle e_3, e_5 \rangle = 1, \quad \langle e_4, e_6 \rangle = 1,
$$
A direct computation shows that this bilinear form is invariant under
the Lie bracket.

We aim to determine the space $J$ and the trilinear form \( T \) associated with a symmetric product \( \mu \). Since $J$ must be a totally isotropic subspace of the center $Z(\mathcal{A})$, it satisfies $J \subseteq \mathrm{span}\{e_1, e_3, e_4\}$. To construct an explicit example, we consider $J = \mathrm{span}\{e_1, e_3, e_4\}$. The condition $T(J^\perp, \cdot, \cdot) = 0$ then implies $\mu(u, v) \in J$ for all $u, v \in \mathcal{A}$.
A direct computation shows that
$$
\mu(e_2, e_2) = \alpha e_1, \quad \mu(e_2, e_5) = \beta e_4, \quad \mu(e_2, e_6) = -\beta e_3, \quad \mu(e_5, e_6) = \beta e_1,
$$
with $(\alpha, \beta) \neq (0, 0)$. Therefore, the corresponding quadratic 2-step nilpotent symmetric Leibniz superalgebra is given by:
$$
\begin{array}{c}
e_2 \bullet e_2 = \alpha e_1, \quad e_2 \bullet e_5 = e_3 + \beta e_4, \quad e_2 \bullet e_6 = e_4 - \beta e_3, \quad e_5 \bullet e_2 = -e_3 + \beta e_4, \\[1mm]
 e_5 \bullet e_5 = e_1, \quad e_5 \bullet e_6 = \beta e_1, \quad e_6 \bullet e_2 = -e_4 - \beta e_3, \quad e_6 \bullet e_5 = \beta e_1, \quad e_6 \bullet e_6 = e_1, \text{ and } \\[1mm]
 \langle e_1, e_2 \rangle = 1, \quad \langle e_3, e_5 \rangle = 1, \quad \langle e_4, e_6 \rangle = 1.
\end{array}
$$
}
\end{Example}
 Next, we introduce an odd-quadratic 2-step nilpotent symmetric Leibniz superalgebra.

\sssbegin{Example}
{\rm 
Let $\mathcal{A} = \mathcal{A}_{\bar 0} \oplus \mathcal{A}_{\bar 1}$ be a four-dimensional Lie superalgebra, where $\mathcal{A}_0 = \mathrm{span}\{e_1, e_2 \}$ and $\mathcal{A}_1 = \mathrm{span}\{f_1, f_2\}$. The non-zero Lie brackets are given by:
$$
[f_1, f_1] = e_1.
$$
Consider the odd symmetric bilinear form $\langle\cdot,\cdot\rangle : \mathcal{A} \times \mathcal{A} \to \mathbb{K}$ defined by:
$$
\langle e_1, f_1\rangle = 1, \quad \langle e_2, f_2\rangle = 1.
$$
A straightforward computation shows that this bilinear form is odd-invariant under the Lie bracket.

We now determine the space $J$ and the trilinear map $T$ associated with a symmetric product $\mu$. Since $J$ must be a totally isotropic subspace of the center $Z(\mathcal{A}) = \mathrm{span}\{e_1, e_2, f_2\}$, the possible choices for $J$ are:
$$
J = \mathrm{span}\{e_1\},\quad \mathrm{span}\{e_2\},\quad \mathrm{span}\{f_2\},\quad \mathrm{span}\{e_1, f_2\},\quad \text{or}\quad \mathrm{span}\{e_1, e_2\}.
$$
We choose $J = \mathrm{span}\{e_1, f_2\}$, which satisfies $J = J^\perp$. Therefore, the condition $T(J, \cdot, \cdot) = 0$ implies
$
T(e_1, \cdot, \cdot) =  T(f_2, \cdot, \cdot) = 0.
$
Since $T$ is an odd symmetric trilinear map, then we have
$
T(e_2, e_2, e_2) = T(f_1, f_1, \cdot) = 0.
$
The only non-zero component we may have is $T(e_2, e_2, f_1) = \alpha$, where $\alpha \in \mathbb{K} \setminus \{0\}$. Consequently, the associated symmetric product is defined by:
$$
\mu(e_2, e_2) = \alpha e_1, \quad \mu(e_2, f_1) = \mu(f_1, e_2) = \alpha f_2.
$$
Therefore, the corresponding odd-quadratic 2-step nilpotent symmetric Leibniz superalgebra is given by:
\[
\begin{array}{c}
 e_2 \bullet e_2 = \alpha e_1, \quad e_2 \bullet f_1 = f_1 \bullet e_2 = \alpha f_2, \quad f_1 \bullet f_1 = e_1, \text{ and } \\[1mm]
 \langle e_1, f_1\rangle = 1, \quad \langle e_2, f_2\rangle = 1.
\end{array}
\]
}
\end{Example}

\sssbegin{Example}{\rm 
Let $\mathcal{A} = \mathcal{A}_{\bar 0} \oplus \mathcal{A}_{\bar 1}$ be a six-dimensional Lie superalgebra, where $\mathcal{A}_0 = \mathrm{span}\{e_1, e_2, e_3\}$ and $\mathcal{A}_1 = \mathrm{span}\{f_1, f_2,f_3\}$. The non-zero Lie brackets are given by:
\[
[e_1, e_2] = e_3, \quad [e_2, f_3] = f_1,\quad [e_1,f_3]=-f_2.
\]
Consider the odd symmetric bilinear form $\prs : \mathcal{A} \times \mathcal{A} \to \mathbb{K}$ defined by:
\[
\langle e_1, f_1\rangle =1, \quad \langle e_2, f_2\rangle=1,\quad \langle e_3, f_3\rangle =1.
\]
A direct computation shows that this bilinear form is odd-invariant under the Lie bracket $\br$.

We aim to determine the space $J$ and the trilinear form \( T \) associated with a symmetric product \( \mu \). Since $J$ must be a totally isotropic subspace of the center $Z(\mathcal{A})=\{e_3, f_1,f_2\}$, it satisfies $J \subseteq \mathrm{span}\{e_3, f_1,f_2\}$. The possible choices for $J$ are:
\[
J = \mathrm{span}\{e_3\},\;  \mathrm{span}\{f_1\},\; \mathrm{span}\{f_2\},\; \mathrm{span}\{e_3, f_1\},\; \mathrm{span}\{e_3, f_2\},\; \mathrm{span}\{f_1, f_2\},\; \text{or}\; \mathrm{span}\{e_3, f_1, f_2\}.
\]
We choose $J = \mathrm{span}\{e_3,f_1, f_2\}$, which satisfies $J = J^\perp$. Therefore, the condition $T(J, \cdot, \cdot)=0$ implies
\[
T(e_3, \cdot, \cdot)=T(f_1, \cdot, \cdot)=T(f_2, \cdot, \cdot)=0.
\]
Since $T$ is odd symmetric, we also have $T(e_i, e_j,e_k)=0$ for all $i,j,k\in \{1,2\}$ and $T(f_3,f_3,\cdot)=0$. Thus, we get:
\[
T(e_1,e_1,f_3)=\alpha,\quad T(e_1,e_2,f_3)=\beta,\quad T(e_2,e_2,f_3)=\lambda,
\]
where $(\alpha, \beta, \lambda)\neq(0,0,0)$. Consequently, the symmetric product $\mu$ is given by:
\[
\mu(e_1,e_1)=\alpha e_3,\;  \mu(e_1,e_2)=\beta e_3,\; \mu(e_2,e_2)=\lambda e_3,\; \mu(e_1,f_3)=\alpha f_1+\beta f_2,\;  \mu(e_2,f_3)=\beta f_1+\lambda f_2.
\]
Therefore, the corresponding odd-quadratic 2-step nilpotent symmetric Leibniz superalgebra is:
\begin{align*}
    & e_1\bullet e_1= \alpha e_3,\quad e_1\bullet e_2= (1+\beta) e_3,\quad e_2\bullet e_1= (\beta-1) e_3,\quad e_2\bullet e_2= \lambda e_3, \\
    & e_1\bullet f_3=(1+\alpha)f_1+\beta f_2,\quad f_3\bullet e_1=(\alpha-1)f_1+\beta f_2, \quad e_2\bullet f_3=\beta f_1+(\lambda -1) f_2,\\
    &f_3\bullet e_2=\beta f_1+(\lambda +1) f_2,\end{align*}
    and
    \[
    \langle e_1, f_1\rangle =1,\quad \langle e_2, f_2\rangle=1,\quad \langle e_3, f_3\rangle =1.
\]
}
\end{Example}

\subsection{$T^*$-extension and $\Pi(T^*)$-extension of flat quadratic Leibniz superalgebras}
This part deals with the $T^*$-extension in the context of even quadratic Leibniz superalgebras and $\Pi(T^*)$-extensions in the case of odd quadratic Leibniz superalgebras. We provide a description of  quadratic 2-step nilpotent symmetric  Leibniz superalgebras. Recall that the notion of the $T^*$-extension was introduced by Bordemann in \cite{Bordemann} in the context of nondegenerate associative bilinear forms on non-associative algebras. See also \cite{BC, BM} on Lie (super)algebras when the bilinear form is antisymmetric, and referred to as Lagrangian extensions.
\sssbegin{Theorem}[\cite{BenFa0}]\label{theo TT-extension}
Let \( (\A, \bullet) \) be a symmetric Leibniz superalgebra and let \( \Omega : \A \times \A \to \A^* \) be an even bilinear map. Then, the \(\mathbb{Z}_2\)-graded vector space \( T_\Omega^*(\A) := \A \oplus \A^* \), endowed with the product:
\[
(u + f)\bullet (v + g) := u\bullet v + \Omega(u,v) + (\Ll^\bullet)^*(u)(g)+(-1)^{|u||v|}(\Rr^\bullet)^*(v)(f) ,
\]
for all \( u + f \in T_\Omega^*(\A)_{|u|} \), \( v + g \in T_\Omega^*(\A)_{|v|} \),
is a symmetric Leibniz superalgebra if and only if \( \Omega \) is an even Leibniz 2-cocycle of \((\A, \bullet) \) relative to the representation \( ((\Ll^\bullet)^*, (\Rr^\bullet)^*) \).
Moreover, the even bilinear form
\[
B : T_\Omega^*(\A) \times T_\Omega^*(\A) \to \mathbb{K}, \quad (u + f,\, v + g) \mapsto f(v) + (-1)^{|u||v|} g(u)
\]
is symmetric and non-degenerate on \( T_\Omega^*(\A) \). The bilinear form \( B \) is invariant on \( T_\Omega^*(\A) \) if and only if
\begin{equation}
\Omega(u,v)(w) = (-1)^{|u|(|v| + |w|)} \Omega(v, w)(u), \quad \forall\, u,v,w\in \A.
\label{cyclic} \end{equation}
In this case, the quadratic Leibniz superalgebra \( (T_\Omega^*(\A), B) \) is called the \emph{\( T^*_\Omega \)-extension of \(( \A, \bullet) \) by means of \( \Omega \)}.
\end{Theorem}

Let $(\mathcal{A}, \bullet)$ be a left symmetric superalgebra, and let $(\Ll^\bullet, \Rr^\bullet)$ denote its left and right adjoint representations, respectively. According to Proposition~\ref{represantaion}, the pair $(\Ll^{\bullet *}, \Rr^{\bullet *})$ defines a representation of $(\mathcal{A}, \bullet)$ on $\mathcal{A}^*$.

We now define another pair of linear maps:

$$
(\Rr^{\bullet *})^\Pi: \mathcal{A} \to \mathrm{End}(\Pi(\mathcal{A}^*)), \qquad
(\Ll^{\bullet *})^\Pi: \mathcal{A} \to \mathrm{End}(\Pi(\mathcal{A}^*)),
$$
given for homogeneous $u \in \mathcal{A}$ by:
$$
(\Rr^{\bullet *})^\Pi(u) =  \Pi \circ \Rr^{\bullet *}(u) \circ \Pi, \qquad
(\Ll^{\bullet *})^\Pi(u) =  \Pi \circ \Ll^{\bullet *}(u) \circ \Pi.
$$

\sssbegin{Proposition}\label{Pirepresantaion}
Let $(\mathcal{A}, \bullet)$ be a symmetric Leibniz superalgebra, and let $(\Rr^\bullet, \Ll^\bullet)$ be its regular representation in $\mathcal{A}$. Then the pair $((\Rr^{\bullet *})^\Pi, (\Ll^{\bullet *})^\Pi)$ defines a representation of $(\mathcal{A}, \bullet)$ on $\Pi(\mathcal{A}^*)$.
\end{Proposition}

\begin{proof}
Let $u, v \in \mathcal{A}$. Using the definition of $(\Ll^{\bullet *})^\Pi$, we have
\begin{align*}
(\Ll^{\bullet *})^\Pi(u \bullet v)
&=  \Pi \circ \Ll^{\bullet *}(u \bullet v) \circ \Pi \\
&=  \Pi \circ \left( \Ll^{\bullet *}(u) \circ \Ll^{\bullet *}(v) - (-1)^{|u||v|} \Ll^{\bullet *}(v) \circ \Ll^{\bullet *}(u) \right) \circ \Pi \\
&= \Pi \circ \Ll^{\bullet *}(u) \circ \Pi \circ \Pi \circ \Ll^{\bullet *}(v) \circ \Pi-(-1)^{|u||v|} \Pi \circ \Ll^{\bullet *}(v) \circ \Pi \circ \Pi \circ \Ll^{\bullet *}(u) \circ \Pi \\
  &= (\Ll^{\bullet *})^\Pi(u) \circ (\Ll^{\bullet *})^\Pi(v) - (-1)^{|u||v|} (\Ll^{\bullet *})^\Pi(v) \circ (\Ll^{\bullet *})^\Pi(u) \\
  &= [(\Ll^{\bullet *})^\Pi(u), (\Ll^{\bullet *})^\Pi(v)].
  \end{align*}
The remaining conditions can be established in a similar manner. Hence, $((\Rr^{\bullet *})^\Pi, (\Ll^{\bullet *})^\Pi)$ defines a representation of $(\mathcal{A}, \bullet)$ on $\Pi(\mathcal{A}^*)$.
  \end{proof}
Now, we introduce an analogue of the notion of $T^*_\Omega$-extension in the odd-quadratic case, which we denote by $\Pi(T^*_\Omega)$.
\sssbegin{Theorem}\label{theo Pi-extension}
Let \( (\A, \bullet) \) be a symmetric Leibniz superalgebra and let \( \Omega : \A \times \A \to \Pi(\A^*) \) be an even bilinear map. Then, the \(\mathbb{Z}_2\)-graded vector space \( \Pi(T_\Omega^*)(\A) := \A \oplus \Pi(\A^*) \), endowed with the product:
\[
(u + \Pi(f))\bullet (v + \Pi(g)) :=u\bullet v +\Omega(u,v) +  (\Ll^{\bullet *})^\Pi(u) (\Pi(g))+ (-1)^{|u||v|}(\Rr^{\bullet *})^\Pi(v) (\Pi(f)),
\]
for all \( u + \Pi(f) \in \Pi(T_\Omega^*)(\A)_{|u|} \), \( v + \Pi(g) \in\Pi(T_\Omega^*)(\A)_{|v|} \),
is a symmetric Leibniz superalgebra if and only if \( \Omega \) is an even Leibniz 2-cocycle of \((\A, \bullet) \) relative to the representation \( ((\Rr^{\bullet *})^\Pi, (\Ll^{\bullet *})^\Pi) \).  Moreover, the even bilinear form
\[
B : \Pi(T_\Omega^*)(\A) \times \Pi(T_\Omega^*)(\A)\to \mathbb{K}, \quad (u + \Pi(f),\, v + \Pi(g)) \mapsto f(v) + (-1)^{|u||v|} g(u)
\]

is symmetric and non-degenerate on \( \Pi(T_\Omega^*)(\A) \). The bilinear form \( B \) is odd-invariant on \( \Pi(T_\Omega^*)(\A) \) if and only if
\begin{equation}
\Pi(\Omega(u, v))(w) = (-1)^{|u|(|v| + |w|)} \Pi(\Omega(v, w))(u)
, \quad \forall\, u,v,w\in \A.
\label{cyclic1}\end{equation}
In this case, the odd-quadratic Leibniz superalgebra \( (\Pi(T_\Omega^*)(\A), B) \) is called the \emph{\( \Pi(T^*_\Omega) \)-extension of \(( \A, \bullet) \) by means of \( \Omega \)}.
\end{Theorem}

\begin{proof}
Let $u + \Pi(f), v + \Pi(g) \in \Pi(T_\Omega^*)(\mathcal{A})$, and let $w \in \mathcal{A}$. By Proposition~\ref{Pirepresantaion}, the pair $((\Rr^{\bullet *})^\Pi, (\Ll^{\bullet *})^\Pi)$ defines a representation of $(\mathcal{A}, \bullet)$ on $\Pi(\mathcal{A}^*)$. Therefore, the terms $u \bullet_\Omega \Pi(g)$ and $\Pi(f) \bullet_\Omega v$ make sense in the Leibniz product on $\Pi(T_\Omega^*)(\mathcal{A})$.

Now, have
\begin{align*}
&\big((u \bullet_\Om v) \bullet_\Om w\big) - \big(u \bullet_\Om (v \bullet_\Om w)\big) + (-1)^{|u||v|} \big(v \bullet_\Om (u \bullet_\Om w)\big) \\
&= \big(u \bullet v + \Omega(u, v)\big) \bullet_\Om w - u \bullet_\Om \big(v \bullet w + \Omega(v, w)\big) + (-1)^{|u||v|} v \bullet_\Om \big(u \bullet w + \Omega(u, w)\big) \\
&= (u \bullet v) \bullet w + \Omega(u \bullet v, w) + (-1)^{|u||v| + |w|(|u|+|v|)} (\Rr^{\bullet *})^\Pi(w)(\Omega(u, v)) \\
&\quad - u \bullet (v \bullet w) - \Omega(u, v \bullet w) - (\Ll^{\bullet *})^\Pi(u)(\Omega(v, w)) \\
&\quad + (-1)^{|u||v|} \big( v \bullet (u \bullet w) + \Omega(v, u \bullet w) + (\Ll^{\bullet *})^\Pi(v)(\Omega(u, w)) \big)\\ &=\big((u \bullet v) \bullet w - u \bullet (v \bullet w) + (-1)^{|u||v|} v \bullet (u \bullet w)\big) \\
&\quad + \Omega(u \bullet v, w) - \Omega(u, v \bullet w) + (-1)^{|u||v|} \Omega(v, u \bullet w) \\
&\quad + (-1)^{|w|(|u|+|v|)} (\Rr^{\bullet *})^\Pi(w)(\Omega(u, v)) - (\Ll^{\bullet *})^\Pi(u)(\Omega(v, w)) + (-1)^{|u||v|} (\Ll^{\bullet *})^\Pi(v)(\Omega(u, w))
\end{align*}
Since $(\mathcal{A}, \bullet)$ is a symmetric Leibniz superalgebra, the first line vanishes. Hence, the product on $\Pi(T_\Omega^*)(\mathcal{A})$ is left Leibniz identity if and only if
$$
\begin{aligned}
&\Omega(u \bullet v, w) - \Omega(u, v \bullet w) + (-1)^{|u||v|} \Omega(v, u \bullet w) \\
&\quad - (\Ll^{\bullet *})^\Pi(u)(\Omega(v, w)) + (-1)^{|u||v|} (\Ll^{\bullet *})^\Pi(v)(\Omega(u, w)) \\
&\quad + (-1)^{|w|(|u| + |v|)} (\Rr^{\bullet *})^\Pi(w)(\Omega(u, v)) = 0.
\end{aligned}
$$
This is precisely the cocycle condition for if even 2-cocycle of the left Leibniz superalgebra $(\A, \bullet)$ relative to \( ((\Rr^{\bullet *})^\Pi, (\Ll^{\bullet *})^\Pi) \). The same reasoning applies to the right Leibniz identity. Therefore, $\Pi(T_\Omega^*)(\mathcal{A})$ is symmetric Leibniz if and only if $\Omega$ is an even Leibniz 2-cocycle relative to \( ((\Rr^{\bullet *})^\Pi, (\Ll^{\bullet *})^\Pi) \).

Now, we show that $B$ is an odd-quadratic form on $(\mathcal{A}, \bullet)$. By definition of $B$, it is clear that $B$ is a  non-degenerate, symmetric and odd bilinear form. On the one hand, we have 
\begin{align*}
   & B( (u + \Pi(f)) \bullet (v + \Pi(g)), w + \Pi(h))\\
    &= B(u\bullet v+\Om(u,v)+ (\Ll^{\bullet *})^\Pi(u) (\Pi(g))+ (-1)^{|u||v|}(\Rr^{\bullet *})^\Pi(v) (\Pi(f)), w + \Pi(h))\\&=\Pi(\Om(u,v))(w)+(f\circ \Ll^\bullet_v)(w)+(-1)^{|u||v|} (g \circ \Rr^\bullet_u)(w)+(-1)^{|w|(|u|+|v|)} h(u\bullet v)\\&= \Pi(\Om(u,v))(w)+f(v\bullet w)+ (-1)^{|u|(|v|+|w|)}g (w\bullet u)+(-1)^{|w|(|u|+|v|)} h(u\bullet v). 
\end{align*}
On the other hand, we have
\begin{align*}
  &  B( u + \Pi(f),  (v + \Pi(g))\bullet( w + \Pi(h)))\\&= B( u + \Pi(f), v \bullet w +\Omega(v, w) +(\Ll^{\bullet *})^\Pi(v) (\Pi(h))+ (-1)^{|v||w|}(\Rr^{\bullet *})^\Pi(w) (\Pi(g)) )\\&=f(v\bullet w)+(-1)^{|u|(|v|+|w|)}(\Pi(\Om(v,w))(u)+(g\circ \Ll^\bullet_w)(u)+ (-1)^{|v||w|}h \circ \Rr^\bullet_v(u))\\&= f(v\bullet w)+(-1)^{|u|(|v|+|w|)}\Pi(\Om(v,w))(u)+(-1)^{|u|(|v|+|w|)}g(w\bullet u)+ (-1)^{|w|(|u|+|v|)}h (u\bullet v). 
\end{align*}

Therefore, we conclude that
$$
B\big( (u + \Pi(f)) \bullet (v + \Pi(g)), w + \Pi(h)\big) = B\big( u + \Pi(f), (v + \Pi(g)) \bullet (w + \Pi(h))\big)
$$
if and only if
$$
\Pi(\Omega(u, v))(w) = (-1)^{|u|(|v| + |w|)} \Pi(\Omega(v, w))(u).\qed
$$
\noqed
\end{proof}

\sssbegin{Definition}
A quadratic (resp. odd-quadratic) Leibniz superalgebra $(\mathcal{A}, \bullet, \prs)$ is said to be reduced if $\mathcal{A} \neq \{0\}$ and $\operatorname{Ann}(\mathcal{A}, \bullet)$ is totally isotropic with respect to  $\prs$.
\end{Definition}

\sssbegin{Lemma}[\cite{BenFa}]\label{orth}
Let $(\mathcal{A}, \bullet)$ be a quadratic (resp. odd-quadratic) symmetric Leibniz superalgebra. Then
$$
(\mathcal{A} \bullet \mathcal{A})^\perp = \operatorname{Ann}(\mathcal{A}, \bullet),
$$
\end{Lemma}

Now, we provide the converse of Theorem \ref{theo TT-extension} in the particular case of quadratic 2-step nilpotent symmetric Leibniz superalgebras.
\sssbegin{Theorem}\label{T-ex}
Let $(\mathcal{A}, \bullet, \prs)$ be a reduced quadratic 2-step nilpotent symmetric Leibniz superalgebra. Then, $(\mathcal{A}, \bullet)$ is a $T^*_\Omega$-extension of a trivial algebra $\mathfrak{h}$ by means of $\Omega$ where $\Om$ is non-degenerate.
\end{Theorem}
\begin{proof}
Since $(\mathcal{A}, \bullet)$ is 2-step nilpotent, we have
$
\mathcal{A} \bullet \mathcal{A} \subseteq \operatorname{Ann}(\mathcal{A}, \bullet).
$
According to Lemma \ref{orth}, we also have
$
(\mathcal{A} \bullet \mathcal{A})^\perp = \operatorname{Ann}(\mathcal{A}, \bullet).
$
Given that $\operatorname{Ann}(\mathcal{A}, \bullet)$ is totally isotropic with respect to  $\prs$, it follows that
$
\mathcal{A} \bullet \mathcal{A} = \operatorname{Ann}(\mathcal{A}, \bullet),
$ and consequently $\dim(\mathcal{A})$ is even. Thus, $\operatorname{Ann}(\mathcal{A}, \bullet)$ is a totally isotropic subspace of dimension $\frac{1}{2} \dim(\mathcal{A})$. We can therefore decompose $\mathcal{A}$ as
$
\mathcal{A} = \operatorname{Ann}(\mathcal{A}, \bullet) \oplus \mathfrak{h},
$
where $\mathfrak{h}$ is a totally isotropic subspace of the same dimension. Define the map
$
\varphi: \operatorname{Ann}(\mathcal{A}, \bullet) \to \mathfrak{h}^*, \, \varphi(a)(u) := \langle a, u \rangle,
$ for all $a\in \operatorname{Ann}(\mathcal{A}, \bullet)$, $u\in \mathfrak{h}$,
which is clearly a linear isomorphism. Since the bilinear form $\prs$ is even, it follows that $\varphi$ is an even isomorphism.  Now, define the bilinear map $\Omega: \mathfrak{h} \times \mathfrak{h} \to \mathfrak{h}^*$ by
$
\Omega(u, v) := \varphi(u \bullet v).
$ Since both $\varphi$ and the product $\bullet$ are even maps, it follows that $\Omega$ is also even. Since \( \varphi \) is bijective, it follows that \( \Omega \) is also bijective.

On the other hand, for all $u, v, w \in \mathfrak{h}$, we have
\begin{align*}
\Omega(u, v)(w) = \langle u \bullet v, w \rangle 
= \langle u, v \bullet w \rangle 
= (-1)^{|u|(|v|+|w|)} \langle v \bullet w, u \rangle = (-1)^{|u|(|v|+|w|)} \Omega(v, w)(u).
\end{align*}
Now, consider $\mathfrak{h}$ as a trivial superalgebra (i.e., the product on $\mathfrak{h}$ is zero). Since $\Omega$ is an even bilinear map satisfying the cyclicity condition given by \eqref{cyclic}, it follows from Theorem \ref{theo TT-extension} that $T^*_\Om(\mathfrak{h})=\mathfrak{h} \oplus \mathfrak{h}^*$ is the $T_\Omega^*$-extension of the trivial superalgebra $\mathfrak{h}$ by means of $\Omega$.

Moreover, define the even map,
$
\Phi: \mathfrak{h} \oplus \operatorname{Ann}(\mathcal{A}, \bullet) \to T^*_\Omega(\mathfrak{h}),$ by $  \Phi(u + a) := u + \varphi(a).
$
This map is clearly linear and bijective. For all $u, v \in \mathfrak{h}$ and $a, b \in \operatorname{Ann}(\mathcal{A}, \bullet)$, we have:
\begin{align*}
\Phi((u + a) \bullet (v + b)) &= \Phi(u \bullet v) = \varphi(u \bullet v) = \Omega(u, v),\\
&= (u + \varphi(a)) \bullet_\Omega (v + \varphi(b)) \\
&= \Phi(u + a) \bullet_\Omega \Phi(v + b),
\end{align*}
so $\Phi$ is an even isomorphism of Leibniz superalgebras. Let $B$ be the even bilinear form on $T^*_\Omega(\mathfrak{h})$. Then for all $u+a \in (\mathfrak{h}\oplus \operatorname{Ann}(\mathcal{A}, \bullet))_{|u|}$ and $v+b \in (\mathfrak{h}\oplus \operatorname{Ann}(\mathcal{A}, \bullet))_{|v|}$ , we have:
\begin{align*}
B(\Phi(u + a), \Phi(v + b)) &= B(u + \varphi(a), v + \varphi(b)) \\
&= \varphi(a)(v) + (-1)^{|u||v|} \varphi(b)(u) \\
&= \langle a, v \rangle + (-1)^{|u||v|} \langle b, u \rangle \\
&= \langle u + a, v + b \rangle,
\end{align*}
so $\Phi$ is an isometric. Therefore, $(\mathcal{A}, \bullet, \prs)$ is isometric and isomorphic to the $T^*_\Omega$-extension of the trivial algebra $\mathfrak{h}$ by means of $\Omega$.
\end{proof}

Let $E = E_{\bar{0}} \oplus E_{\bar{1}}$ be a $\mathbb{Z}_2$-graded vector space. Let
$
\Omega : E \times E \to E^*
$ be an even bilinear map. We say that $\Omega$ is cyclic if, for all  $u,v,w \in E$, the following identity holds
\begin{align}
\Omega(u, v)(w) = (-1)^{|u|(|v|+|w|)} \Omega(v, w)(u).
\label{cycl}
\end{align}
We define the associated even trilinear form $\tilde{\Omega} : E \times E \times E \to \mathbb{K}$ by
$$
\tilde{\Omega}(u,v,w) := \Omega(u,v)(w), \quad \forall u,v,w \in E.
$$
Due to the cyclicity of $\Omega$, the form $\tilde{\Omega}$ satisfies the following identities:
\begin{align}
\tilde{\Omega}(u,v,w) = (-1)^{|u|(|v|+|w|)} \tilde{\Omega}(v,w,u) = (-1)^{|w|(|u|+|v|)} \tilde{\Omega}(w,u,v).
\label{tri}
\end{align}
We denote by $\mathcal{T}^3(E)$ the set of all even trilinear forms $\tilde{\Omega} : E^{\otimes 3} \to \mathbb{K}$ that satisfy identity \eqref{tri}. That is,
$$
\mathcal{T}^3(E) := \left\{ \tilde{\Omega} \in \mathrm{Hom}(E^{\otimes 3}, \mathbb{K})_{\bar{0}} \ \middle| \ \tilde{\Omega} \text{ satisfies \eqref{tri}} \right\}.
$$

Consequently, we have a restricted definition  for the reduced quadratic 2-step nilpotent symmetric Leibniz case as follows:

\sssbegin{Definition}
Let $\mathfrak{h}$ be a $\mathbb{Z}_2$-graded vector space and let $\Om : \mathfrak{h} \times \mathfrak{h} \to \mathfrak{h}^*$ be a non-degenerate cyclic even bilinear map. Define the $\mathbb{Z}_2$-graded vector space $\A := \mathfrak{h} \oplus \mathfrak{h}^*$ and endow it with the product
$$
(u + f)\bullet ( v+ g) := \Om(u, v),
$$
and the symmetric bilinear form
$$
B(u + f, v + g) := f(v) + (-1)^{|u||v|}g(u),
$$
for all $u+f \in \A_{|u|}$ and $v+g \in \A_{|v|}$. Then $(\A, B)$ is called a quadratic  2-step nilpotent symmetric Leibniz  algebra. In this case, we say that $\A$ is the $T_\Om^*$-extension of $\mathfrak{h}$ by $\Om$.
\end{Definition}

\sssbegin{Theorem}
Let $(\A, \bullet_{\Om_1})$ and $(\mathcal{B}, \bullet_{\Om_2})$ be $T_\Om^*$-extensions of $\mathfrak{h}$ by   $\Om_1$ and $\Om_2$, respectively. Then:
\begin{enumerate}
    \item[\textup{(}i\textup{)}]  There exists an even isomorphism between $(\A, \bullet_{\Om_1})$ and $(\mathcal{B}, \bullet_{\Om_2})$ as symmetric Leibniz superalgebras if and only if there exist an even isomorphism $\varphi_1 : \mathfrak{h} \to \mathfrak{h}$ and an even isomorphism $\varphi_2 : \mathfrak{h}^* \to \mathfrak{h}^*$ such that
$$
   \varphi_2(\Om_1(u,v)) = \Om_2(\varphi_1(u), \varphi_1(v)), \quad \forall u,v \in \mathfrak{h}.
   $$
\item[\textup{(}ii\textup{)}]  There exists an isometrically isomorphic (i.e., an isometric isomorphism preserving the natural bilinear form) between $\A$ and $\mathcal{B}$ if and only if there exists an isomorphism $\varphi_1 : \mathfrak{h} \to \mathfrak{h}$ such that
$$
   \Om_1(u,v) = \Om_2(\varphi_1(u), \varphi_1(v)) \circ \varphi_1, \quad \forall u,v \in \mathfrak{h}.
   $$
   \end{enumerate}
\end{Theorem}
\begin{proof}
(i) Let $\varphi: \A \to \mathcal{B}$ be an even isomorphism symmetric Leibniz algebra .  Since the annihilators of both $\A$ and $\mathcal{B}$ satisfy $\operatorname{Ann}(\A, \bullet_{\Om_1}) = \operatorname{Ann}(\mathcal{B}, \bullet_{\Om_2})  = \mathfrak{h}^*$, and these are stable under $\varphi$, there exist even linear maps
$
\varphi_1 : \mathfrak{h} \to \mathfrak{h}, \quad \varphi_1' : \mathfrak{h} \to \mathfrak{h}^*, \quad \text{and} \quad \varphi_2 : \mathfrak{h}^* \to \mathfrak{h}^*
$
such that for all $u \in \mathfrak{h}$, $f \in \mathfrak{h}^*_{|u|}$,
$
\varphi(u + f) = \varphi_1(u) + \varphi_1'(u) + \varphi_2(f).
$ 
It is immediate that $\varphi_2$ is an even linear isomorphism of $\mathfrak{h}^*$. We now show that $\varphi_1$ is also an even isomorphism of $\mathfrak{h}$. Suppose there exists $u \in \mathfrak{h}$ such that $\varphi_1(u) = 0$. Then for any $Y \in \mathcal{B}$,
$
0=\varphi(u)\bullet_{\Om_2} Y = \varphi(u\bullet_{\Om_1} \varphi^{-1}(Y) ).
$
Since $\varphi$ is an  isomorphism, this implies that
$
u\bullet \varphi^{-1}(Y) = 0 \quad \text{for all } Y \in \mathcal{B}.
$
Thus, $u\bullet v = 0$, for any $v\in \A$ and hence  $u \in \operatorname{Ann}(\A, \bullet)  = \mathfrak{h}^*$, which contradicts $u \in \mathfrak{h}$. Therefore, $u = 0$, and $\varphi_1$ is injective. As $\dim \mathfrak{h} = \dim \varphi_1(\mathfrak{h})$, it follows that $A_1$ is an even isomorphism.

Now, for all $u,v \in \mathfrak{h}$ and $f \in \mathfrak{h}^*_{|u|}$,  $g \in \mathfrak{h}^*_{|v|}$, we have
\begin{align*}
&\varphi((u+f)\bullet_{\Om_1} (v+g) )= \varphi(\Om_1(u,v)) = \varphi_2(\Om_1(u,v)),
\esp\\&
\varphi(u+f)\bullet_{\Om_2} \varphi(v+g) = (\varphi_1(u) + \varphi_2(f))\bullet_{\Om_2} (\varphi_1(v) + \varphi_2(g)) = \Om_2(\varphi_1(u), \varphi_1(v)).
\end{align*}
Hence,
$
\varphi_2(\Om_1(u,v)) = \Om_2(\varphi_1(u), \varphi_1(v)),$ for all $ u,v\in \mathfrak{h}.
$

Conversely, assume that there exist even isomorphisms $\varphi_1 : \mathfrak{h} \to \mathfrak{h}$ and $\varphi_2 : \mathfrak{h}^* \to \mathfrak{h}^*$ such that
$
\varphi_2(\Om_1(u,v)) = \Om_2(\varphi_1(u), \varphi_1(v)), $ for all $ u,v \in \mathfrak{h}.
$ Define $\varphi : \A \to \mathcal{B}$ by $
\varphi(u + f) = \varphi_1(u) + \varphi_2(f),$ for all $ u \in \mathfrak{h},\ f \in \mathfrak{h}_{|u|}^*.
$
It is straightforward to verify that $\varphi$ is an even isomorphism symmetric Leibniz algebra .

(ii) Suppose $\varphi: \A \to \mathcal{B}$ is an isometry. Then, as in (i), there exist even maps $\varphi_1 : \mathfrak{h} \to \mathfrak{h}$ and $\varphi_2 : \mathfrak{h}^* \to \mathfrak{h}^*$ such that
$
\varphi(u + f) = \varphi_1(u) + \varphi_2(f),$ for all $ u \in \mathfrak{h},\ f \in \mathfrak{h}^*_\al.
$
For all $u \in \mathfrak{h},\ f \in \mathfrak{h}^*_\al$, we have
$$
B'(\varphi(u), \varphi(f)) = B(u, f) \Rightarrow (-1)^{\al|u|}\varphi_2(f)(\varphi_1(u)) = (-1)^{\al|u|} f(u).
$$
Therefore,
$
\varphi_2(f) = f \circ \varphi_1^{-1},$ $f \in \mathfrak{h}^*_\al.
$
Now, since
$
\varphi_2(\Om_1(u, v)) = \Om_2(\varphi_1(u), \varphi_1(v)),
$
we obtain
$$
\Om_1(u,v) = \Om_2(\varphi_1(u), \varphi_1(v)) \circ \varphi_1,$$ for all $u, v \in \mathfrak{h}.
$

Conversely, suppose there exists an even isomorphism $\varphi_1 : \mathfrak{h} \to \mathfrak{h}$ such that
$$
\Om_1(u,v) = \Om_2(\varphi_1(u), \varphi_1(v)) \circ \varphi_1,$$ for all $u,v \in \mathfrak{h}.
$
Define $\varphi: \A \to \mathcal{B}$ by
$
\varphi(u+ f) = \varphi_1(u) + f \circ \varphi_1^{-1},$ for all $u \in \mathfrak{h},\ f \in (\mathfrak{h}^*)_\al.$ Then, it is clear that $\varphi$ is an even isomorphism of symmetric Leibniz algebras.

Moreover, for all $u+f\in (\mathfrak{h}\oplus \mathfrak{h}^*)_{|u|}$ and $v+g\in (\mathfrak{h}\oplus \mathfrak{h}^*)_{|v|}$ 
\begin{align*}
  B'(\varphi(u+f), \varphi(v+g))&=B'(\varphi_1(u) + f \circ \varphi_1^{-1}, \varphi_1(v) + g \circ \varphi_1^{-1} )\\&= f \circ \varphi_1^{-1}(\varphi(v))+(-1)^{|u||v|} g\circ \varphi_1^{-1}(\varphi(u))\\&= f(u)+(-1)^{|u||v|} g(u)\\&= B(u+f, v+g).
\end{align*}
Therefore, $\varphi$ is an  isometry.
\end{proof}

Let $\mathfrak{h}$ be a $\mathbb{Z}_2$-graded vector space and let $\Om : \mathfrak{h} \times \mathfrak{h} \to \mathfrak{h}^*$ be a non-degenerate cyclic even bilinear map. We denote by $\tilde{\Omega}$ the associated even trilinear form defined by
$$
\tilde{\Omega}(u, v, w) := \Omega(u, v)(w), \quad \forall u, v, w \in \mathfrak{h}.
$$
The non-degeneracy of $\Omega$ is equivalent to the condition
$
\iota_u \tilde{\Omega} \neq 0, \quad \text{for all } u \in \mathfrak{h} \setminus \{0\},
$
where $\iota_u \tilde{\Omega}$ denotes the interior product (contraction) of $\tilde{\Omega}$ with $u$.
Conversely, let $\mathfrak{h}$ be a $\mathbb{Z}_2$-graded vector space and $\tilde{\Om} \in \mathcal{T}^3(\mathfrak{h})$ such that $\iota_u \tilde{\Om} \neq 0$ for every non-zero $u \in \mathfrak{h}$. Define a bilinear map $\Om : \mathfrak{h} \times \mathfrak{h} \to \mathfrak{h}^*$ by
$$
\Om(u,v)(w) := \tilde{\Om}(u, v,w), \quad \forall u,v,w \in \mathfrak{h}.
$$
Then $\Om$ is non-degenerate. Moreover, since $\tilde{\Om}\in \mathcal{T}(\mathfrak{h})$, then $\Om$ is cyclic, and thus defines a reduced  quadratic 2-step nilpotent symmetric Leibniz  algebra $T^*_\Om(\mathfrak{h})$.
Hence, we obtain a one-to-one correspondence between the set of reduced $T_\Om^*$-extensions of $\mathfrak{h}$ and the set
\begin{equation}
\left\{ \tilde{\Om} \in \mathcal{T}^3(\mathfrak{h}) \,\middle|\, \iota_u \tilde{\Om} \neq 0 \text{ for all } u \in \mathfrak{h} \setminus \{0\} \right\}.
\label{trino}\end{equation}
This leads to the following corollary:

\sssbegin{Corollary}
Let $(\A, \bullet_{\Om_1})$ and $(\mathcal{B}, \bullet_{\Om_2})$ be $T^*_\Omega$-extensions of a $\mathbb{Z}_2$-graded vector space $\mathfrak{h}$ associated with even trilinear forms $\tilde{\Omega}_1$ and $\tilde{\Omega}_2$, respectively. Then $(\A, \bullet_{\Om_1})$ and $(\mathcal{B}, \bullet_{\Om_2})$ are isometric-isomorphic if and only if there exists an even isomorphism $\varphi : \mathfrak{h} \to \mathfrak{h}$ such that
$$
\tilde{\Omega}_1(u,v,w) = \tilde{\Omega}_2(\varphi(u), \varphi(v), \varphi(w)), $$

for all $ u,v,w \in \mathfrak{h}$. In this case, we say that $\tilde{\Omega}_1$ and $\tilde{\Omega}_2$ are equivalent.
\end{Corollary}

We conclude with the following classification result:

\sssbegin{Corollary}
There is a bijective correspondence between the isometry classes of reduced quadratic 2-step nilpotent symmetric Leibniz algebras and the equivalence classes of even trilinear forms belonging to the set

$$
\left\{ \tilde{\Omega} \in \mathcal{T}^3(\mathfrak{h}) \,\middle|\, \iota_u \tilde{\Omega} \neq 0 \text{ for all } u \in \mathfrak{h} \setminus \{0\} \right\}.
$$

\end{Corollary}

\sssbegin{Corollary}
Let $(\A, \bullet, \prs)$ be a quadratic 2-step nilpotent symmetric Leibniz algebra. Then $(\A, \bullet, \prs)$ is isometrically isomorphic to $(\mathfrak{A} := I \oplus T^*_\Omega(U), \prs)$, where $I \subset \operatorname{Ann}(\mathfrak{A}, \bullet)$,
 $(I, \prs_I)$ is a pseudo-Euclidean $\mathbb{Z}_2$-graded vector space,
 and $\prs$ is defined by
$$
\langle x + u + f,\; y + v + g \rangle = \langle x, y \rangle_I + B(u + f,\; v + g),
$$
for all $x, y \in I$ and $u + f,\; v + g \in T^*_\Omega(U)$, where $B$ is the canonical invariant bilinear form on the $T^*_\Omega(U)$-extension.
\end{Corollary}

\begin{proof}
Since $(\A, \bullet)$ is 2-step nilpotent, we have
\(
\A \bullet \A\subseteq \operatorname{Ann}(\A, \bullet).
\)
According to Lemma~\ref{orth}, we also have
$
\A \bullet \A = \operatorname{Ann}(\A, \bullet)^\perp,
$
which implies that $\A \bullet \A$ is totally isotropic. Let $U$ be a complementary subspace of $\A \bullet \A$ in $\operatorname{Ann}(\A, \bullet)$. Since $\prs$ is non-degenerate, the restriction $\prs_{\mid U \times U}$ is also non-degenerate. Hence, $\A$ decomposes as
\[
\A = U \oplus U^\perp.
\]
Since $U \subseteq \operatorname{Ann}(\A, \bullet)$, it follows that $\A \bullet \A \subseteq U^\perp$, which implies that $U^\perp$ is a graded ideal of $(\A, \bullet)$. Consequently, the pair $(U^\perp, \prs_{U^\perp})$ forms a quadratic 2-step nilpotent symmetric Leibniz superalgebra,  where $\prs_{U^\perp} := \prs_{\mid U^\perp \times U^\perp}$  . Moreover, we have
\[
\operatorname{Ann}(U^\perp, \bullet) = U^\perp \bullet U^\perp,
\]
so $(U^\perp, \prs_{U^\perp})$ is reduced. According to Theorem~\ref{T-ex}, this reduced algebra is isometrically isomorphic to the $T_\Om^*$-extension $T^*_\Omega(\mathfrak{h})$ of $\mathfrak{h}$ by $\Om$. Therefore,  $(\A, \bullet, \prs)$ is isometrically isomorphic to $(U\oplus T^*_\Omega(\mathfrak{h}), \prs)$ with the given structure.
\end{proof}

Now, we provide the converse of Theorem \ref{theo Pi-extension} in the particular case of odd-quadratic 2-step nilpotent symmetric Leibniz superalgebras.

\sssbegin{Theorem}\label{Pi-ex}
Let $(\mathcal{A}, \bullet, \prs)$ be a reduced odd-quadratic 2-step nilpotent symmetric Leibniz superalgebra. Then, $(\mathcal{A}, \bullet)$ is a $\Pi(T^*_\Omega)$-extension of a trivial algebra $\mathfrak{h}$ by means of $\Omega$ where $\Om$ is non-degenerate.
\end{Theorem}
\begin{proof}
We show, similarly to the first part of the proof of Theorem \ref{T-ex}, that $\A \bullet \A = \mathrm{Ann}(\A, \bullet),
$ and that this space is totally isotropic of dimension $\frac{1}{2} \dim \A$. We can therefore decompose $\mathcal{A}$ as
$
\mathcal{A} = \operatorname{Ann}(\mathcal{A}, \bullet) \oplus \mathfrak{h},
$
where $\mathfrak{h}$ is a totally isotropic subspace of the same dimension. Define the map
$
\varphi: \operatorname{Ann}(\mathcal{A}, \bullet) \to \mathfrak{h}^*, \, \varphi(a)(u) := \langle a, u \rangle,
$ for all $a\in \operatorname{Ann}(\mathcal{A}, \bullet)$, $u\in \mathfrak{h}$,
which is clearly a linear isomorphism. Since the bilinear form $\prs$ is odd, it follows that $\varphi$ is a odd isomorphism. 

 Now, define the bilinear map $\Omega: \mathfrak{h} \times \mathfrak{h} \to \Pi(\mathfrak{h}^*)$ by
$
\Omega(u, v) := \Pi\circ \varphi(u \bullet v).
$ Since  $\varphi$ is odd and the product $\bullet$ is even maps, it follows that $\Omega$ is even. Since \( \varphi \) is bijective, then \( \Omega \) is also bijective.

On the other hand, for all $u, v, w \in \mathfrak{h}$, we have
\begin{align*}
\Pi(\Omega(u, v))(w)  = \langle u \bullet v, w \rangle  = \langle u, v \bullet w \rangle = (-1)^{|u|(|v|+|w|)} \langle v \bullet w, u \rangle = (-1)^{|u|(|v|+|w|)} \Pi(\Omega(v, w))(u).
\end{align*}
Now, consider $\mathfrak{h}$ as a trivial superalgebra (i.e., the product on $\mathfrak{h}$ is zero). Since $\Omega$ is an even bilinear map satisfying the cyclicity condition given by \eqref{cyclic1}, it follows from Theorem \ref{theo Pi-extension} that $\Pi(T^*_\Om(\mathfrak{h}))=\mathfrak{h} \oplus \Pi(\mathfrak{h}^*)$ is the $\Pi(T_\Omega^*)$-extension of the trivial superalgebra $\mathfrak{h}$ by means of $\Omega$.

Moreover, define the map 
$
\Phi: \mathfrak{h} \oplus \operatorname{Ann}(\mathcal{A}, \bullet) \to \Pi(T^*_\Omega(\mathfrak{h})),\,  \Phi(u + a) := u + \Pi(\varphi(a)).
$
This map is clearly linear and bijective. For all $u, v \in \mathfrak{h}$ and $a, b \in \operatorname{Ann}(\mathcal{A}, \bullet)$, we have:
\begin{align*}
\Phi((u + a) \bullet (v + b)) &= \Phi(u \bullet v) = \Pi(\varphi(u \bullet v)) = \Omega(u, v)\\
&= (u + \Pi(\varphi(a)) \bullet_\Omega (v + \Pi(\varphi(b))) \\
&= \Phi(u + a) \bullet_\Omega \Phi(v + b),
\end{align*}
so $\Phi$ is an isomorphism of Leibniz superalgebras. Let $B$ be the odd bilinear form on $\Pi(T^*_\Omega(\mathfrak{h}))$. Then for all $u+a \in (\mathfrak{h}\oplus \operatorname{Ann}(\mathcal{A}, \bullet))_{|u|}$ and $v +b \in (\mathfrak{h}\oplus \operatorname{Ann}(\mathcal{A}, \bullet))_{|v|},$ we have:
\begin{align*}
B(\Phi(u + a), \Phi(v + b)) &= B(u + \Pi(\varphi(a)), v + \Pi(\varphi(b))) \\
&= \varphi(a)(v) + (-1)^{|u||v|} \varphi(b)(u) \\
&= \langle a, v \rangle + (-1)^{|u||v|} \langle b, u \rangle \\
&= \langle u + a, v + b \rangle,
\end{align*}
so $\Phi$ is an isometric. Therefore, $(\mathcal{A}, \bullet, \prs)$ is isometric and isomorphic to the $\Pi(T^*_\Omega)$-extension of the trivial algebra $\mathfrak{h}$ by means of $\Omega$.
\end{proof}

Let $E = E_{\bar{0}} \oplus E_{\bar{1}}$ be a $\mathbb{Z}_2$-graded vector space. Let
$
\Omega : E \times E \to \Pi(E^*)
$ be an even bilinear map. We say that $\Omega$ is cyclic if, for all  $u,v,w \in E$, the following identity holds
\begin{align}
\Pi(\Omega(u, v))(w) = (-1)^{|u|(|v|+|w|)} \Pi(\Omega(v, w))(u).
\label{cycl}
\end{align}
We define the associated odd trilinear form $\tilde{\Omega} : E \times E \times E \to \mathbb{K}$ by
$$
\tilde{\Omega}(u,v,w) := \Pi(\Omega(u,v))(w), \quad \forall u,v,w \in E.
$$
Due to the cyclicity of $\Omega$, the form $\tilde{\Omega}$ satisfies the following identities:
\begin{align}
\tilde{\Omega}(u,v,w) = (-1)^{|u|(|v|+|w|)} \tilde{\Omega}(v,w,u) = (-1)^{|w|(|u|+|v|)} \tilde{\Omega}(w,u,v).
\label{tri1}
\end{align}
We denote by $\mathcal{T}^3(E)_{\bar{1}}$ the set of all odd trilinear forms $\tilde{\Omega} : E^{\otimes 3} \to \mathbb{K}$ that satisfy identity \eqref{tri1}. That is,
$$
\mathcal{T}^3(E)_{\bar{1}} := \left\{ \tilde{\Omega} \in \mathrm{Hom}(E^{\otimes 3}, \mathbb{K})_{\bar{1}} \ \middle| \ \tilde{\Omega} \text{ satisfies \eqref{tri1}} \right\}.
$$

Consequently, we have a restricted definition for the reduced odd-quadratic 2-step nilpotent symmetric Leibniz case as follows:

\sssbegin{Definition}
Let $\mathfrak{h}$ be a $\mathbb{Z}_2$-graded vector space and let $\Om : \mathfrak{h} \times \mathfrak{h} \to \Pi(\mathfrak{h}^*)$ be a non-degenerate cyclic even bilinear map. Define the $\mathbb{Z}_2$-graded vector space $\A := \mathfrak{h} \oplus \Pi(\mathfrak{h}^*)$ and endow it with the product
$$
(u + \Pi(f))\bullet ( v+ \Pi(g)) := \Om(u, v),
$$
and the symmetric bilinear form
$$
B(u + \Pi(f), v + \Pi(g)) := f(v) + (-1)^{|u||v|}g(u),
$$
for all $u+\Pi(f) \in \A_{|u|}$ and $v+\Pi(g) \in \A_{|v|}$. Then $(\A, B)$ is called a odd-quadratic  2-step nilpotent symmetric Leibniz  algebra. In this case, we say that $\A$ is the $\Pi(T_\Om^*)$-extension of $\mathfrak{h}$ by $\Om$.
\end{Definition}

\sssbegin{Theorem}
Let $(\A, \bullet_{\Om_1})$ and $(\mathcal{B}, \bullet_{\Om_2})$ be $\Pi(T_\Om^*)$-extensions of $\mathfrak{h}$ by   $\Om_1$ and $\Om_2$, respectively. Then:
\begin{enumerate}
    \item[\textup{(}i\textup{)}]  There exists an even isomorphism symmetric Leibniz algebra between $(\A, \bullet_{\Om_1})$ and $(\mathcal{B}, \bullet_{\Om_2})$ if and only if there exist an even isomorphism $\varphi_1 : \mathfrak{h} \to \mathfrak{h}$ and an even isomorphism $\varphi_2 : \Pi(\mathfrak{h}^*) \to \Pi(\mathfrak{h}^*)$ such that
$$
   \varphi_2(\Om_1(u,v)) = \Om_2(\varphi_1(u), \varphi_1(v)), \quad \forall u,v \in \mathfrak{h}.
   $$
\item[\textup{(}ii\textup{)}]  There exists an isometrically isomorphic (i.e., an isometric isomorphism preserving the natural bilinear form) between $\A$ and $\mathcal{B}$ if and only if there exists an even isomorphism $\varphi_1 : \mathfrak{h} \to \mathfrak{h}$ such that
$$
  \Pi(\Om_1(u,v)) = \Pi(\Om_2(\varphi_1(u), \varphi_1(v))) \circ \varphi_1, \quad \forall u,v \in \mathfrak{h}.
   $$
   \end{enumerate}
\end{Theorem}

\begin{proof}
(i)  Let $\varphi: \mathcal{A} \to \mathcal{B}$ be an even isomorphism of symmetric Leibniz algebras. Since the annihilators of both $\mathcal{A}$ and $\mathcal{B}$ satisfy
$     \operatorname{Ann}(\mathcal{A}, \bullet_{\Omega_1}) = \operatorname{Ann}(\mathcal{B}, \bullet_{\Omega_2}) = \Pi(\mathfrak{h}^*),
    $
and since these spaces are stable under $\varphi$, there exist even linear maps
$     \varphi_1 : \mathfrak{h} \to \mathfrak{h}, \; \varphi_1' : \mathfrak{h} \to \Pi(\mathfrak{h}^*), \; \text{and} \; \varphi_2 : \Pi(\mathfrak{h}^*) \to \Pi(\mathfrak{h}^*)
    $
such that, for all $u \in \mathfrak{h}_{|u|}$, $\Pi(f) \in \Pi(\mathfrak{h}^*)_{|u|}$,
$    \varphi(u + \Pi(f)) = \varphi_1(u) + \varphi_1'(u) + \varphi_2(\Pi(f)).
    $
It is immediate that $\varphi_2$ is an even linear isomorphism. We now show that $\varphi_1$ is also an even isomorphism of $\mathfrak{h}$.

Suppose \(\varphi_1(u) = 0\) for some \(u \in \mathfrak{h}\). Then, for any \(Y \in \mathcal{B}\),
\(
0 = \varphi(u) \bullet_{\Omega_2} Y = \varphi(u \bullet_{\Omega_1} \varphi^{-1}(Y)).
\)
Since \(\varphi\) is an isomorphism, this implies that
\(
u \bullet_{\Omega_1} \varphi^{-1}(Y) = 0 \; \text{for all } Y \in \mathcal{B}.
\)
Therefore, \(u \bullet_{\Omega_1} v = 0\) for all \(v \in \mathcal{A}\), hence \(u \in \operatorname{Ann}(\mathcal{A}, \bullet) = \Pi(\mathfrak{h}^*)\), which contradicts \(u \in \mathfrak{h}\). Thus, \(u = 0\), and \(\varphi_1\) is injective. As \(\dim \mathfrak{h} = \dim \varphi_1(\mathfrak{h})\), it follows that \(\varphi_1\) is an even isomorphism.

Now, for all \(u, v \in \mathfrak{h}\), and \(\Pi(f) \in \Pi(\mathfrak{h}^*)_{|u|}, \Pi(g) \in \Pi(\mathfrak{h}^*)_{|v|}\), we have:
\begin{align*}
    \varphi((u+\Pi(f)) \bullet_{\Omega_1} (v+\Pi(g))) &= \varphi(\Omega_1(u, v)) = \varphi_2(\Omega_1(u, v)), \\
    \varphi(u+\Pi(f)) \bullet_{\Omega_2} \varphi(v+\Pi(g)) &= (\varphi_1(u) + \varphi_1'(u) + \varphi_2(\Pi(f))) \bullet_{\Omega_2} (\varphi_1(v) + \varphi_1'(v) + \varphi_2(\Pi(g))) \\
    &= \Omega_2(\varphi_1(u), \varphi_1(v)).
\end{align*}
Hence,
\[
\varphi_2(\Omega_1(u, v)) = \Omega_2(\varphi_1(u), \varphi_1(v)) \quad \text{for all } u, v \in \mathfrak{h}.
\]
Conversely, suppose there exist even isomorphisms \(\varphi_1 : \mathfrak{h} \to \mathfrak{h}\) and \(\varphi_2 : \Pi(\mathfrak{h}^*) \to \Pi(\mathfrak{h}^*)\) such that
\[
\varphi_2(\Omega_1(u, v)) = \Omega_2(\varphi_1(u), \varphi_1(v)) \quad \text{for all } u, v \in \mathfrak{h}.
\]
Define \(\varphi : \mathcal{A} \to \mathcal{B}\) by
\(
\varphi(u + \Pi(f)) = \varphi_1(u) + \varphi_2(\Pi(f))\), for all $ u \in \mathfrak{h}, \, \Pi(f) \in \Pi(\mathfrak{h}^*)_{|u|}.$

It is straightforward to verify that \(\varphi\) is an even isomorphism of symmetric Leibniz algebras.

(ii) Suppose \(\varphi : \mathcal{A} \to \mathcal{B}\) is an isometric isomorphism. Then, as in (1), there exist even isomorphisms \(\varphi_1 : \mathfrak{h} \to \mathfrak{h}\) and \(\varphi_2 : \Pi(\mathfrak{h}^*) \to \Pi(\mathfrak{h}^*)\) such that
\[
\varphi(u + \Pi(f)) = \varphi_1(u) + \varphi_2(\Pi(f)) \quad \text{for all } u \in \mathfrak{h},\, \Pi(f) \in \Pi(\mathfrak{h}^*)_\alpha.
\]
For all \(u \in \mathfrak{h},\, \Pi(f) \in \Pi(\mathfrak{h}^*)_\alpha\), we have
\[
B'(\varphi(u), \varphi(\Pi(f))) = B(u, \Pi(f)) \Rightarrow (-1)^{\alpha |u|} \Pi(\varphi_2(\Pi(f)))(\varphi_1(u)) = (-1)^{\alpha |u|} f(u).
\]
Therefore,
\(
\varphi_2(\Pi(f)) = \Pi(f \circ \varphi_1^{-1}) \quad \text{for all } \Pi(f) \in \Pi(\mathfrak{h}^*)_\alpha.
\)
Since
\(
\varphi_2(\Omega_1(u, v)) = \Omega_2(\varphi_1(u), \varphi_1(v)),
\)
we obtain
\[
\Pi(\Omega_1(u, v)) = \Pi(\Omega_2(\varphi_1(u), \varphi_1(v))) \circ \varphi_1, \quad \forall u, v \in \mathfrak{h}.
\]

Conversely, suppose there exists an even isomorphism \(\varphi_1 : \mathfrak{h} \to \mathfrak{h}\) such that
\[
\Pi(\Omega_1(u, v)) = \Pi(\Omega_2(\varphi_1(u), \varphi_1(v))) \circ \varphi_1, \quad \forall u, v \in \mathfrak{h}.
\]
Define \(\varphi : \mathcal{A} \to \mathcal{B}\) by
\(
\varphi(u + \Pi(f)) = \varphi_1(u) + \Pi(f \circ \varphi_1^{-1}), \quad \forall u \in \mathfrak{h},\; \Pi(f) \in \Pi(\mathfrak{h}^*)_\alpha.
\)
Then it is clear that \(\varphi\) is an even isomorphism of symmetric Leibniz algebras.

Moreover, for all \(u + \Pi(f) \in \mathcal{A}_{|u|}\) and \(v + \Pi(g) \in \mathcal{A}_{|v|}\), we have:
\begin{align*}
    B'(\varphi(u + \Pi(f)), \varphi(v + \Pi(g))) &= B'(\varphi_1(u) + \Pi(f \circ \varphi_1^{-1}), \varphi_1(v) + \Pi(g \circ \varphi_1^{-1})) \\
    &= f(\varphi_1^{-1}(v)) + (-1)^{|u||v|} g(\varphi_1^{-1}(u)) \\
    &= f(v) + (-1)^{|u||v|} g(u) \\
    &= B(u + \Pi(f), v + \Pi(g)).
\end{align*}
Therefore, \(\varphi\) is an isometric isomorphism.
\end{proof}

Let $\mathfrak{h}$ be a $\mathbb{Z}_2$-graded vector space and let $\Om : \mathfrak{h} \times \mathfrak{h} \to \Pi(\mathfrak{h}^*)$ be a non-degenerate cyclic even bilinear map. We denote by $\tilde{\Omega}$ the associated odd trilinear form defined by
$$
\tilde{\Omega}(u, v, w) := \Pi(\Omega(u, v))(w), \quad \forall u, v, w \in \mathfrak{h}.
$$
The non-degeneracy of $\Omega$ is equivalent to the condition
$
\iota_u \tilde{\Omega} \neq 0, \quad \text{for all } u \in \mathfrak{h} \setminus \{0\},
$
where $\iota_u \tilde{\Omega}$ denotes the interior product (contraction) of $\tilde{\Omega}$ with $u$.
Conversely, let $\mathfrak{h}$ be a $\mathbb{Z}_2$-graded vector space and $\tilde{\Om} \in \mathcal{T}^3(\mathfrak{h})_{\bar{1}}$ such that $\iota_u \tilde{\Om} \neq 0$ for every non-zero $u \in \mathfrak{h}$. Define a bilinear map $\Om : \mathfrak{h} \times \mathfrak{h} \to \mathfrak{h}^*$ by
$$
\Om(u,v)(w) := \tilde{\Om}(u, v,w), \quad \forall u,v,w \in \mathfrak{h}.
$$
Then $\Om$ is non-degenerate. Moreover, since $\tilde{\Om}\in \mathcal{T}(\mathfrak{h})_{\bar{1}}$, then $\Om$ is cyclic, and thus defines a reduced  quadratic 2-step nilpotent symmetric Leibniz  algebra $\Pi(T^*_\Om(\mathfrak{h}))$.
Hence, we obtain a one-to-one correspondence between the set of reduced $\Pi(T_\Om^*)$-extensions of $\mathfrak{h}$ and the set
\begin{equation}
\left\{ \tilde{\Om} \in \mathcal{T}^3(\mathfrak{h})_{\bar{1}} \,\middle|\, \iota_u \tilde{\Om} \neq 0 \text{ for all } u \in \mathfrak{h} \setminus \{0\} \right\}.
\label{trino}\end{equation}

This leads to the following corollary:

\sssbegin{Corollary}
Let $(\A, \bullet_{\Om_1})$ and $(\mathcal{B}, \bullet_{\Om_2})$ be $\Pi(T^*_\Omega)$-extensions of a $\mathbb{Z}_2$-graded vector space $\mathfrak{h}$ associated with odd trilinear forms $\tilde{\Omega}_1$ and $\tilde{\Omega}_2$, respectively. Then $(\A, \bullet_{\Om_1})$ and $(\mathcal{B}, \bullet_{\Om_2})$ are isometric-isomorphic if and only if there exists an even isomorphism $\varphi : \mathfrak{h} \to \mathfrak{h}$ such that
$$
\tilde{\Omega}_1(u,v,w) = \tilde{\Omega}_2(\varphi(u), \varphi(v), \varphi(w)), $$
for all $ u,v,w \in \mathfrak{h}$. In this case, we say that $\tilde{\Omega}_1$ and $\tilde{\Omega}_2$ are equivalent.
\end{Corollary}

We conclude with the following classification result:

\sssbegin{Corollary}
There is a bijective correspondence between the isometry classes of reduced odd-quadratic 2-step nilpotent symmetric Leibniz algebras and the equivalence classes of odd trilinear forms belonging to the set

$$
\left\{ \tilde{\Omega} \in \mathcal{T}^3(\mathfrak{h})_{\bar{1}} \,\middle|\, \iota_u \tilde{\Omega} \neq 0 \text{ for all } u \in \mathfrak{h} \setminus \{0\} \right\}.
$$

\end{Corollary}

\sssbegin{Corollary}
Let $(\A, \bullet, \prs)$ be a odd-quadratic 2-step nilpotent symmetric Leibniz algebra. Then $(\A, \bullet, \prs)$ is isometrically isomorphic to $(\mathfrak{A} := I \oplus \Pi(T^*_\Omega(U)), \prs)$, where $I \subset \operatorname{Ann}(\mathfrak{A}, \bullet)$,
 $(I, \prs_I)$ is a odd pseudo-Euclidean $\mathbb{Z}_2$-graded vector space,
 and $\prs$ is defined by
$$
\langle x + u + \Pi(f),\; y + v + \Pi(g) \rangle = \langle x, y \rangle_I + B(u + \Pi(f),\; v + \Pi(g)),
$$
for all $x, y \in I$ and $u + \Pi(f),\; v + \Pi(g)\in \Pi(T^*_\Omega(U))$, where $B$ is the canonical invariant bilinear form on the $\Pi(T^*_\Omega(U))$-extension.
\end{Corollary}

\begin{proof}
Since $(\A, \bullet)$ is 2-step nilpotent, we have
\(
\A \bullet \A\subseteq \operatorname{Ann}(\A, \bullet).
\)
According to Lemma~\ref{orth}, we also have
$
\A \bullet \A = \operatorname{Ann}(\A, \bullet)^\perp,
$
which implies that $\A \bullet \A$ is totally isotropic. Let $U$ be a complementary subspace of $\A \bullet \A$ in $\operatorname{Ann}(\A, \bullet)$. Since $\prs$ is non-degenerate, the restriction $\prs_{\mid U \times U}$ is also non-degenerate. Hence, $\A$ decomposes as
\[
\A = U \oplus U^\perp.
\]
Since $U \subseteq \operatorname{Ann}(\A, \bullet)$, it follows that $\A \bullet \A \subseteq U^\perp$, which implies that $U^\perp$ is a graded ideal of $(\A, \bullet)$. Consequently, the pair $(U^\perp, \prs_{U^\perp})$ forms a odd-quadratic 2-step nilpotent symmetric Leibniz superalgebra,  where $\prs_{U^\perp} := \prs_{\mid U^\perp \times U^\perp}$  . Moreover, we have
\[
\operatorname{Ann}(U^\perp, \bullet) = U^\perp \bullet U^\perp,
\]
so $(U^\perp, \prs_{U^\perp})$ is reduced. According to Theorem~\ref{Pi-ex}, this reduced algebra is isometrically isomorphic to the $\Pi(T^*_\Om)$-extension $\Pi(T^*_\Omega(\mathfrak{h}))$ of $\mathfrak{h}$ by a $\Om$. Therefore,  $(\A, \bullet, \prs)$ is isometrically isomorphic to $(U\oplus \Pi(T^*_\Omega(\mathfrak{h})), \prs)$ with the given structure.
\end{proof}

\section{Double extensions of flat pseudo-Euclidean Leibniz superalgebras}
In this section, we introduce the notion of a double extension of flat even (resp. odd) pseudo-Euclidean left-Leibniz (resp. Lie) superalgebras by a one-dimensional $\mathbb{Z}_2$-graded vector space. This construction is carried out in two distinct steps: a central extension, followed by a semi-direct product. In the first step, we define the central extension of pre-left Leibniz superalgebra by a one-dimensional $\mathbb{Z}_2$-graded vector space. In the second step, we introduce the semi-direct product of pre-left Leibniz superalgebra  by a one-dimensional $\mathbb{Z}_2$-graded vector space. After establishing these two operations, we present the general notion of a double extension of flat pseudo-Euclidean left Leibniz superalgebras. This construction is inspired by the double extension method developed by  A. Medina  and P. Revoy in their study of quadratic Lie algebras \cite{Medina3}. We show that every flat even (resp. odd) pseudo-Euclidean non-Lie left Leibniz superalgebra can be obtained by a double extension procedure.

\subsection{Central extensions}

We now introduce the notion of central extension of pre-left Leibniz superalgebras.
\sssbegin{Proposition}
    Let $(\mathcal{A}, \star, \circ)$ be a pre-left Leibniz superalgebra, and let $V = \mathbb{K}e$ be a one-dimensional $\mathbb{Z}_2$-graded vector space. Let $\mu, \ga: \mathcal{A} \times \mathcal{A} \to \mathbb{K}$ be two homogeneous bilinear maps of given degrees. Consider the $\mathbb{Z}_2$-graded vector space $\tilde{\mathcal{A}} = \mathcal{A} \oplus \mathbb{K}e$, and define two bilinear operations $(\star_\mu, \circ_\ga)$ on $\tilde{\mathcal{A}}$ by
$$
(u + \alpha e) \star_\mu (v + \beta e) = u \star v + \mu(u, v)e, \quad 
(u + \alpha e) \circ_\ga (v + \beta e) = u \circ v + \ga(u, v)e,
$$
for all $u, v \in \mathcal{A}$ and $\alpha, \beta \in \mathbb{K}$. 
Then, $(\tilde{\mathcal{A}}, \star_\mu, \circ_\ga)$ defines a pre-left Leibniz superalgebra if and only if the bilinear maps $\mu$ and $\ga$ satisfy the following identities:
\begin{align}
\mu(u \star v, w) - \mu(u, v \star w) &= -(-1)^{|u||v|} \left( \mu(v \circ u, w) + \mu(v, u \star w) \right), \label{E1} \\
\ga(u \star v, w) - \mu(u, v \circ w) &= -(-1)^{|u||v|} \left( \ga(v \circ u, w) - \ga(v, u \circ w) \right), \label{E2} \\
\ga(u, v \circ w) &= -\ga(u, v \star w), \label{E3}
\end{align}
for all  $u, v, w \in \mathcal{A}$.

In this case, $(\tilde{\mathcal{A}}, \star_\mu, \circ_\nu)$ is called  the central extension of $\mathcal{A}$ by $\mathbb{K}e$ by means of  $(\mu, \nu)$.
\end{Proposition}
\begin{proof}
    Straightforward calculation.
\end{proof}
\sssbegin{Remark}{\rm 
     When the pre-left Leibniz superalgebra $(\mathcal{A}, \star, \circ)$ satisfies $\star = -\circ$, i.e., when $\mathcal{A}$ is a left-symmetric superalgebra, the conditions \eqref{E1}-\eqref{E3} reduce to a single identity, as shown in the corollary below.}
\end{Remark}

Given that left symmetric superalgebras constitute a subclass of pre-left Leibniz superalgebra, we also consider the particular case of central extension of left symmetric superalgebars.

\sssbegin{Corollary}
   Let $(\mathcal{A}, \star)$ be a left symmetric superalgebra, and let $V = \mathbb{K}e$ be a one-dimensional $\mathbb{Z}_2$-graded vector space. Let $\mu: \mathcal{A} \times \mathcal{A} \to \mathbb{K}$ be a homogeneous bilinear map of given degree. Consider the $\mathbb{Z}_2$-graded vector space $\tilde{\mathcal{A}} = \mathcal{A} \oplus \mathbb{K}e$, endowed with the product defined by 
$$
(u + \alpha e) \star_\mu (v + \beta e) = u \star v + \mu(u, v)e,
$$
for all $u, v \in \mathcal{A}$ and $\alpha, \beta \in \mathbb{K}$. 
Then, $(\tilde{\mathcal{A}}, \star_\mu)$ is a left symmetric superalgebra if and only if the bilinear map $\mu$ and  satisfies the following identity (for all  $u, v, w \in \mathcal{A}$):
\begin{align}
\mu(u \star v, w) - \mu(u, v \star w) &= (-1)^{|u||v|} \left( \mu(v \star u, w) - \mu(v, u \star w) \right).
\label{E11} 
\end{align}
In this case, $(\tilde{\mathcal{A}}, \star_\mu)$ is called the central extension of $\mathcal{A}$ by $\mathbb{K}e$ by means of $\mu$.   
\end{Corollary}
\begin{proof}
  A left-symmetric superalgebra is a particular case of a pre-left Leibniz superalgebra $(\mathcal{A}, \star, \circ)$ in which the two products are related by $\star = -\circ$. In this case, the relations \eqref{E1}-\eqref{E3} reduce to a single identity. Thus, the proof follows.  
\end{proof}

\sssbegin{Proposition}\label{Auconidition}
Let $(\mathcal{A}, \bullet, \prs)$ be a flat pseudo-Euclidean left Leibniz superalgebra, and  let $(\star, \circ)$ be its associated Levi-Civita products. Let $\mu, \ga : \mathcal{A} \times \mathcal{A} \to \mathbb{K}$ be two bilinear maps of degrees $\alpha$ and $\beta$, respectively. Then, there exist two homogeneous endomorphisms $G, D: \mathcal{A} \to \mathcal{A}$ of degrees $\alpha$ and $\beta$, such that

$$
\mu(u, v) = \langle G(u), v \rangle, \quad \ga(u, v) = \langle D(u), v \rangle, \quad \text{for all } u, v \in \mathcal{A}.
$$
The bilinear maps $\mu$ and $\ga$ satisfy relations \eqref{E1}–\eqref{E3} if and only if the following conditions hold for all  $u, v \in \mathcal{A}$
\begin{itemize}
\item[\textup{(}i\textup{)}] $G(u \bullet v) = -(-1)^{|u|\al} u \circ G(v) + (-1)^{|v|(\alpha + |u|)} v \circ G(u)$,
\item[\textup{(}ii\textup{)}] $D(u \bullet v) = (-1)^{|u|\beta} u \star D(v) + (-1)^{|v|(\al+|u| )} v \star G(u)$,
\item[\textup{(}iii\textup{)}] $v \star D(u) = -v \circ D(u)$.
\end{itemize}
\end{Proposition}

\begin{proof}
We prove assertion (i), the other relations follow by similar computations.

Let $u, v, w \in \mathcal{A}$. Then,
\begin{align*}
&\mu(u \star v, w) - \mu(u, v \star w) + (-1)^{|u||v|} \left( \mu(v \circ u, w) + \mu(v, u \star w) \right) \\
&= \langle G(u \star v), w \rangle - \langle G(u), v \star w \rangle + (-1)^{|u||v|} \left( \langle G(v \circ u), w \rangle + \langle G(v), u \star w \rangle \right).
\end{align*}
Using the \eqref{compatibl} and symmetry of the bilinear form $\prs$
\begin{align*}
\langle G(u), v \star w \rangle &= (-1)^{(|v| + |w|)(\alpha + |u|)} \langle v \star w, G(u) \rangle \\
&= (-1)^{|v||w| + |v|(\alpha + |u|)} \langle w, v \circ G(u) \rangle = (-1)^{|v|(\alpha + |u|)} \langle v \circ G(u), w \rangle.
\end{align*}
Similarly,
$$
\langle G(v), u \star w \rangle = (-1)^{|u|(\alpha + |v|)} \langle u \circ G(v), w \rangle.
$$
Combining all terms, we obtain
\begin{align*}
&\langle G(u \star v), w \rangle - (-1)^{|v|(\alpha + |u|)} \langle v \circ G(u), w \rangle + (-1)^{|u||v|} \langle G(v \circ u), w \rangle \\
&\qquad + (-1)^{|u|\alpha} \langle u \circ G(v), w \rangle \\
&= \left\langle G(u \bullet v) - (-1)^{|v|(\alpha + |u|)} v \circ G(u) + (-1)^{|u|\alpha} u \circ G(v), w \right\rangle.
\end{align*}
Since this equality holds for all $w \in \mathcal{A}$, it follows that

$$
G(u \bullet v) = -(-1)^{|v|(\alpha + |u|)} u \circ G(v) + (-1)^{|u|\alpha} v \circ G(u).
$$
The proofs of (ii) and (iii) are similar and left to the reader.
\end{proof}
\sssbegin{Remark}{\rm \label{Remark1}
    When the left Leibniz superalgebra $(\mathcal{A}, \bullet)$ is anticommutative, i.e., when $(\mathcal{A}, \bullet)$ is a Lie superalgebra, then the Levi-Civita products associated with the pseudo-Euclidean Lie superalgebra satisfy $\star = -\circ$. In this case, the conditions (i)-(iii) reduce to a single identity, as stated in the corollary below.}
\end{Remark}
Given that flat Lie superalgebras form a subclass of flat left Leibniz superalgebras, we now consider the particular case of Proposition \ref{Auconidition}.
\sssbegin{Corollary}
  Let $(\mathcal{A}, \br, \prs)$ be a flat pseudo-Euclidean Lie  superalgebra, and  let $\star$ be its associated Levi-Civita product. Let $\mu: \mathcal{A} \times \mathcal{A} \to \mathbb{K}$ be a bilinear map of degree $\alpha$. Then, there exists a homogeneous endomorphism $G: \mathcal{A} \to \mathcal{A}$ of degree $\alpha$, such that
$
\mu(u, v) = \langle G(u), v \rangle$,   for all  $u, v \in \mathcal{A}.
$
The bilinear map $\mu$ satisfies relation \eqref{E11} if and only if the following condition holds for all  $u, v \in \mathcal{A}$
 \begin{equation}G([u , v]) = (-1)^{|u|\alpha } u \star G(v) - (-1)^{|v|(\alpha + |u|)} v \star G(u).\label{cocycle}\end{equation}
\end{Corollary}
\begin{proof}
From  Remark \ref{Remark1}, it follows that $\star = -\circ$, and hence, by Proposition \ref{Auconidition}, we obtain the result.
\end{proof}
\sssbegin{Proposition}\label{Auconidition1}
Let $(\mathcal{A}, \bullet, \prs)$ be a flat even pseudo-Euclidean left Leibniz superalgebra, and  let $(\star, \circ)$ be its associated Levi-Civita products. Let $\mu, \ga : \mathcal{A} \times \mathcal{A} \to \mathbb{K}$ be two bilinear maps of degrees $\alpha$ and $\beta$, respectively. Then, there exist two homogeneous endomorphisms $G, D: \mathcal{A} \to \mathcal{A}$ of degrees $\alpha$ and $\beta$, such that
$$
\mu(u, v) =(-1)^{|v|} \langle G(u), v \rangle, \quad \ga(u, v) = -(-1)^{|u|}\langle D(u), v \rangle, \quad \text{for all } u, v \in \mathcal{A}.
$$
The bilinear maps $\mu$ and $\ga$ satisfy relations \eqref{E1}-\eqref{E3} if and only if the following conditions hold for all  $u, v \in \mathcal{A}$
\begin{itemize}
\item[\textup{(}i\textup{)}] $G(u \bullet v) = -(-1)^{|u|(1+\alpha )} u \circ G(v) + (-1)^{|v|(1+\alpha+|u|)} v \circ G(u)$,
\item[\textup{(}ii\textup{)}] $\langle D(u \bullet v), w\rangle = (-1)^{|u|(\beta+1)}\langle u \star D(v), w\rangle - (-1)^{|u|+|v|(|u| + \alpha)+|w|} \langle v \star G(u), w\rangle$,
\item[\textup{(}iii\textup{)}] $v \star D(u) = -v \circ D(u)$.
\end{itemize}
\end{Proposition}

\begin{proof}
The proof of this proposition is similar to that of Proposition~\ref{Auconidition}.
\end{proof}

\begin{Remark}
Assuming the same hypotheses as in the preceding proposition, if the degrees of $G$ and $D$ are odd, then condition \textup{(ii)} is equivalent to the following identity:
\[
 D(u \bullet v)  =  u \star D(v) + (-1)^{|u||v|}  v \star G(u),
\quad \text{for all } u, v, w \in \mathcal{A}.
\]
Indeed, we distinguish two cases:
\begin{itemize}
    \item If $w$ is even, then for $\langle v \star G(u), w \rangle \neq 0$, we must have $|v \star G(u)| + |w| = \bar{0}$, i.e., $|v| + |G(u)| + |w| = \bar{0}$. Since $G$ is odd, this implies $|u| + |v| + |w| = \bar{1}$.
    \item If $w$ is odd, then similarly, the condition $\langle v \star G(u), w \rangle \neq 0$ leads to $|u| + |v| + |w| = \bar{1}$.
\end{itemize}
Thus, in both cases, the parity of $|u| + |v| + |w|$ is $\bar{1}$, which yields the stated identity.
\end{Remark}

\subsection{Semi-direct product for pre-left Leibniz superalgebras }
We now introduce a new type of extension for pre-left Leibniz superalgebras.

Let $(\mathcal{A}, \star, \circ)$ be a pre-left Leibniz superalgebra, and let $V = \mathbb{K}d$ be a one-dimensional $\mathbb{Z}_2$-graded vector space. Let $\delta, \xi, D, G: \mathcal{A} \to \mathcal{A}$ be homogeneous endomorphisms of degree $|d|$, and let $a_0, b_0 \in  \mathcal{A}_{\bar{0}}$ and $\al, \beta \in \K$.

We consider the $\mathbb{Z}_2$-graded vector space $\bar{\mathcal{A}} = \mathcal{A} \oplus V$, and define two bilinear operations $\bar{\star}$ and $\bar{\circ}$ on $\bar{\mathcal{A}}$ by (for all $u, v \in \mathcal{A}$):
\begin{equation} \label{eq:semi-extension}
\begin{cases}
d \bar{\star}  d = a_0+\al d, \quad 
u \bar{\star}  v = u \star v, \\[1mm]
d  \bar{\star}  u = \delta(u), \quad \quad\;\;
u  \bar{\star} d = D(u),
\end{cases}
\quad \quad \text{and} \quad \quad
\begin{cases}
d  \bar{\circ}  d = b_0+\beta d, \quad 
u  \bar{\circ}  v = u \circ v, \\[1mm]
d  \bar{\circ}  u = \xi(u), \quad  \quad\;\;
u  \bar{\circ}  d = G(u).
\end{cases}
\end{equation}
\sssbegin{Proposition}
The $\mathbb{Z}_2$-graded vector space $\bar{\mathcal{A}} = \mathcal{A} \oplus \mathbb{K}d$, equipped with the bilinear products defined in \eqref{eq:semi-extension}, is a pre-left Leibniz superalgebra if and only if the following conditions are satisfied for all $u, v \in \mathcal{A}$:
\begin{align*}
D(u \bullet v) &= u \star D(v)  - (-1)^{|u||v|}v \star D(u),\\
D(u) \star v - u \star \delta(v) &= -(-1)^{|u||d|} \left(\xi(u) \star v + \delta(u \star v)\right), \\
D^2(u) - u \star a_0-\al D(u) &= -(-1)^{|u||d|} \left(D(\xi(u)) + \delta(D(u))\right), \\
\delta(u) \star v - \delta(u \star v) &= -(-1)^{|u||d|} \left(G(u) \star v + u \star \delta(v)\right), \\
D(\delta(u)) - \delta(D(u)) &= -(-1)^{|d||u|} \left(D(G(u)) + u \star a_0+\al D(u)\right), \\
(1-(-1)^{|d|}) \delta^2(u)-\de(\al+(-1)^{|d|}\beta) &=  a_0 \star u +(-1)^{|d|} b_0\star u, \\
D(a_0+(-1)^{|d|}b_0) -(1-(-1)^{|d|}) \delta(a_0) &= -(-1)^{|d|}(\al+\beta)a_0, \\
G(u \bullet v)&=  u \star G(v)  +(-1)^{|u||v|}  v \circ G(u), \\
D(u) \circ v - u \star \xi(v) &= -(-1)^{|u||d|} \left(\xi(u) \circ v - \xi(u \circ v)\right), \\
G(D(u)) - u \star b_0-\beta D(u) &= -(-1)^{|u||d|} \left(G(\xi(u)) - \xi(G(u))\right), \\
\delta(u) \circ v - \delta(u \circ v) &= -(-1)^{|u||d|} \left(G(u) \circ v - u \circ \xi(v)\right), \\
G(\delta(u)) - \delta(G(u)) &= -(-1)^{|u||d|} \left(G^2(u) - u \circ b_0-\beta G(u)\right), \\
a_0 \circ u+\al\xi(u) - \delta(\xi(u)) &= -(-1)^{|d|} \left(b_0 \circ u+\beta \xi(u) - \xi^2(u)\right), \\
G(a_0) - \delta(b_0)+\al b_0-\beta a_0 &= -(-1)^{|d|} \left(G(b_0) - \xi(b_0)\right), \\
u \circ G(v) = -u \circ D(v), \quad u \circ \xi(v) &= -u \circ \delta(v), \quad u \circ a_0 +u \circ b_0=-(\al+\beta)G(u), \\
\xi(u \circ v) = -\xi(u \star v), \, \xi(G(u)) = -\xi(&D(u)), \; \xi^2(u) = -\xi(\delta(u)), \; \xi(b_0 + a_0) =-(\al+\beta) b_0,\\ &\al^2=\beta^2=-\al\beta.
\end{align*}
where, $u\bullet v=u\star v+(-1)^{|u||v|} v\circ u$, for all $u,v\in \A$.

In this case, the  $(\delta, D, \xi, G, a_0, b_0, \al, \beta)$ is called an ``admissible", and the pre-left Leibniz superalgebra $\bar{\mathcal{A}}$ is called the semi-direct products of $\mathcal{A}$ by $\mathbb{K}d$ by means $(\delta, D, \xi, G, a_0, b_0, \al, \beta)$.
\end{Proposition}

\begin{proof}
The result follows by direct verification of the pre-left Leibniz superidentities using the products definitions in \eqref{eq:semi-extension}.
\end{proof}
\sssbegin{Remark}
In the preceding proposition, we have $\alpha^2 = \beta^2 = -\alpha \beta$, from which we deduce that either $\alpha = -\beta$ or $\alpha = \beta = 0$. Thus, we get $\al=-\beta$.
\end{Remark}

\subsection{Semi-direct product for left-symmetric superalgebras} We now present the semi-direct product in the case of left-symmetric superalgebras.

Let $(\mathcal{A}, \star)$ be a left-symmetric superalgebra, and let $V = \mathbb{K}d$ be a one-dimensional $\mathbb{Z}_2$-graded vector space. Let $\delta,  D: \mathcal{A} \to \mathcal{A}$ be homogeneous endomorphisms of degree $|d|$, and let $a_0 \in  \mathcal{A}_{\bar{0}}$ and $\al \in \K$.

We consider the $\mathbb{Z}_2$-graded vector space $\bar{\mathcal{A}} = \mathcal{A} \oplus V$, and define a bilinear operation $\bar{\star}$  on $\bar{\mathcal{A}}$, for all $u, v \in \mathcal{A}$ by:
\begin{equation} \label{eq:semi-extensionsymmetric}
d \bar{\star}  d = a_0+\al d, \qquad   
u \bar{\star}  v = u \star v, \qquad 
d  \bar{\star}  u = \delta(u), \qquad   
u  \bar{\star} d = D(u).
\end{equation}
\sssbegin{Corollary}
   The $\mathbb{Z}_2$-graded vector space $\bar{\mathcal{A}} = \mathcal{A} \oplus \mathbb{K}d$, equipped with the bilinear product defined in \eqref{eq:semi-extensionsymmetric}, is a left-symmetric superalgebra if and only if the following conditions are satisfied for all $u, v \in \mathcal{A}$ \textup{(}where $[u,v]=u\star v-(-1)^{|u||v|} v\star u$, for all $u,v\in \A$\textup{)}:
\begin{equation}
\begin{cases}
D([u , v]) = u \star D(v)  - (-1)^{|u||v|}v \star D(u),\\
D(u) \star v = u \star \delta(v)  +(-1)^{|u||d|} \de(u) \star v - (-1)^{|u||d|}\delta(u \star v), \\
D\circ \de(u)-\de\circ D(u)=(-1)^{|u||d|}(
D^2(u) - u \star a_0-\al D(u)), \\
(1-(-1)^{|d|})(\de^2(u)-\al \de(u)-a_0\star u)=0,\quad 
(1-(-1)^{|d|})(D-\de)(a_0)=0.
\end{cases}\end{equation}
In this case, the  $(\delta, D, a_0, \al)$ is called an ``admissible", and the left-symmetric superalgebra $\bar{\mathcal{A}}$ is called the semi-direct products of $\mathcal{A}$ by $\mathbb{K}d$ by means $(\delta, D, a_0, \al)$. 
\end{Corollary}
\begin{proof}
A left-symmetric superalgebra is a particular case of a pre-left Leibniz superalgebra $(\mathcal{A}, \star, \circ)$ in which the two products are related by $\star = -\circ$. From this, we deduce the following result.
\end{proof}
\subsection{Double extension of flat even pseudo-Euclidean left-Leibniz superalgebras} 
Now, we are in a position to introduce the double extensions of flat pseudo-Euclidean left Leibniz superalgebras.
We start by presenting the concept of even double extension of these superalgebras by the one-dimensional Lie algebra.

\sssbegin{Theorem}\label{Thdbeven}
Let $(\mathcal{A}, \bullet, \prs_\A)$ be a flat even pseudo-Euclidean left Leibniz superalgebra, and let $(\star, \circ)$ be the associated Levi-Civita products. Let $V = \mathbb{K} d$ be a one-dimensional  vector space and $V^* = \mathbb{K} e$ its dual. Suppose $(\delta, D, \delta^*, G, a_0, b_0, \al, -\al)$ is an admissible of $\mathcal{A}$ satisfying the compatibility conditions: 
\begin{equation}\label{desapp}
\begin{array}{lclrcl}
  \Rr_{c_0}^\circ + (R_{c_0}^\star)^* &=&  D^* \circ G- G^*\circ D,& [G, \de^*] &=&-G\circ D+\al G-\Rr_{b_0}^\circ,\; \\G\circ \de+\de^*\circ G &=& -G^2-\al G+\Rr_{b_0}^\circ,& [D, \de] & = &-D\circ G-\al D+\Rr_{b_0}^\star,\\ D\circ \de^*-\de\circ G & = & -D^2+\al D+\Rr_{a_0}^\star,& G^*\circ D & = &-D^*\circ D,\\ \de\circ D&=&-\de^*\circ D, &
\Rr_{a_0}^\star+\Rr_{a_0}^\circ & = & -\alpha (D+G), \\ \Rr_{c_0}^\circ + (\Rr_{c_0}^\star)^* & = &  D^* \circ G+ G^* \circ G, &
G(a_0+b_0)& = &0,\\  D^*(a_0+b_0)& =&0 , &(G^*+D^*)(a_0)& = &-2\al c_0,\\  (G+\de)(a_0+b_0)& =&0 , & (\de+\de^*)(a_0)& = &-\al(a_0+b_0) , \\
\langle a_0+b_0, c_0\rangle_\A& =  &0,& \langle a_0+b_0, a_0\rangle_\A & =  &0,\\ 
  \langle a_0+b_0, b_0\rangle_\A& = &0,& 
    (\de+D^*)(c_0)& =& G^*(a_0+b_0) , \\ (\de^*+G^*)(c_0)&=&-G^*(a_0)+D^*(b_0) ,&  (\de^*+G^*)(c_0)& =& (G^*+D)(b_0)-2\al c_0,\\ D^*(b_0+c_0)&=&2\al c_0-\de(c_0)+G^*(a_0), & D^* \circ G- G^*\circ D&=& G^* \circ G+D^*\circ G,\\ u\star G(v)&=&u\circ D(v)=-u\star D(v)\\&=&-u\circ G(v).
    \end{array}
\end{equation}
Define the $\mathbb{Z}_2$-graded vector space $\bar{\mathcal{A}} = \K e\oplus \mathcal{A} \oplus \K d$, equipped with the bilinear products $\bar{\star}$ and $\bar{\circ}$ defined for all $u, v \in \mathcal{A}$ by:
\begin{equation}\label{eq:extension}
\begin{cases}
d \bar{\star} d = a_0 + \lambda e+\al d, \\
u \bar{\star} v = u \star v + \langle G(u), v \rangle_\A e, \\
d \bar{\star} u = \delta(u) + \langle b_0, u \rangle_\A e, \\
u \bar{\star} d = D(u) + \langle c_0, u \rangle_\A e,\\
d\bar{\star} e=-\al e, \; e\bar{\star} d= 0,
\end{cases}
\quad \quad \text{and} \quad \quad
\begin{cases}
d \bar{\circ} d = b_0 + \lambda e-\al d, \\
u \bar{\circ} v = u \circ v + \langle D(u), v \rangle_\A e, \\
d \bar{\circ} u = \delta^*(u) + \langle a_0, u \rangle_\A e, \\
u \bar{\circ} d = G(u) + \langle c_0, u \rangle_\A e,\\
d\bar{\circ} e=\al e, \; e\bar{\circ} d=0,
\end{cases}
\end{equation}
where $\al, \lambda \in \mathbb{K}$ and $c_0 \in \mathcal{A}_{\bar{0}}$.
Then, $(\bar{\mathcal{A}}, \bar{\star}, \bar{\circ})$ is a pre-left Leibniz superalgebra. Moreover, the even bilinear form $\prs : \bar{\mathcal{A}} \times \bar{\mathcal{A}} \to \mathbb{K}$ defined by
$$
\prs|_{\mathcal{A} \times \mathcal{A}} = \prs_\A, \quad \langle d, e\rangle = \langle e, d\rangle= 1,
$$
is a pseudo-Euclidean scalar product on $\bar{\mathcal{A}}$.

The flat even pseudo-Euclidean left Leibniz superalgebra $(\bar{\mathcal{A}},\bar{\bullet}, \prs)$ is called the even double extension of $(\mathcal{A}, \bullet, \prs_\A)$ by means of $(\delta, D, \delta^*, G, a_0, b_0, c_0,\al,  \la)$, where $u\bar{\bullet} v=u\bar{\star} v+(-1)^{|u||v|} v\bar{\circ} u$, for all $u,v\in \bar{\A}$.
\end{Theorem}

\begin{proof}
Since $(\delta, D, \delta^*, G, a_0, b_0, \al, -\al)$ is an admissible of $\mathcal{A}$, and $u\star G(v)=u\circ D(v)=-u\star D(v)=-u\circ G(v)$, for all $u,v\in\A$, then the pair $(G, D)$ satisfies conditions (i), (ii), and (iii) of Proposition~\ref{Auconidition}, it follows from Proposition \ref{Auconidition} that the two even bilinear maps
$$
\mu,\, \gamma : \A \times \A \longrightarrow \mathbb{K}, \quad \mu(u, v) = \langle G(u), v \rangle_\A, \quad \gamma(u, v) = \langle D(u), v \rangle_\A,
$$
satisfy the relations \eqref{E1}–\eqref{E3}. Consequently, the $\mathbb{Z}_2$-graded vector space $\tilde{\mathcal{A}} = \A \oplus \mathbb{K}e$, equipped with the products
$$
(u + a e)\, \bar{\star}\, (v + b e) = u \star v + \langle G(u), v \rangle_\A e, \quad (u + a e)\, \bar{\circ}\, (v + b e) = u \circ v + \langle D(u), v \rangle_\A e,
$$
for all $u, v \in \mathcal{A}$ and $a, b \in \mathbb{K}$, is a pre-left Leibniz superalgebra, central extension of $\mathcal{A}$ by means of  $(\mu, \gamma)$.

We now define even endomorphisms $\delta', D', \xi', G' : \tilde{\mathcal{A}} \to \tilde{\mathcal{A}}$ by: 
$$
\begin{aligned}
&\delta'(u + a e) = \delta(u) + (\langle b_0, u \rangle_\A-a\al) e, \quad
D'(u + a e) = D(u) +\langle c_0, u \rangle_\A  e, \\
&\xi'(u + a e) = \delta^*(u) + (\langle a_0, u \rangle_\A+a\al ) e, \quad
G'(u + a e) = G(u) + \langle c_0, u \rangle_\A e,
\end{aligned}
$$
and set
$$
a_0' = a_0 + \lambda e, \quad b_0' = b_0 + \lambda e,
$$

for all $u \in \A$, $a \in \mathbb{K}$.

To show that $(\delta', D', \xi', G', a_0', b_0', \al, -\al)$ is an admissible  of $\tilde{\mathcal{A}}$, we verify one of the key admissibility conditions. Let $u, v \in \mathcal{A}$. We have,
$$
\begin{aligned}
D'&(u \bar{\star} v) - u \bar{\star} D'(v) + (-1)^{|u||v|} \left( D'(v \bar{\circ} u) + v \bar{\star} D'(u) \right) \\
&= D(u \star v) + \langle c_0, u \star v \rangle_\A e - u \star D(v) - \langle G(u), D(v) \rangle_\A e \\
&\quad + (-1)^{|u||v|} \left( D(v \circ u) + \langle c_0, v \circ u \rangle_\A e + v \star D(u) + \langle G(v), D(u) \rangle_\A e \right) \\
&= \left( D(u \star v) - u \star D(v) + (-1)^{|u||v|}( D(v \circ u) + v \star D(u) ) \right) \\
&\quad + \text{ ( } \langle  u \star v, c_0 \rangle_\A +(-1)^{|u||v|} \langle v\star  c_0, u  \rangle_\A - \langle G(u), D(v) \rangle_\A +  (-1)^{|u||v|} \langle G(v), D(u) \rangle_\A \text{)} e.
\end{aligned}
$$

Since $(\delta, D, \delta^*, G, a_0, b_0, \al, -\al)$ is an admissible of $\mathcal{A}$, the expression in the first parentheses vanishes. Therefore, we have

$$
D'(u \bar{\star} v) - u \bar{\star} D'(v) + (-1)^{|u||v|} \left( D'(v \bar{\circ} u) + v \bar{\star} D'(u) \right)
= \Theta(u, v)\, e,
$$
where
$$
\Theta(u, v) = \langle  u \star v ,c_0\rangle_\A + (-1)^{|u||v|} \langle v\star c_0,  u \rangle_\A - \langle G(u), D(v) \rangle_\A + (-1)^{|u||v|} \langle G(v), D(u) \rangle_\A.
$$
Now, using dualities $\langle G(v), D(u) \rangle_\A = \langle v, G^* D(u) \rangle_\A=(-1)^{|u||v|}\langle  G^* D(u), v \rangle_\A$ and similar for the others, we rewrite
$$
\Theta(u, v) = \langle  \Rr_{c_0}^\circ(u) + (R_{c_0}^\star)^*(u) + G^* D(u) - D^* G(u), v \rangle_\A.
$$
Since $\prs_\A$ is nondegenerate, it follows that $$\Rr_{c_0}^\circ(u) + (R_{c_0}^\star)^*(u) + G^* D(u) - D^* G(u) = 0$$ if and only if $\Theta(u, v) =0$, for all $u,v\in \A$. This proves the desired compatibility condition for $D'$. The remaining admissibility conditions for $(\delta', \xi', G')$ can be verified similarly using the same type of computation.
Hence, $(\delta', D', \xi', G', a_0', b_0', \al, -\al)$ is an admissible of $\Tilde{\mathcal{A}}$. We may therefore consider the pre-left Leibniz superalgebra $\bar{\mathcal{A}} = V^* \oplus \mathcal{A} \oplus V$, which is the semi-direct product of the central extension $\tilde{\mathcal{A}} = \mathcal{A} \oplus V^*$ by the one-dimensional vector space $V = \mathbb{K}d$.

Furthermore, it is straightforward to verify that the bilinear form
$$
\prs : \bar{\mathcal{A}} \times \bar{\mathcal{A}} \to \mathbb{K}, \quad \text{defined by } \prs|_{\mathcal{A} \times \mathcal{A}} = \prs_\A, \quad \langle e, d \rangle = 1,
$$
is an even pseudo-Euclidean scalar product on $\bar{\mathcal{A}}$. Moreover, the compatibility condition
$$
\langle u \bar{\star} v, w \rangle = (-1)^{|u||v|} \langle v, u \bar{\circ} w \rangle, \quad \forall u, v, w \in \bar{\mathcal{A}},
$$
is satisfied.

Hence, $(\bar{\mathcal{A}}, \bar{\bullet}, \prs)$ is a flat even pseudo-Euclidean left Leibniz superalgebra, where  $u\bar{\bullet} v=u\bar{\star} v+(-1)^{|u||v|} v\bar{\circ} u$, for all $u,v\in \bar{\A}$.
\end{proof}

Since flat even pseudo-Euclidean Lie superalgebras constitute a particular case of flat even pseudo-Euclidean left Leibniz superalgebras, we now present the procedure of double extension for flat even pseudo-Euclidean Lie superalgebras. Let us recall that the notion of double extension of flat pseudo-Euclidean Lie algebras was introduced by A. Aubert and A. Medina in \cite{Medina2}.
\sssbegin{Corollary} \label{dbexLie1}
Let $(\mathcal{A}, \br, \prs_\A)$ be a flat even pseudo-Euclidean Lie superalgebra, and let $\star$ be its associated Levi-Civita product. Let $V = \mathbb{K} d$ be a one-dimensional  vector space and $V^* = \mathbb{K} e$ its dual. Suppose two endomorphisms $D$ and $\delta$ of $\mathcal{A}$, and an element $a_0 \in \mathcal{A}$, satisfy the following relations for all $u, v \in \mathcal{A}$, 
\begin{equation} \label{mapsLie}
\begin{aligned}
[D, \delta] &= D^2 - \alpha D - \Rr_{a_0}^\star, &
D([u, v]) &= u \star D(v) - v \star D(u), \\
D(u) \star v &= u \star \delta(v) + \delta(u) \star v - \delta(u \star v),
\end{aligned}
\end{equation}
where $\delta$ is anti-symmetric with respect to $\prs$.

Define the $\mathbb{Z}_2$-graded vector space $\bar{\mathcal{A}} = \K e\oplus \mathcal{A} \oplus \K d$, equipped with the bilinear products $\bar{\star}$ and $\bar{\circ}$ defined for all $u, v \in \mathcal{A}$ by:
\begin{equation}\label{eq:extensionLie}
\begin{cases}
d \bar{\star} d = a_0 +\al d, &
u \bar{\star} v = u \star v - \langle D(u), v \rangle_\A e, \\
d \bar{\star} u = \delta(u) - \langle a_0, u \rangle_\A e, &
u \bar{\star} d = D(u),\\
d\bar{\star} e=-\al e, &  u\bar{\star} e=e\bar{\star} u=e\bar{\star} e= e\bar{\star} d= 0,
\end{cases}
\end{equation}
where $\al,  \lambda \in \mathbb{K}$.
Then, $(\bar{\mathcal{A}}, \bar{\star})$ is a left-symmetric superalgebra. Moreover, the even bilinear form $\prs : \bar{\mathcal{A}} \times \bar{\mathcal{A}} \to \mathbb{K}$ defined by
$$
\prs|_{\mathcal{A} \times \mathcal{A}} = \prs_\A, \quad \langle d, e\rangle = \langle e, d\rangle= 1,
$$
is a pseudo-Euclidean scalar product on $\bar{\mathcal{A}}$.

The flat even pseudo-Euclidean Lie superalgebra $(\bar{\mathcal{A}},\bar{\bullet}, \prs)$ is called the even double extension of $(\mathcal{A}, \br_{\bar{\A}}, \prs_\A)$ by means of $(\delta, D,  a_0, \al)$, where $[u\, v]_{\bar{\A}}=u\bar{\star} v-(-1)^{|u||v|} v\bar{\star} u$,  for all $u,v\in \bar{\A}$.

\end{Corollary}

\begin{proof}
When the left Leibniz superalgebra $(\mathcal{A}, \br)$ is anticommutative, that is, when $(\mathcal{A}, \br)$ is a Lie superalgebra, the Levi-Civita products associated with the pseudo-Euclidean Lie superalgebra satisfy $\star = -\circ$. In this case, according to the structure of the products \eqref{eq:extension}, we deduce that
$$
 \delta = -\delta^*, \quad G=-D, \quad a_0 = -b_0, \quad c_0 = 0, \quad \lambda = 0.
$$
By substituting these into the system \eqref{desapp}, we obtain  \eqref{mapsLie}.

\end{proof}
We now introduce the notion of an odd double extension of flat even pseudo-Euclidean  left Lebniz superalgebras by a purely odd one-dimensional vector superspace.
\sssbegin{Theorem}\label{Thdbodd}
 Let $(\mathcal{A}, \bullet, \prs_\A)$ be a flat even pseudo-Euclidean left Leibniz superalgebra, and let $(\star, \circ)$ be its associated Levi-Civita products. Let $V=V_{\bar{1}} = \mathbb{K} d$ be a purely odd one-dimensional $\mathbb{Z}_2$-graded  vector space  and $V^* = \mathbb{K} e$ its dual. Suppose $(\delta, D, \delta^*, G, a_0, b_0, 0,0)$ is an admissible of $\mathcal{A}$ satisfying the compatibility 
 conditions:

$$
\begin{array}{lclrcl}
  R_{c_0}^\circ + (R_{c_0}^\star)^* &=& G^*\circ D+D^*\circ G,&  \de\circ D+ \de^*\circ D & = & 0,\\   
  D^*(b_0-c_0) & = &-G^*(a_0) -\de(c_0),& 
  G^*(a_0+c_0)+D^*(b_0) - \de^*(c_0)& = & 0,\\  G(a_0-b_0) & = &-2\de^*(b_0),&
  ( \de+D^*)(c_0) - G^*(a_0+b_0)& =&0  ,\\  (\de^*+G^*)(c_0) & = & (G^*-D^*)(b_0),& 
  \langle a_0, a_0+b_0\rangle_\A & = & 0,\\ \langle a_0-b_0,c_0\rangle_\A & = & 2\langle a_0, b_0\rangle_\A,&
   G^*\circ D + D^*\circ D& = &0,\\ G^*\circ D + D^*\circ G & = & D^*\circ G-G^*\circ G,&
  (D-\de)(a_0-b_0)&=&0,\\   (\de+\de^*)(a_0)&=&0,&
   D^*(a_0+b_0) & = & 0, \\ \Rr_{a_0}^\star +\Rr_{a_0}^\circ& = & 0,&
    (D^*+G^*)(a_0)& = &0, \\ u\star G(v)&=&u\circ D(v)=-u\star D(v),& 
    \langle a_0-b_0,b_0-c_0\rangle_\A &= &0,\\
    u\circ D(v)&=&-u\circ G(v).
  \end{array}
$$
and \textup{(}for all $u\in \A$\textup{)}
$$
\begin{array}{lcl}
(\de^*\circ G-G\circ \de^*)(u) & = & (-1)^{|u|}(G\circ D+\Rr^\circ_{b_0})(u),\\ 
  (G\circ \de+\de^*\circ G)(u) & = & (-1)^{|u|}(\Rr^\circ_{b_0}-G^2)(u),\\
  (D\circ \de^*-\de\circ G)(u) & = & -(-1)^{|u|}(D^2+\Rr^\star_{a_0})(u),\\ (D\circ \de-\de\circ D)(u) & = & (-1)^{|u|}(\Rr^\star_{b_0}-D\circ G)(u).
\end{array}
$$
Define the $\mathbb{Z}_2$-graded vector space $\bar{\mathcal{A}} = \K e\oplus \mathcal{A} \oplus \K d$, equipped with the bilinear products $\bar{\star}$ and $\bar{\circ}$ defined for all $u, v \in \mathcal{A}$ by:
\begin{equation}\label{eq:extensionodd}
\begin{cases}
d \bar{\star} d = a_0, \\
u \bar{\star} v = u \star v + (-1)^{|v|}\langle G(u), v \rangle_\A e, \\
d \bar{\star} u = \delta(u) +(-1)^{|u|} \langle b_0, u \rangle_\A e, \\
u \bar{\star} d = D(u) + \langle c_0, u \rangle_\A e,\\
u\bar{\star} e=e\bar{\star} u=e\bar{\star} e= e\bar{\star} d=d \bar{\star} e=0,
\end{cases}
\quad \quad \text{and} \quad \quad
\begin{cases}
d \bar{\circ} d = b_0, \\
u \bar{\circ} v = u \circ v -(-1)^{|u|} \langle D(u), v \rangle_\A e, \\
d \bar{\circ} u = \delta^*(u) + \langle a_0, u \rangle_\A e, \\
u \bar{\circ} d = G(u) -(-1)^{|u|}  \langle c_0, u \rangle_\A e,\\
u\bar{\circ} e=e\bar{\circ} u=e\bar{\circ} e= e\bar{\circ} d=d \bar{\circ} e=0,
\end{cases}
\end{equation}
where $\al, \beta, \lambda \in \mathbb{K}$ and $c_0 \in \mathcal{A}_{\bar{0}}$.
Then, $(\bar{\mathcal{A}}, \bar{\star}, \bar{\circ})$ is a pre-left Leibniz superalgebra. Moreover, the even bilinear form $\prs : \bar{\mathcal{A}} \times \bar{\mathcal{A}} \to \mathbb{K}$ defined by
$$
\prs|_{\mathcal{A} \times \mathcal{A}} = \prs_\A, \quad \langle e, d\rangle =- \langle d, e\rangle= 1,
$$
is a pseudo-Euclidean scalar product on $\bar{\mathcal{A}}$.

The flat even pseudo-Euclidean left Leibniz superalgebra $(\bar{\mathcal{A}},\bar{\bullet}, \prs)$ is called the odd double extension of $(\mathcal{A}, \bullet, \prs_\A)$ by means of $(\delta, D, \delta^*, G, a_0, b_0, c_0)$, where $u\bar{\bullet} v=u\bar{\star} v+(-1)^{|u||v|} v\bar{\circ} u$, for all $u,v\in \bar{\A}$.   
\end{Theorem}

\begin{proof}
Since $(\delta, D, \delta^*, G, a_0, b_0, 0,0)$ is an admissible of $\mathcal{A}$, and $u\star G(v)=u\circ D(v)=-u\star D(v)=-u\circ G(v)$, for all $u,v\in\A$, then the pair $(G, D)$ satisfies conditions (i), (ii), and (iii) of Proposition~\ref{Auconidition1}, it follows from that proposition that the two even bilinear maps
$$
\mu, \gamma : \A \times \A \longrightarrow \mathbb{K}, \quad \mu(u, v) = (-1)^{|v|} \langle G(u), v \rangle_\A, \quad \gamma(u, v) = -(-1)^{|v||u|} \langle D(u), v \rangle_\A,
$$
satisfy the relations \eqref{E1}-\eqref{E3}. Consequently, the $\mathbb{Z}_2$-graded vector space $\tilde{\mathcal{A}} = \A \oplus \mathbb{K}e$, endowed with the products
$$
(u + a e)\, \bar{\star}\, (v + b e) = u \star v + (-1)^{|v|} \langle G(u), v \rangle_\A e, \quad (u + a e)\, \bar{\circ}\, (v + b e) = u \circ v - (-1)^{|u|} \langle D(u), v \rangle_\A e,
$$
for all $u, v \in \mathcal{A}$ and $a, b \in \mathbb{K}$, is a pre-left Leibniz superalgebra, which is a central extension of $\mathcal{A}$ by means of $(\mu, \gamma)$.

We now define even endomorphisms $\delta', D', \xi', G' : \tilde{\mathcal{A}} \to \tilde{\mathcal{A}}$ by: 
$$
\begin{array}{ll}
\delta'(u + a e) = \delta(u) +(-1)^{|u|} \langle b_0,  u \rangle_\A  e, &
D'(u + a e) =  D(u) +\langle c_0, u \rangle_\A e, \\[1mm]
\xi'(u + a e)   =  \delta^*(u) + \langle a_0,  u \rangle_\A e, &  
G'(u + a e)   =  G(u) -(-1)^{|u|} \langle c_0, u \rangle_\A e,
\end{array}
$$
and set
$$
a_0' = a_0 , \quad b_0' = b_0 ,
$$
for all $u \in \A$, $a \in \mathbb{K}$.

To show that $(\delta', D', \xi', G', a_0', b_0', 0, 0)$ is an admissible of $\tilde{\mathcal{A}}$, we verify one of the key admissibility conditions. Let $u, v \in \mathcal{A}$. We have,
$$
\begin{aligned}
D'&(u \bar{\star} v) - u \bar{\star} D'(v) + (-1)^{|u||v|} \left( D'(v \bar{\circ} u) + v \bar{\star} D'(u) \right) \\
&= D(u \star v) + \langle c_0, u \star v \rangle_\A e - u \star D(v) -(-1)^{|d|+|v|} \langle G(u), D(v) \rangle_\A e \\
&\quad + (-1)^{|u||v|} \left( D(v \circ u) + \langle c_0, v \circ u \rangle_\A e + v \star D(u) +(-1)^{|d|+|u|} \langle G(v), D(u) \rangle_\A e \right) \\
&= \left( D(u \star v) - u \star D(v) + (-1)^{|u||v|}( D(v \circ u) + v \star D(u) ) \right) \\
&\quad + \text{(} \langle  u \star v, c_0 \rangle_\A + (-1)^{|u||v|} \langle v\star  c_0, u  \rangle_\A   \\ &\quad-(-1)^{|d|+|v|} \langle G(u), D(v) \rangle_\A + (-1)^{|u||v|+|d|+|u|} \langle G(v), D(u) \rangle_\A \text{)} e.
\end{aligned}
$$
Since $(\delta, D, \delta^*, G, a_0, b_0, 0, 0)$ is an admissible of $\mathcal{A}$, the expression in the first parentheses vanishes. Therefore, we have

$$
D'(u \bar{\star} v) - u \bar{\star} D'(v) + (-1)^{|u||v|} \left( D'(v \bar{\circ} u) + v \bar{\star} D'(u) \right)
= \Theta(u, v)\, e,
$$
where
$$
\begin{aligned}
\Theta(u, v) &= \langle  u \star v ,c_0\rangle_\A + (-1)^{|u||v|} \langle v\star c_0,  u \rangle_\A -(-1)^{|d|+|v|} \langle G(u), D(v) \rangle_\A + (-1)^{|u||v|+|d|+|u|} \langle G(v), D(u) \rangle_\A.
\end{aligned}
$$
Now, applying duality properties, we observe
$$
\begin{aligned}
(-1)^{|u||v| + |d| + |u|} \langle G(v), D(u) \rangle_\mathcal{A}
&= (-1)^{|u||v| + |d| + |u| + |v|} \langle v, G^* D(u) \rangle_\mathcal{A} \\
&= (-1)^{|d| + |u| + |v|} \langle G^* D(u), v \rangle_\mathcal{A}.
\end{aligned}
$$
Since $d$ is odd, we deduce
$$
(-1)^{|u||v| + |d| + |u|} \langle G(v), D(u) \rangle_\mathcal{A} = -(-1)^{|u| + |v|} \langle G^* D(u), v \rangle_\mathcal{A}.
$$
In particular, if $\langle G^* D(u), v \rangle_\mathcal{A} \neq 0$, this expression vanishes only when $|u| + |v| = \bar{0}$, i.e., when $u$ and $v$ are both even or both odd. Similarly, the other terms can be rewritten using the adjoint maps.

Moreover, we can rewrite $\Theta(u, v)$ as
$$
\Theta(u, v) = \langle  R_{c_0}^\circ(u) + (R_{c_0}^\star)^*(u) - G^* D(u) - D^* G(u), v \rangle_\mathcal{A}.
$$
Since $\prs_\mathcal{A}$ is nondegenerate, we conclude that the identity
$$
R_{c_0}^\circ(u) + (R_{c_0}^\star)^*(u) = G^* D(u) + D^* G(u)
$$
holds if and only if $\Theta(u, v) = 0$ for all $u, v \in \mathcal{A}$. This proves the desired compatibility condition for $D'$. The remaining admissibility conditions for $(\delta', \xi', G')$ can be verified similarly using the same type of computation.
Hence, $(\delta', D', \xi', G', a_0', b_0',0,0)$ is an admissible of $\Tilde{\mathcal{A}}$. We may therefore consider the pre-left Leibniz superalgebra $\bar{\mathcal{A}} = V^* \oplus \mathcal{A} \oplus V$, which is the semi-direct product of the central extension $\tilde{\mathcal{A}} = \mathcal{A} \oplus V^*$ by the one-dimensional vector space $V = \mathbb{K}d$.

Furthermore, it is straightforward to verify that the bilinear form
$$
\prs : \bar{\mathcal{A}} \times \bar{\mathcal{A}} \to \mathbb{K}, \quad \text{defined by } \prs|_{\mathcal{A} \times \mathcal{A}} = \prs_\A, \quad \langle e, d \rangle =-\langle d, e \rangle= 1,
$$
is an even pseudo-Euclidean scalar product on $\bar{\mathcal{A}}$. Moreover, the compatibility condition
$$
\langle u \bar{\star} v, w \rangle = (-1)^{|u||v|} \langle v, u \bar{\circ} w \rangle, \quad \forall u, v, w \in \bar{\mathcal{A}},
$$
is satisfied.

Hence, $(\bar{\mathcal{A}}, \bar{\bullet}, \prs)$ is a flat even pseudo-Euclidean left Leibniz superalgebra, where $u\bar{\bullet} v=u \bar \star v+(-1)^{|u||v|} v \bar \circ u$, for all $u,v\in \bar{\A}$.

\end{proof}

We now introduce the notion of an odd double extension of flat even pseudo-Euclidean Lie superalgebras by a purely odd one-dimensional vector superspace.
\sssbegin{Corollary} \label{eq:extensionLie-evod}
Let $(\mathcal{A}, \bullet, \prs_\A)$ be a flat even pseudo-Euclidean Lie superalgebra, and let $\star$ be its associated Levi-Civita product. Let $V=V_{\bar{1}} = \mathbb{K} d$ be a purely odd one-dimensional $\mathbb{Z}_2$-graded  vector space  and $V^* = \mathbb{K} e$ its dual.  Suppose two endomorphisms $D$ and $\delta$ of $\mathcal{A}$, and an element $a_0 \in \mathcal{A}$, satisfy the following relations for all $u, v \in \mathcal{A}$, 
\begin{equation} 
\begin{cases}
D([u , v]) = u \star D(v)  - (-1)^{|u||v|}v \star D(u),\\
D(u) \star v = u \star \delta(v)  +(-1)^{|u|} \de(u) \star v - (-1)^{|u|}\delta(u \star v), \;\;  
\de^2-\Ll_{a_0}^\star =0,\\  (R_{c_0}^\star)^* -R_{c_0}^\star = D^*\circ D=\Rr^\star_{a_0}=0,\; \;(\de\circ D-D\circ \de)(u)=(-1)^{|u|} D^2(u),\\ (D-\de)(a_0)=D^*(c_0)=\de(c_0)=D^*(a_0)=0,\; \;  \langle a_0, a_0+c_0\rangle_\A=0.
  \end{cases}
\end{equation}
where $\delta$ is anti-symmetric with respect to $\prs$.

Define the $\mathbb{Z}_2$-graded vector space $\bar{\mathcal{A}} = \K e\oplus \mathcal{A} \oplus \K d$, equipped with the bilinear products $\bar{\star}$ and $\bar{\circ}$ defined for all $u, v \in \mathcal{A}$ by:
\begin{equation}\label{eq:extensionLie1}
\begin{cases}
d \bar{\star} d = a_0, &
u \bar{\star} v = u \star v -(-1)^{|v|} \langle D(u), v \rangle_\A e, \\
d \bar{\star} u = \delta(u) - (-1)^{|u|}\langle a_0, u \rangle_\A e, &
u \bar{\star} d = D(u)+\langle c_0, u\rangle_\A e,\\
u\bar{\star} e=e\bar{\star} u=e\bar{\star} e= e\bar{\star} d=d\bar{\star} e= 0,
\end{cases}
\end{equation}
where $\al,  \lambda \in \mathbb{K}$.
Then, $(\bar{\mathcal{A}}, \bar{\star})$ is a left-symmetric superalgebra. Moreover, the even bilinear form $\prs : \bar{\mathcal{A}} \times \bar{\mathcal{A}} \to \mathbb{K}$ defined by
$$
\prs|_{\mathcal{A} \times \mathcal{A}} = \prs_\A, \quad \langle d, e\rangle = \langle e, d\rangle= 1,
$$
is a pseudo-Euclidean scalar product on $\bar{\mathcal{A}}$.

The flat even pseudo-Euclidean Lie superalgebra $(\bar{\mathcal{A}},\br_{\bar{\A}}, \prs)$ is called the odd double extension of $(\mathcal{A}, \br, \prs_\A)$ by means of $(\delta, D,  a_0)$, where $[u\, v]_{\bar{\A}}=u\bar{\star} v-(-1)^{|u||v|} v\bar{\star} u$,  for all $u,v\in \bar{\A}$.

\end{Corollary}
\begin{proof}
The proof of this corollary is similar to the proof of Corollary \ref{dbexLie1}.
\end{proof}

We now aim to establish the converses of Theorem~\ref{Thdbeven} and Theorem~\ref{Thdbodd}.
\sssbegin{Theorem}\label{Th-pr}
Let $(\mathcal{A}, \bullet, \prs)$ be a flat even pseudo-Euclidean non-Lie left Leibniz superalgebra of superdimension $n|m$. Then $(\mathcal{A}, \bullet, \prs)$ is either:
\begin{itemize}
\item[(i)] an even double extension of a flat even pseudo-Euclidean left Leibniz superalgebra $(\mathcal{B}, \bullet_\mathcal{B}, \prs_\mathcal{B})$ of superdimension $n-2|m$ by means of $(\delta, G, \delta^*, G, a_0, b_0, c_0, 0,\lambda)$, or
\item[(ii)] an odd double extension of a flat even pseudo-Euclidean left Leibniz superalgebra $(\mathcal{B}, \bullet_\mathcal{B}, \prs_\mathcal{B})$ of superdimension $n|m-2$ by means of $(\delta, G, \delta^*, G, a_0, b_0, c_0)$.
\end{itemize}
\end{Theorem}

\begin{proof}
According to Proposition~\ref{ideal}, there exists a $\mathbb{Z}_2$-graded subspace $I := \mathbb{K}e$ which is a one-dimensional, totally isotropic two-sided ideal of $(\mathcal{A}, \star, \circ)$ and $e\star u=u\star e=0$, for all $u\in \A$. Moreover, its orthogonal complement $I^\perp$ is also a two-sided ideal of $(\mathcal{A}, \star, \circ)$. We distinguish two cases depending on the parity of the generator $e$ of $I$: either $I \cap \mathcal{A}_{\bar{0}} \neq \{0\}$ or $I \cap \mathcal{A}_{\bar{1}} \neq \{0\}$.

We begin with the first case and assume that $I \subset \mathcal{A}_{\bar{0}}$, i.e., $e \in \mathcal{A}_{\bar{0}}$. Since the bilinear form $\prs$ is even and non-degenerate, there exists an even element $d \in \mathcal{A}_{\bar{0}} \setminus \{0\}$ such that
$
\langle e, d \rangle = \langle d, e \rangle = 1.
$ Define the subspaces $V := \mathbb{K}d$ and $
H := (I \oplus V)^\perp.
$
Then, we have the orthogonal direct sum decomposition
$
\mathcal{A} = \mathbb{K}e \oplus H \oplus \mathbb{K}d,$
and $I^\perp = \mathbb{K}e \oplus H.
$ Let $u, v \in H$. Since $I^\perp$ is a two-sided ideal of $(\mathcal{A}, \star, \circ)$, the products in $\mathcal{A}$ restrict as follows
$$
u \star v = u \star_0 v + \mu(u, v)e, \quad
u \circ v = u \circ_0 v + \gamma(u, v)e,
$$
where $\star_0$ and $\circ_0$ denote the restrictions of $\star$ and $\circ$ to $H$, and $\mu, \gamma : H \times H \to \mathbb{K}$ are even bilinear forms. It is easy to verify that $(H, \star_0, \circ_0)$ is a pre-left Leibniz superalgebra and that the restriction $\prs_H := \prs|_{H \times H}$ is an even non-degenerate symmetric bilinear form on $H$, satisfying the identity
$
\langle u \star_0 v, w \rangle_H = (-1)^{|u||v|} \langle v, u \circ_0 w \rangle_H,$
\; for all $u, v, w \in H.
$
Hence, $(H, \bullet_0, \prs_H)$ is a flat even pseudo-Euclidean left Leibniz superalgebra, where $u\bullet_0v=u\star_0 v+(-1)^{|u||v|} v\circ_0 u$, for all $u,v\in \A$. Moreover, the even bilinear forms $\mu$ and $\gamma$ satisfy relations~\eqref{E1}-\eqref{E3}.

Now, since $I^\perp = \mathbb{K}e \oplus H$ is a two-sided ideal of $(\mathcal{A}, \star, \circ)$, it follows that 

\begin{equation}
\begin{cases}
d \star d = a_0 + \lambda e+\al d,  \\
d \star u = \delta(u) +f(u) e, \\
u \star d = D(u) + g(u) e,\\  d\star e=e\star d=u\star e=0,\\e\star u=e\star e=0,
\end{cases}
\quad \quad \text{and} \quad \quad
\begin{cases}
d \circ d = b_0 + \lambda_1 e+\beta d, \\
d \circ  u = \xi(u) + h(u)e, \\
u \circ  d = G(u) +k(u) e,\\
d\circ e=e\circ d=u\circ e=0,\\
e\circ u=e\circ e=0
,\end{cases}
\end{equation}
where $\al,\beta, \lambda, \lambda_1 \in \mathbb{K}$, $f, g, h, k$ are even linear forms on $H$, and $\delta, D, \xi, G : H \to H$ are even linear maps.

The identity $
\langle X \star Y, Z \rangle = (-1)^{|X||Y|} \langle Y, X \circ Z \rangle,$ for all  $X, Y, Z \in \mathcal{A},
$ implies the following relations:
$$
\lambda = \lambda_1,\, \al=\beta=0 \quad
\mu(u, v) = \langle G(u), v \rangle_H, \quad
\gamma(u, v) = \langle D(u), v \rangle_H, \quad
\langle \delta(u), v \rangle_H = \langle u, \xi(v) \rangle_H,
$$
$$
f(u) = \langle b_0, u \rangle_H, \quad
h(u) = \langle a_0, u \rangle_H, \quad
g(u) = k(u).
$$
Since $\prs_H$ is non-degenerate and $g$ is an even linear form on $H$, there exists an even element $c_0 \in H_{\bar{0}}$ such that
$
g(u) = \langle c_0, u \rangle_H,$ for all $ u \in H.
$

Since $(\mathcal{A}, \bullet, \prs)$ is a flat even  pseudo-Euclidean left Leibniz superalgebra, a straightforward computation shows that  $(\delta, D, \delta^*, G, a_0, b_0, 0,0)$ is an admissible of $H$, and  $(\delta, D, \delta^*, G,  a_0, b_0, c_0, 0,\lambda)$ satisfies the system of equations~\eqref{desapp}.

Therefore, we conclude that $(\mathcal{A}, \bullet, \prs)$ is an even double extension of a flat even pseudo-Euclidean left Leibniz superalgebra $(H, \bullet_0, \prs_H)$ by means of $(\delta, G, \delta^*, G, a_0, b_0, c_0,0, \lambda)$.

 The second case is established in the same manner as the first. This completes the proof of Theorem.
 \end{proof}

\sssbegin{Corollary} \label{A partir-Lie}
Let $(\mathcal{A}, \bullet, \prs)$ be a flat even pseudo-Euclidean non-Lie left Leibniz superalgebra of superdimension $n|m$. Then $(\mathcal{A}, \bullet, \prs)$  is obtained from a flat even pseudo-Euclidean  Lie superalgebra  by a finite number of even or odd double extensions.    
\end{Corollary}
\begin{proof}
According to Theorem \ref{Th-pr}, the pseudo-Euclidean left Leibniz superalgebra $(\mathcal{A}, \bullet, \prs)$ is a double extension of $(\mathcal{B}_1, \bullet_{\mathcal{B}_1}, \prs_{\mathcal{B}_1})$. Therefore, $\mathcal{B}_1$ is either a Lie superalgebra or a non-Lie superalgebra. If $\mathcal{B}_1$ is a non-Lie superalgebra, then it is itself a double extension of another flat pseudo-Euclidean left Leibniz superalgebra $(\mathcal{B}_2, \bullet_{\mathcal{B}_2}, \prs_{\mathcal{B}_2})$. Since the superdimension of $\mathcal{A}$ is finite, there exists $k \in \mathbb{N}^*$ such that $\mathcal{B}_k$ is a Lie superalgebra or abelian Lie superalgebra.    
\end{proof}

We now aim to establish the converses of Corollary \ref{eq:extensionLie} and Corollary \ref{eq:extensionLie-evod}.
\sssbegin{Lemma} \label{ideal-Lie}
Let $(\A, \br, \prs)$ be a flat pseudo-Euclidean Lie superalgebra and let $\star$ be its associated Levi-Civita product. Let $I$ be a  totally isotropic two-sided ideal of $(\A, \star)$. Then,
\begin{itemize}
\item [$(i)$] The product $\star$ in $I$ is null, $I \star I^\bot= 0$, and $I^\bot$ is a left ideal of $(\A, \star)$.
\item [$(ii)$] $I^\bot $ is a right ideal if and only if  $I^\bot \star I = 0.$
\end{itemize}
\end{Lemma}
\begin{proof}
\begin{itemize}
\item[$(i)$] Let $u \in I$, $v \in I^\bot$, and $ w\in \A$.  The fact that $I$ is totally isotropic two-sided ideal, it follows that
$$
\langle  u\star v, w\rangle= -(-1)^{|u||v|} \langle   v, u\star w\rangle = 0.
$$
Thus, $I \star I^\bot = 0$. In particular, $I \star I = 0$. Moreover, since
$$
\langle  w\star v, u\rangle= -(-1)^{|w||v|} \langle   v, w\star u\rangle = 0,
$$
we deduce that $\A \star I^\bot \subset I^\bot$.
\item[$(ii)$]  For $u\in I$, $v \in I^\bot$, and $w \in \A$, we obtain
$$
\langle v \star w, u\rangle = 0 \quad \Leftrightarrow \quad \langle w, v\star u\rangle = 0,
$$
that is, $I^\bot \star \A \subset I^\bot$ if and only if $I^\bot\star I = 0$.\qed
\end{itemize}
\noqed
\end{proof}

\sssbegin{Theorem}\label{conv-Lie}
Let $(\mathcal{A}, \br, \prs)$ be a flat even pseudo-Euclidean Lie superalgebra, and let $\star$ be its associated Levi-Civita product. Then
\begin{itemize}
\item[$(i)$] If $I$ is a one-dimensional totally isotropic two-sided ideal of $(\mathcal{A}, \star)$, with $I \subseteq \mathcal{A}_0$, and if $I^\perp$ is a totally isotropic two-sided ideal of $(\mathcal{A}, \star)$, then $(\mathcal{A}, \br, \prs)$ is an even double extension of a flat even pseudo-Euclidean Lie superalgebra $(\mathcal{B}, \br_\mathcal{B}, \prs_\mathcal{B})$ of superdimension 
$n-2|m$ by means of the data $(\delta, D, a_0, \alpha)$.
\item[$(ii)$] If $I$ is a one-dimensional totally isotropic two-sided ideal of $(\mathcal{A}, \star)$, with $I \subseteq \mathcal{A}_1$, and if $I^\perp$ is a totally isotropic two-sided ideal of $(\mathcal{A}, \star)$, then $(\mathcal{A}, \br, \prs)$ is an odd double extension of a flat even pseudo-Euclidean Lie superalgebra $(\mathcal{B}, \br_\mathcal{B}, \prs_\mathcal{B})$ of superdimension $n|m-2$ by means of the data $(\delta, D, a_0, \alpha)$.
\end{itemize}
\end{Theorem}
\begin{proof}
Assume that $I$  is a one-dimensional totally isotropic two-sided ideal of $(\mathcal{A}, \star)$, with $I \subseteq \mathcal{A}_0$ and  $I^\perp$ is a totally isotropic two-sided ideal of $(\mathcal{A}, \star)$. Put $I:=\K e$ where $e\in I$.  Since the bilinear form $\prs$ is even and non-degenerate, there exists an even element $d \in \mathcal{A}_{\bar{0}} \setminus \{0\}$ such that
$
\langle e, d \rangle = \langle d, e \rangle = 1.
$ Define the subspaces $V := \mathbb{K}d$ and $
H := (I \oplus V)^\perp.
$
Then, we have the orthogonal direct sum decomposition
$
\mathcal{A} = \mathbb{K}e \oplus H \oplus \mathbb{K}d,$
and $I^\perp = \mathbb{K}e \oplus H.
$ Let $u, v \in H$. Since $I^\perp$ is a two-sided ideal of $(\mathcal{A}, \star)$, the product in $\mathcal{A}$ restrict as follows
$$
u \star v = u \star_0 v + \mu(u, v)e, 
$$
where $\star_0$  denote the restriction of $\star$  to $H$, and $\mu,  : H \times H \to \mathbb{K}$ are even bilinear form. It is easy to verify that $(H, \star_0)$ is a left-symmetric superalgebra and that the restriction $\prs_H := \prs|_{H \times H}$ is an even non-degenerate symmetric bilinear form on $H$, satisfying the identity
$
\langle u \star_0 v, w \rangle_H = -(-1)^{|u||v|} \langle v, u \star_0 w \rangle_H,$
\; for all $u, v, w \in H.
$
Hence, $(H, \br_0, \prs_H)$ is a flat even pseudo-Euclidean Lie superalgebra, where $[u, v]_0=u\star_0 v-(-1)^{|u||v|} v\star_0 u$, for all $u,v\in \A$. Moreover, the even bilinear forms $\mu$ and $\gamma$ satisfy relation~\eqref{cocycle}. Since $I$ and $I^\bot$ are two-sided of $(\A, \star)$, according  to Lemma \ref{ideal-Lie}, $I\star I^\bot=I^\bot\star I=\{0\}$ and we have 
$$
\begin{cases}
d \star d = a_0 + \al d+ \lambda e,  & d \star u = \delta(u) +f(u) e, \\
u \star d = D(u) + g(u) e,& 
d\star e=\beta e,\\
e\star d=\ga e,& e\star e=e\star u=u\star e=0,
\end{cases} $$
where $\al, \la, \beta, \ga\in \K$, $f, g$ are even linear forms on $H$, and $\delta, D : H \to H$ are even linear maps.

The identity $
\langle X \star Y, Z \rangle = -(-1)^{|X||Y|} \langle Y, X \circ Z \rangle,$ for all  $X, Y, Z \in \mathcal{A},
$ implies the following relations:
$$
\al = -\beta,\; g=0,\;\la=\ga=0, \;
\mu(u, v) = -\langle D(u), v \rangle_H, \; f(u) =- \langle a_0, u \rangle_H,\; \esp
$$
$$
\langle \de(u), v\rangle_H =-\langle u, \de(v)\rangle_H.
$$
Thus, we get that 
\begin{equation}
\begin{cases}
d \star d = a_0 + \al d, & 
u\star v=u\star_0 v-\langle D(u), v\rangle_H e,\\
d \star u = \delta(u) -\langle a_0, u\rangle_H e, & 
u \star d = D(u) ,\\
e\star d=-\al e,&  d\star e= e\star e=e\star u=u\star e=0,
\end{cases} \label{prleft}
\end{equation}
Finally, the product \eqref{prleft} is left symmetric if and only if for all $u,v\in H$ we have that
\begin{itemize}
\item $(u\star v)\star d-u\star (v\star d)=(-1)^{|u|v|}(v\star u)\star d-v\star (u\star d)$. This equality holds if and only if
$$D([u,v]_0)=u\star_0 D(v)-(-1)^{|u||v|} v\star_0 D(u).$$
\item $(u\star d)\star v-u\star (d\star v)=(d\star u)\star v-d\star (u\star v)$. This equality holds if and only if $$[D, \de]=D^2-\al D-\Rr_{a_0}^\star,\esp D(u)\star_0 v=u\star_0 \de(v)+\de(u)\star_0 v-\de(u\star_0 v).$$
\item $(u\star d)\star d-u\star (d\star d)=(d\star u)\star d-d\star (u\star d)$. This equality holds if and only if
$$[D, \de]=D^2-\al D-\Rr_{a_0}^\star.$$
\end{itemize}
Thus, we get that $(\mathcal{A}, \br, \prs)$ is an odd double extension of a flat even pseudo-Euclidean Lie superalgebra $(H, \br_H, \prs_H)$ of superdimension $n|m-2$ by means of the data $(\delta, D, a_0, \alpha)$.

The second case is established in the same way as the first. This completes the proof
of Theorem.
\end{proof}
\subsection{Double extension of flat odd pseudo-Euclidean left-Leibniz superalgebras}

We start by presenting the concept of even double extension of flat odd pseudo-Euclidean superalgebras by the one-dimensional vector space.

\sssbegin{Theorem}\label{Thdbeven2}
Let $(\mathcal{A}, \bullet, \prs_\A)$ be a flat odd pseudo-Euclidean left Leibniz superalgebra, and let $(\star, \circ)$ be its associated Levi-Civita products. Let $V = \mathbb{K} d$ be a one-dimensional  vector space and $\Pi(V^*) = (\mathbb{K} e)_{\bar{1}}$, where $V^*$ its dual. Suppose $(\delta, D, \delta^*, G, a_0, b_0, \al, -\al)$ is an admissible of $\mathcal{A}$ satisfying the compatibility conditions: 
\begin{equation}\label{desapp2}
\begin{array}{lclrcl}
  \Rr_{c_0}^\circ + (R_{c_0}^\star)^* &=&  D^* \circ G- G^*\circ D,& [G, \de^*] &=&-G\circ D+\al G-\Rr_{b_0}^\circ,\; \\G\circ \de+\de^*\circ G &=& -G^2-\al G+\Rr_{b_0}^\circ,& [D, \de] & = &-D\circ G-\al D+\Rr_{b_0}^\star,\\ D\circ \de^*-\de\circ G & = & -D^2+\al D+\Rr_{a_0}^\star,& G^*\circ D & = &-D^*\circ D,\\ \de\circ D&=&-\de^*\circ D, &
\Rr_{a_0}^\star+\Rr_{a_0}^\circ & = & -\alpha (D+G), \\ \Rr_{c_0}^\circ + (\Rr_{c_0}^\star)^* & = &  D^* \circ G+ G^* \circ G, &
G(a_0+b_0)& = &0,\\  D^*(a_0+b_0)& =&0 , &(G^*+D^*)(a_0)& = &-2\al c_0,\\  (G+\de)(a_0+b_0)& =&0 , & (\de+\de^*)(a_0)& = &-\al(a_0+b_0) , \\
\langle a_0+b_0, c_0\rangle_\A& =  &0,& \langle a_0+b_0, a_0\rangle_\A & =  &0,\\ 
  \langle a_0+b_0, b_0\rangle_\A& = &0,& 
    (\de+D^*)(c_0)& =& G^*(a_0+b_0) , \\ (\de^*+G^*)(c_0)&=&-G^*(a_0)+D^*(b_0) ,&  (\de^*+G^*)(c_0)& =& (G^*+D)(b_0)-2\al c_0,\\ D^*(b_0+c_0)&=&2\al c_0-\de(c_0)+G^*(a_0), & D^* \circ G- G^*\circ D&=& G^* \circ G+D^*\circ G, \\ u\star G(v)&=&u\circ D(v)=-u\star D(v),\\
    u\star D(v)&=&u\circ G(v)
    \end{array}
\end{equation}
Define the $\mathbb{Z}_2$-graded vector space $\bar{\mathcal{A}} = \Pi(V^*)\oplus \mathcal{A} \oplus V$, equipped with the bilinear products $\bar{\star}$ and $\bar{\circ}$ defined for all $u, v \in \mathcal{A}$ by:
\begin{equation}\label{eq:extension-odd}
\begin{cases}
d \bar{\star} d = a_0 +\al d, \\
u \bar{\star} v = u \star v + \langle G(u), v \rangle_\A e, \\
d \bar{\star} u = \delta(u) + \langle b_0, u \rangle_\A e, \\
u \bar{\star} d = D(u) + \langle c_0, u \rangle_\A e,\\
d\bar{\star} e=-\al e, \; e\bar{\star} d= 0,
\end{cases}
\quad \quad \text{and} \quad \quad
\begin{cases}
d \bar{\circ} d = b_0 -\al d, \\
u \bar{\circ} v = u \circ v + \langle D(u), v \rangle_\A e, \\
d \bar{\circ} u = \delta^*(u) + \langle a_0, u \rangle_\A e, \\
u \bar{\circ} d = G(u) + \langle c_0, u \rangle_\A e,\\
d\bar{\circ} e=\al e, \; e\bar{\circ} d=0,
\end{cases}
\end{equation}
where $\al, \lambda \in \mathbb{K}$ and $c_0 \in \mathcal{A}_{\bar{0}}$.
Then, $(\bar{\mathcal{A}}, \bar{\star}, \bar{\circ})$ is a pre-left Leibniz superalgebra. Moreover, the odd bilinear form $\prs : \bar{\mathcal{A}} \times \bar{\mathcal{A}} \to \mathbb{K}$ defined by
$$
\prs|_{\mathcal{A} \times \mathcal{A}} = \prs_\A, \quad \langle d, e\rangle = \langle e, d\rangle= 1,
$$
is a pseudo-Euclidean scalar product on $\bar{\mathcal{A}}$.

The flat odd pseudo-Euclidean left Leibniz superalgebra $(\bar{\mathcal{A}},\bar{\bullet}, \prs)$ is called the even double extension of $(\mathcal{A}, \bullet, \prs_\A)$ by means of $(\delta, D, \delta^*, G, a_0, b_0, c_0, \al)$, where $u\bar{\bullet} v=u\bar{\star} v+(-1)^{|u||v|} v\bar{\circ} u$,  for all $u,v\in \bar{\A}$.
\end{Theorem}
\begin{proof}
 The proof of this Theorem is similar to that of Theorem \ref{Thdbeven}.   
\end{proof}

Since flat odd pseudo-Euclidean Lie superalgebras constitute a particular case of flat odd pseudo-Euclidean left Leibniz superalgebras, we now present the procedure of double extension for flat odd pseudo-Euclidean Lie superalgebras. 
\sssbegin{Corollary} \label{dbexLie3}
Let $(\mathcal{A}, \br, \prs_\A)$ be a flat odd pseudo-Euclidean Lie superalgebra, and let $\star$ be its associated Levi-Civita product. Let $V = \mathbb{K} d$ be a one-dimensional  vector space and $\Pi(V^*) = (\mathbb{K} e)_{\bar{1}}$, where $V^*$ its dual. Suppose two endomorphisms $D$ and $\delta$ of $\mathcal{A}$, and an element $a_0 \in \mathcal{A}$, satisfy the following relations for all $u, v \in \mathcal{A}$,
\begin{equation} \label{mapsLie-odd}
\begin{aligned}
[D, \delta] &= D^2 - \alpha D - \Rr_{a_0}^\star, &
D([u, v]) &= u \star D(v) - v \star D(u), \\
D(u) \star v &= u \star \delta(v) + \delta(u) \star v - \delta(u \star v), \\ 
\end{aligned}
\end{equation}
where $\delta$ is anti-symmetric with respect to $\prs$.

Define the $\mathbb{Z}_2$-graded vector space $\bar{\mathcal{A}} = \K e\oplus \mathcal{A} \oplus \K d$, equipped with the bilinear products $\bar{\star}$ and $\bar{\circ}$ defined for all $u, v \in \mathcal{A}$ by:
\begin{equation}\label{eq:extensionLie-odd}
\begin{cases}
d \bar{\star} d = a_0 +\al d, &
u \bar{\star} v = u \star v - \langle D(u), v \rangle_\A e, \\
d \bar{\star} u = \delta(u) - \langle a_0, u \rangle_\A e, &
u \bar{\star} d = D(u),\\
d\bar{\star} e=-\al e, &  u\bar{\star} e=e\bar{\star} u=e\bar{\star} e= e\bar{\star} d= 0,
\end{cases}
\end{equation}
where $\al,  \lambda \in \mathbb{K}$.
Then, $(\bar{\mathcal{A}}, \bar{\star})$ is a left-symmetric superalgebra. Moreover, the odd bilinear form $\prs : \bar{\mathcal{A}} \times \bar{\mathcal{A}} \to \mathbb{K}$ defined by
$$
\prs|_{\mathcal{A} \times \mathcal{A}} = \prs_\A, \quad \langle d, e\rangle = \langle e, d\rangle= 1,
$$
is a pseudo-Euclidean scalar product on $\bar{\mathcal{A}}$.

The flat odd pseudo-Euclidean Lie superalgebra $(\bar{\mathcal{A}},\br_{\bar{\A}}, \prs)$ is called the even double extension of $(\mathcal{A}, \br, \prs_\A)$ by means of $(\delta, D,  a_0, \al)$, where $[u,  v]_{\bar{\A}}=u\bar{\star} v-(-1)^{|u||v|} v\bar{\star} u$,  for all $u,v\in \bar{\A}$.

\end{Corollary}
\begin{proof}
The proof of this Corollary is similar to that of Corollary \ref{dbexLie1}.
\end{proof}

We now introduce the notion of an odd double extension of flat even pseudo-Euclidean  left Lebniz superalgebras by a purely odd one-dimensional vector space.

\sssbegin{Theorem}\label{Thdbodd2}
 Let $(\mathcal{A}, \bullet, \prs_\A)$ be a flat odd pseudo-Euclidean left Leibniz superalgebra, and let $(\star, \circ)$ be its associated Levi-Civita products. Let $V=V_{\bar{1}} = (\mathbb{K} d)_{\bar{1}}$ be a purely odd one-dimensional $\mathbb{Z}_2$-graded  vector space  and $\Pi(V^*) = \mathbb{K} e$ where $V^*$ its dual. Suppose $(\delta, D, \delta^*, G, a_0, b_0, 0, 0)$ is an admissible of $\mathcal{A}$ satisfying the compatibility conditions:
$$
\begin{array}{lclrcl}
  R_{c_0}^\circ + (R_{c_0}^\star)^* &=& -G^*\circ D-D^*\circ G,&  \de\circ D & = & -\de^*\circ D,\\  
  D^*(b_0-c_0) & = &-G^*(a_0) +\de(c_0),& 
  G^*(a_0+c_0)+D^*(b_0) & = &- \de^*(c_0),\\  G(a_0+b_0) & = &-2\de^*(b_0),&
  ( \de-D^*)(c_0) & =&G^*(a_0+b_0),\\  (\de^*-G^*)(c_0) & = & (G^*-D^*)(b_0), & 
  \langle a_0, a_0+b_0\rangle_\A & = & 0,\\ \langle a_0-b_0,c_0\rangle_\A & = & 2\langle a_0, b_0\rangle_\A, &
   G^*\circ D & = &-D^*\circ D,\\ G^*\circ D + D^*\circ G & = & -D^*\circ G-G^*\circ G,&
  (D-\de)(a_0-b_0)&=&0,\\   (\de-\de^*)(a_0)&=&0,&
   D^*(a_0+b_0) & = & 0,\\  \Rr_{a_0}^\star & = & - \Rr_{a_0}^\circ, &
    (D^*+G^*)(a_0)& = &0,\\ u\star G(v)&=&u\circ D(v)=-u\star D(v),& 
    \langle a_0,b_0-c_0\rangle_\A &= &\langle b_0, b_0+c_0\rangle_\A,\\
    u\star D(v)&=&u\circ G(v)
  \end{array}
$$
and (for all $u\in \A$)
$$
\begin{array}{lcl}
(\de^*\circ G-G\circ \de^*)(u) & = & (-1)^{|u|}(G\circ D+\Rr^\circ_{b_0})(u),\\ 
  (G\circ \de+\de^*\circ G)(u) & = & (-1)^{|u|}(\Rr^\circ_{b_0}-G^2)(u),\\
  (D\circ \de^*-\de\circ G)(u) & = & -(-1)^{|u|}(D^2-\Rr^\star_{a_0})(u),\\ (D\circ \de-\de\circ D)(u) & = & (-1)^{|u|}(\Rr^\star_{b_0}-D\circ G)(u).
\end{array}
$$
Define the $\mathbb{Z}_2$-graded vector space $\bar{\mathcal{A}} = \K e\oplus \mathcal{A} \oplus \K d$, equipped with the bilinear products $\bar{\star}$ and $\bar{\circ}$ defined for all $u, v \in \mathcal{A}$ by:
\begin{equation}\label{eq:extensionodd}
\begin{cases}
d \bar{\star} d = a_0 + \lambda e, \\
u \bar{\star} v = u \star v + (-1)^{|v|}\langle G(u), v \rangle_\A e, \\
d \bar{\star} u = \delta(u) +(-1)^{|u|} \langle b_0, u \rangle_\A e, \\
u \bar{\star} d = D(u) + \langle c_0, u \rangle_\A e,\\
u\bar{\star} e=e\bar{\star} u=d\bar{\star} e= e\bar{\star} d=0,
\end{cases}
\quad \quad \text{and} \quad \quad
\begin{cases}
d \bar{\circ} d = b_0 - \lambda e, \\
u \bar{\circ} v = u \circ v +(-1)^{|u|} \langle D(u), v \rangle_\A e, \\
d \bar{\circ} u = \delta^*(u) - \langle a_0, u \rangle_\A e, \\
u \bar{\circ} d = G(u) +(-1)^{|u|}  \langle c_0, u \rangle_\A e,\\
u\bar{\circ} e=e\bar{\circ} u=d\bar{\circ} e= e\bar{\circ} d=0,
\end{cases}
\end{equation}
where $ \lambda \in \mathbb{K}$ and $c_0 \in \mathcal{A}_{\bar{0}}$.
Then, $(\bar{\mathcal{A}}, \bar{\star}, \bar{\circ})$ is a pre-left Leibniz superalgebra. Moreover, the odd bilinear form $\prs : \bar{\mathcal{A}} \times \bar{\mathcal{A}} \to \mathbb{K}$ defined by
$$
\prs|_{\mathcal{A} \times \mathcal{A}} = \prs_\A, \quad \langle e, d\rangle = \langle d, e\rangle= 1,
$$
is a pseudo-Euclidean scalar product on $\bar{\mathcal{A}}$.

The flat odd pseudo-Euclidean left Leibniz superalgebra $(\bar{\mathcal{A}},\bar{\bullet}, \prs)$ is called the odd double extension of $(\mathcal{A}, \bullet, \prs_\A)$ by means of $(\delta, D, \delta^*, G, a_0, b_0, c_0, \la)$, where $u\bar{\bullet} v=u\bar{\star} v+(-1)^{|u||v|} v\bar{\circ} u$, for all $u,v\in \bar{\A}$.   
\end{Theorem}

\begin{proof}
The proof of this Theorem is similar to that of Theorem \ref{Thdbodd}.

\end{proof}
\sssbegin{Corollary} \label{A partir-Lie2}
Let $(\mathcal{A}, \bullet, \prs)$ be a flat odd pseudo-Euclidean non-Lie left Leibniz superalgebra of superdimension $n|m$. Then $(\mathcal{A}, \bullet, \prs)$  is obtained from a flat odd pseudo-Euclidean  Lie superalgebra  by a finite number of even or odd double extensions.    
\end{Corollary}

\begin{proof}
The proof of this Corollary is similar to that of Corollary \ref{A partir-Lie}. \end{proof}
\sssbegin{Corollary} \label{dbexLie4}
 Let $(\mathcal{A}, \bullet, \prs_\A)$ be a flat odd pseudo-Euclidean Lie superalgebra, and let $\star$ be its associated Levi-Civita product. Let $V=V_{\bar{1}} = (\mathbb{K} d)_{\bar{1}}$ be a purely odd one-dimensional $\mathbb{Z}_2$-graded  vector space  and $\Pi(V^*) = \mathbb{K} e$ where $V^*$ its dual. Suppose $(\delta, D,  a_0, 0)$ is an admissible of $\mathcal{A}$ satisfying the compatibility conditions:

  \begin{equation} 
\begin{cases}
D([u , v]) = u \star D(v)  - (-1)^{|u||v|}v \star D(u),\\
D(u) \star v = u \star \delta(v)  +(-1)^{|u|} \de(u) \star v - (-1)^{|u|}\delta(u \star v), \;\; \de^2-\Ll_{a_0}^\star =0,\\  (R_{c_0}^\star)^* -R_{c_0}^\star = D^*\circ D=0,\;\; (D\circ \de-\de\circ D)(u)=(-1)^{|u|} (D^2(u)-\Rr^\star_{a_0}),\\ D(a_0)=\de(a_0)=2D^*(a_0)+(\de+D^*)(c_0)=(\de-D^*)(c_0)=0,\\   
\langle a_0, a_0\rangle_\A=\langle a_0, c_0\rangle_\A=0.
  \end{cases}
\end{equation}
where $\delta$ is anti-symmetric with respect to $\prs$.

Define the $\mathbb{Z}_2$-graded vector space $\bar{\mathcal{A}} = \K e\oplus \mathcal{A} \oplus \K d$, equipped with the bilinear products $\bar{\star}$ and $\bar{\circ}$ defined for all $u, v \in \mathcal{A}$ by:
\begin{equation}\label{eq:extensionodd}
\begin{cases}
d \bar{\star} d = a_0 + \lambda e &
u \bar{\star} v = u \star v - (-1)^{|v|}\langle D(u), v \rangle_\A e, \\
d \bar{\star} u = \delta(u) -(-1)^{|u|} \langle a_0, u \rangle_\A e, &
u \bar{\star} d = D(u) + \langle c_0, u \rangle_\A e,\\
u \bar{\star} e=e\bar{\star} e=e\bar{\star} e=d\bar{\star} e=e\bar{\star} d=0,
\end{cases}
\end{equation}
where $ \lambda \in \mathbb{K}$ and $c_0 \in \mathcal{A}_{\bar{0}}$.
Then, $(\bar{\mathcal{A}}, \bar{\star}, \bar{\circ})$ is a left-symmetric superalgebra. Moreover, the odd bilinear form $\prs : \bar{\mathcal{A}} \times \bar{\mathcal{A}} \to \mathbb{K}$ defined by
$$
\prs|_{\mathcal{A} \times \mathcal{A}} = \prs_\A, \quad \langle e, d\rangle = \langle d, e\rangle= 1,
$$
is a pseudo-Euclidean scalar product on $\bar{\mathcal{A}}$.

The flat odd pseudo-Euclidean Lie superalgebra $(\bar{\mathcal{A}},\br_{\bar{\A}}, \prs)$ is called the odd double extension of $(\mathcal{A}, \bullet, \prs_\A)$ by means of $(\delta, D,  a_0, c_0, \la)$, where $[u,  v]_{\bar{\mathcal{\A}}}=u\bar{\star} v-(-1)^{|u||v|} v\bar{\star} u$, for all $u,v\in \bar{\A}$. 
 \end{Corollary}

\begin{proof}
The proof of this Corollary is similar to that of Corollary \ref{dbexLie1}.
\end{proof}

We now aim to establish the converses of Theorem~\ref{Thdbeven2} and Theorem~\ref{Thdbodd2}.
\sssbegin{Theorem}
Let $(\mathcal{A}, \bullet, \prs)$ be a flat odd pseudo-Euclidean non-Lie left Leibniz superalgebra of superdimension $n|n$. Then $(\mathcal{A}, \bullet, \prs)$ is either:
\begin{itemize}
\item[(i)] an even double extension of a flat odd pseudo-Euclidean left Leibniz superalgebra $(\mathcal{B}, \bullet_\mathcal{B}, \prs_\mathcal{B})$ of superdimension $n-1|n-1$ by means of $(\delta, G, \delta^*, G, a_0, b_0, c_0, \al)$, or
\item[(ii)] an odd double extension of a flat odd pseudo-Euclidean left Leibniz superalgebra $(\mathcal{B}, \bullet_\mathcal{B}, \prs_\mathcal{B})$ of superdimension $n-1|n-1$ by means of $(\delta, G, \delta^*, G, a_0, b_0, c_0, \lambda)$.
\end{itemize}
\end{Theorem}

\begin{proof}
The proof of this Theorem is similar to that of Theorem \ref{Th-pr}.
\end{proof}

We now aim to establish the converses of Corollary \ref{dbexLie3} and Corollary \ref{dbexLie4}.
\sssbegin{Theorem}
Let $(\mathcal{A}, \br, \prs)$ be a flat odd pseudo-Euclidean Lie superalgebra, and let $\star$ be its associated Levi-Civita product. Then
\begin{itemize}
\item[$(i)$] If $I$ is a one-dimensional totally isotropic two-sided ideal of $(\mathcal{A}, \star)$, with $I \subseteq \mathcal{A}_1$, and if $I^\perp$ is a totally isotropic two-sided ideal of $(\mathcal{A}, \star)$, then $(\mathcal{A}, \br, \prs)$ is an even double extension of a flat odd pseudo-Euclidean Lie superalgebra $(\mathcal{B}, \br_\mathcal{B}, \prs_\mathcal{B})$ of superdimension 
$n-1|m-1$ by means of the data $(\delta, D, a_0, \alpha)$.
\item[$(ii)$] If $I$ is a one-dimensional totally isotropic two-sided ideal of $(\mathcal{A}, \star)$, with $I \subseteq \mathcal{A}_0$, and if $I^\perp$ is a totally isotropic two-sided ideal of $(\mathcal{A}, \star)$, then $(\mathcal{A}, \br, \prs)$ is an odd double extension of a flat odd pseudo-Euclidean Lie superalgebra $(\mathcal{B}, \br_\mathcal{B}, \prs_\mathcal{B})$ of superdimension $n-1|m-1$ by means of the data $(\delta, D, a_0,c_0, \la)$.
\end{itemize}
\end{Theorem}
\begin{proof}
The proof of this Theorem is similar to that of Theorem \ref{conv-Lie}.
\end{proof}
\section{Conclusion and outlook}
This article studies flat pseudo-Euclidean Leibniz superalgebras and introduces the concept of \emph{double extension} for these algebras.  We prove that every flat pseudo-Euclidean {\bf non-Lie} Leibniz superalgebra can be constructed via successive double extensions starting from a flat pseudo-Euclidean {\bf Lie} superalgebra, see Corollary 
\ref{A partir-Lie} and \ref{A partir-Lie2}.  As a consequence, the characterization of flat pseudo-Euclidean Leibniz superalgebras is reduced to the study of flat pseudo-Euclidean Lie superalgebras.

The study of flat pseudo-Euclidean Lie algebras is well-established in several specific contexts. The foundational work began with Milnor's characterization of the flat Euclidean case \cite{Milnor1}. This was later extended to Lorentzian 2-step nilpotent algebras by Guediri \cite{guediri}. A significant generalization was achieved by Aubert and Medina \cite{Medina2}, who introduced the powerful concept of the double extension to provide a full characterization of flat Lorentzian nilpotent Lie algebras.
A complete classification of such algebras was given in~\cite{Bajo}, and further results can be found in~\cite{ABL, BL1}.  
In particular, all flat pseudo-Euclidean nilpotent Lie algebras of signature $(2,n-2)$ can be described via double extensions, see \cite{BL2}.
Motivated by these results, we plan to investigate other classes of flat pseudo-Euclidean Lie algebras and Lie superalgebras in future work.
For instance, we intend to study flat even pseudo-Euclidean Lie superalgebras whose even part $(\mathfrak{g}_{\bar 0}, [\cdot,\cdot]_{\bar 0}, \langle \cdot,\cdot \rangle_{\bar 0})$ is a flat Euclidean Lie algebra.
Such structures are expected to yield new families of flat pseudo-Euclidean Lie superalgebras and Leibniz superalgebras.  


\end{document}